\documentclass[12pt,a4paper]{article}
\usepackage{amssymb,amsmath,amsthm,euscript, color}
\usepackage{hyperref}

\newcommand*{\hm}[1]{#1\nobreak\discretionary{}%
            {\hbox{$\mathsurround=0pt #1$}}{}}

\newcommand{\End}{\mathop{\mathrm{End}}\nolimits}
\newcommand{\Aut}{\mathop{\mathrm{Aut}}\nolimits}
\newcommand{\Mat}{\mathop{\mathrm{Mat}}\nolimits}
\newcommand{\tr}{\mathop{\mathrm{tr}}\nolimits}

\newcommand{\Char}{\mathop{\mathrm{char}}\nolimits}
\newcommand{\qstar}{\mathop{\mathrm{*}}\nolimits_q}

\numberwithin{equation}{section}

\renewcommand{\det}{\mathop{\mathrm{det}}\nolimits}
\newcommand{\qdet}{\mathop{\mathrm{qdet}}\nolimits}
\renewcommand{\le}{\leqslant}
\renewcommand{\ge}{\geqslant}
\newcommand{\wt}{\widetilde}
\newcommand{\wh}{\widehat}
\newcommand{\inv}{\mathop{\mathrm{inv}}\nolimits}
\newcommand{\sg}{\mathfrak S} 

\newcommand{\pref}[1]{\ref{#1} page \pageref{#1} }

\def\qMMs{$q$-Manin\ matrices}
\def\qMM{$q$-Manin\ matrix}
\def\qDet{$q$-determinant{~}}

\def\qMMs{$q$-Manin\ matrices}
\def\qMM{$q$-Manin\ matrix}

\def\boxit#1#2{\setbox1=\hbox{\kern#1{#2}\kern#1}%
\dimen1=\ht1 \advance\dimen1 by #1 \dimen2=\dp1 \advance\dimen2 by #1
\setbox1=\hbox{\vrule height\dimen1 depth\dimen2\box1\vrule}%
\setbox1=\vbox{\hrule\box1\hrule}%
\advance\dimen1 by .4pt \ht1=\dimen1
\advance\dimen2 by .4pt \dp1=\dimen2 \box1\relax}

\def\bo#1{\boxit{1pt}{$#1$}}

\newcommand{\BEQ}{\begin{equation}}
\newcommand{\EEQ}{\end{equation}}
\newtheorem{Th}{Theorem}[section]
\newtheorem{Lem}[Th]{Lemma}
\newtheorem{Prop}[Th]{Proposition}
\newtheorem{Cor}{Corollary}[Th]

\newtheorem{Ex}{Example}[section]
\newtheorem{Def}{Definition}
\newtheorem{Rem}{Remark}[section]
\def\nn{\nonumber}
\def\one#1{#1^{\raise5pt\hbox{$\scriptstyle\!\!\!\!1$}}\,{}}
\def\two#1{#1^{\raise5pt\hbox{$\scriptstyle\!\!\!\!2$}}\,{}}

\def\bea{\begin{eqnarray}}
\def\eea{\end{eqnarray}}

\newcommand{\CC}{{{\mathbb{C}}}}

\def\T{\mathcal{T}}

\def\p{\partial_z}
\def\tpsi{\tilde \psi}

\def\MMs{Manin\ matrices}
\def\MM{Manin\ matrix}

\def\BX{$\Box$}
\def\PRF{ {\bf Proof.} }

\def\kol{{\mathfrak R}}
\def\s{sgn}

\def\GR97{\href{http://arxiv.org/abs/q-alg/9705026}{[GR97]}}
\def\KL94{\href{http://arxiv.org/abs/q-alg/9411194}{[KL94]}}

\def\Ok961{\href{http://arxiv.org/abs/q-alg/9602028}{[Ok96]}}

\def\T04{\href{http://arxiv.org/abs/hep-th/0404153}{[Ta04]}}

\def\BT07{\href{http://arxiv.org/abs/hep-th/0703124}{[BT07]}}
\def\DubrovinM2{\href{http://arxiv.org/abs/math/0311261}{[DM03]}}
\def\Molev02{\href{http://arxiv.org/abs/math.QA/0211288}{[Mo02]}}

\def\EnriquezRubtsov01{\href{http://arxiv.org/abs/math.AG/0112276}{[ER01]}}
\def\BabelonTalon02{\href{http://arxiv.org/abs/hep-th/0209071}{[BT02]}}
\def\Nazarov91{\href{http://www.springerlink.com/content/rt5r3313732p48j7}{[Na91]}}
\def\GelfandGelfandRetakhWilson02{\href{http://arxiv.org/abs/math.QA/0208146}{[GGRW02]}}

\def\EtingofPak06{\href{http://arxiv.org/abs/math/0608005}{[EP06]}}

\def\HaiLorenz06{\href{http://arxiv.org/abs/math/0603169}{[HL06]}}

\def\KonvalinkaPak06{\href{http://arxiv.org/abs/math/0607737}{[KPak06]}}

\def\FoataHan06{\href{http://arxiv.org/abs/math/0603463}{[FH06]}}

\def\Sklyanin92{\href{http://arxiv.org/abs/hep-th/9212076}{[Sk92b]}}

\def\AHH96{\href{http://arxiv.org/abs/solv-int/9603001}{[AHH96]}}

\def\Sah06{\href{http://arxiv.org/abs/math-ph/0609067}{[Sa06]}}

\def\GGM97{\href{http://arxiv.org/abs/hep-th/9707120}{[GGM97]}}

\def\BartocciFalquiPedroni03{\href{http://lanl.arxiv.org/abs/nlin/0307021}{[BaFP03]}}

\def\PetreraRagnisco07{\href{http://lanl.arxiv.org/abs/math-ph/0703044}{[PR07]}}

\def\AdlerMoerbeke80{\href{http://www.sciencedirect.com}{[AvM80]}}

\def\ReymanSemenov79{\href{http://www.springerlink.com/content/n7360721350h7621}{[RS79]}}

\def\Beaville90{\href{http://www.springerlink.com/content/jn12352983882un4}{[Be90]}}

\def\AdamsHarnadPreviato88{\href{http://projecteuclid.org/DPubS/Repository/1.0/Disseminate?view=body&id=pdf_1&handle=euclid.cmp/1104161743}{[AHP88]}}

\def\Frenkel95{\href{http://lanl.arxiv.org/abs/q-alg/9506003}{[Fr95]}}

\def\Reshetikhin92{\href{http://www.springerlink.com/content/t718750378552740}{[Re92]}}

\def\BT{\href{http://arxiv.org/abs/hep-th/0703124}{[BTa07]}}

\def\OrtizSommaDukelskyRombouts04{\href{http://arxiv.org/abs/cond-mat/0407429}{[OSDR04]}}

\def\FaddeevTakhdadjanSklyanin79{\href{http://www.springerlink.com/content/uk1268x006041782/?p=6121269a4db64b09bba3c959e06b3a30&pi=1}
{[FTS79]}}

\def\AdamsHarnadHurtubise93{\href{http://projecteuclid.org/DPubS?service=UI&version=1.0&verb=Display&handle=euclid.cmp/1104253285}
{[AHH93]}}

\def\KulishSklyanin81{\href{http://www.springerlink.com/content/750u2k2xg12w8886/}{[KS81]}}

\def\ArnaudonAvanCrampeDoikouFrappatRagoucy04{\href{http://arxiv.org/abs/math-ph/0406021}
{[AACDFR04]}}

\def\ACDFR05{\href{http://arxiv.org/abs/math-ph/0512037v3}
{[ACDFR05]}}

\def\Hayashi88{\href{http://www.springerlink.com/content/wp5nwm73g35147u5/}{[Ha88]}}

\def\Kontsevich97{\href{http://arxiv.org/abs/q-alg/9709040}{[K97]}}

\def\Nekrasov95{\href{http://lanl.arxiv.org/abs/hep-th/9503157}{[Ne95]}}
\def\EnriquezRubtsovA95{\href{http://lanl.arxiv.org/abs/alg-geom/9503010}{[ER95]}}

\def\OkounkovB96{\href{http://lanl.arxiv.org/abs/q-alg/9602027}{[Ok96B]}}

\def\NazarovB96{\href{http://lanl.arxiv.org/abs/q-alg/9601027}{[Na96]}}

\def\Hitchin87{\href{http://projecteuclid.org/euclid.dmj/1077305506}{[Hi87]}}

\def\TV07{\href{http://arxiv.org/abs/math.QA/0702277v1}{[TV07]}}
\def\OPS07{\href{http://arxiv.org/abs/0711.2821v2}{[OPS07]}}
\def\HM06{\href{http://arxiv.org/abs/math.QA/0606121v2}{[HM06]}}

\begin{document}
\begin{titlepage}
{}
\vskip 1.5cm

\centerline{\Large\bf Algebraic properties of \MMs{} II:}

\centerline{\Large\bf $q$-analogues and integrable systems.}

\vskip 1.5cm
 \centerline{\large A. Chervov${}^1$,
G. Falqui${}^2$,
V. Rubtsov${}^{1,3}$,
A. Silantyev${}^4$.
}
\bigskip\bigskip
\noindent
{${}^1$ Institute for Theoretical and Experimental Physics, 
25, Bol. Cheremushkinskaya,  117218, Moscow - Russia. E--mail:
chervov@itep.ru}\\
{${}^2$ Dipartimento di Matematica e Applicazioni,  Universit\`a di Milano-Bicocca, via R. Cozzi, 53, 20125 Milano - Italy. E--mail: gregorio.falqui@unimib.it }
\\ 
{${}^{1}$  LAREMA, UMR6093 du CNRS,  Universit\'e D'Angers, Angers, France. \\E--mail:  volodya@univ-angers.fr}
\\
{${}^4$ Graduate School of Mathematical Sciences, University of Tokyo, 
Komaba, Tokyo, 123-8914, Japan. E--mail: aleksejsilantjev@gmail.com }

\vskip 1.0cm
\centerline{\large \bf  Abstract}
\vskip 1cm\noindent
We study a natural $q$-analogue of a class of matrices with noncommutative entries, which were first
considered by Yu. I. Manin in 1988 in relation with quantum group theory, (called {\em Manin Matrices} in \cite{CFR08}) . These matrices we shall  call 
\qMMs (qMMs). They are defined, in the $2\times 2$ case, by the relations
\[
M_{21}M_{12}=qM_{12}M_{21}, \> M_{22}M_{12}=qM_{12}M_{22}, \> [M_{11}, M_{22}]=q^{-1}M_{21}M_{12}-qM_{12}M_{21}.
\]
They were already considered in the literature, especially in connection with the $q$-Mac Mahon master theorem \cite{GLZ}, 
and the $q$-Sylvester identities \cite{KoPa07}.
The main aim of the present paper is to give a full list and detailed proofs of algebraic properties of qMMs\  
known up to the moment and, in particular, to show that most of the basic theorems of linear algebras (e.g., Jacobi ratio 
theorems, Schur complement, the Cayley-Hamilton theorem and so on and so forth) have a straightforward  
counterpart for \qMMs. We also show how this classs of matrices 
fits within the theory of quasi-determninants of Gel'fand-Retakh and collaborators (see, e.g., \cite{GR91}).   
In the last sections of the paper, we frame our definitions within the tensorial 
approach to non-commutative matrices of the Leningrad school, and we show how the notion of \qMM\ is related to 
theory of Quantum Integrable Systems. 
\end{titlepage}

\tableofcontents

\section{Introduction }
It is well-known that  matrices with generically noncommutative
elements play a basic role in the theory of quantum integrability as well as other fields of Mathematical Physics.
In particular  the best known occurrences of such matrices (see, e,g, \cite{CP94, FT87} and the references quoted therein) 
are the theory of quantum groups, the theory of critical and off-critical 
phenomena in statistical mechanics models, and various instances of combinatorial problems.
More recently, applications of such structures in the realm of the theory of Painlev\`e equations was considered \cite{ReRu10}.
The present paper builds, at least in its bulk, on the results of \cite{CFR08}. 
In that paper, a specific class of non-commutative matrices 
(that is, matrices with entries in a non commutative algebra -- or ring) were considered, namely the so-called \MMs.
This class of matrices, which are nothing but matrices representing linear transformations  
of polynomial algebra generators, were first introduced in the seminal paper \cite{Manin} by Yu. I . Manin, 
but they attracted sizable attention only recently (especially in the problem of the so-called quantum spectral curve in the theory 
of Gaudin models \cite{CT06-1}).  
The defining relations for a (column) \MM\  $M_{ij}$ are:
\begin{enumerate}
\item Elements in the same column commute; 
\item Commutators of the cross terms are equal: $[M_{ij} ,M_{kl}] = [M_{kj} ,M_{il}]$
(e.g. $[M_{11},M_{22}] = [M_{21},M_{12}]$). 
\end{enumerate}
The basic claim, fully proven in \cite{CFR08} is that
most theorems of linear algebra hold true
for \MMs~ (in a form {\em identical} to that of the commutative case) although the set of \MMs\ is  fairly different from that of ordinary matrices -- for instance they do not form a ring over the base field.
Moreover in some examples the converse is also true, that is,
\MMs~ are the most general class of matrices
such that ``ordinary`` linear algebra holds
true for them. Apart from that, it is important to remark that the \MM\ structure appears, as it was pointed out in 
\cite{CF07} in many instances related with the theory of Quantum integrable systems, namely the theory of integrable spin chains 
associated with {\em rational} $R$ matrices.
The basic aim of the present work is to extend our analysis to a class of matrices, which we call {\em \qMMs,} whose defining 
relations are in some sort the $q$-analogues of the ones for \MMs. Namely, they read as:
\begin{enumerate}
 \item Entries of the same column $q$-commute, i.e., 
$M_{ij}M_{kj} = q^{-1}  M_{kj} M_{ij }$
\item The $q$-cross commutation relations:
$M_{ij}M_{kl} -  M_{kl} M_{ij}  = q^{-1} M_{kj} M_{il} - q M_{il} M_{kj}, (i<k, j<l)$ hold.
\end{enumerate}
These defining relations, which, as an expert reader has surely already remarked, can be considered as a ''half" of 
the defining relations for the {\em Quantum Matrices}, that is, for the elements of the quantum matrix group $GL_q(N)$. 
These matrices were already introduced in the literature; indeed they were considered in the papers \cite{GLZ, LT07} (under 
the name of {\em right quantum matrices}), where the $q$-generalization of the Mac-Mahon master theorem was proven 
in connection with the boson-fermion correspondence of quantum physics. Their results (as well as the analysis of the 
$q$-Cartier Foata type of matrices) were generalized in \cite{KoPa07}. A multiparameter right-quantum analogue 
of Sylvester's identity was found in \cite{Konvalinka07-2}, also in connection with \cite{KL94}. It is worthwhile to remark that a super-version of the Mac-Mahon master theorem as well as other interesting consequences thereof was discussed in \cite{MR09}.

We would like also to remark that in works of Gel'fand, Retakh and their co-authors
(see, e.g., \cite{GR91, GR92, GGRW02}) a comprehensive linear algebra theory for   
{\em generic} matrices with noncommutative entries was extensively developed. Their approach is based on the fundamental 
notion of {\em quasi-determinant}.
Indeed, in the "general non-commutative case" there is no natural definition of a ``determinant`` (one has $n^2$  "quasi-determinants" 
instead) and, quite often, analogues of the proposition of ''ordinary`` 
linear algebra are sometimes formulated in a completely different way. Nevertheless it is clear that
their results can be fruitfully specialized and applied to some questions here; indeed, 
in Section \ref{sec3} we make contact with this theory, and, in Appendix \ref{AppPropLDJLC} we give an example of how this formalism can be applied to \qMMs.   However, in the rest of the paper, we shall 
not use the formalism of Gel'fand and Retakh. preferring to stress the similarities between our case of \qMMs\ with the case 
of ordinary linear algebra. 

Here we shall give a systematic analysis of  \qMMs, starting from their basic definitions, with an eye towards 
possible application to the theory of $q$-deformed quantum integrable systems 
(that is, in the Mathematical Physics parlance, those associated with {\em trigonometric} $R$-matrices).
The detailed layout of the present paper is the following:

We shall start with the description of some elementary properties directly stemming from the definition of \qMM\ in Section \ref{sez2}. In particular, 
we shall point out that \qMMs\ can be considered both as defining a (left) action on the generators of a $q$-algebra, and as 
defining a (right) action on the generators of a $q$-Grassmann algebra. We would like to stress that these properties, albeit elementary,
will be instrumental in the direct extension of the ''ordinary'' proofs of some determinantal properties to our case.
In Section  \ref{sec3} and \ref{sec4}, after having considered the relations/differences of the \qMMs\ with Quantum group theory,  
we shall recall in Section the definition of the $q$-determinant and of $q$-minor of a \qMM. In particular, 
we shall state and prove the analogue of the Cauchy-Binet theorem about the multiplicative property of the 
$q$-determinant. We shall limit ourselves to consider the simplest case (of two \qMMs\ with the elements of the 
first commuting with those of the second). We would like to point out that a more ample discussion of
(suitable variants of) the Cauchy Binet formula 
and the Capelli identity in the noncommutative case can be found, besides \cite{CFR08}, 
in~\cite{CSS08} and in the recent preprint~\cite{AC12}.

Moreover, we shall discuss how the Cramer's rule can be generalized in our case, and prove 
the Jacobi ratio theorem and the Lagrange-Desnanot-Jacobi-Lewis Carroll formula (also known as Dodgson  {\em condensation} formula). 
Also, we shall anticipate how the notion of $q$-characteristic polynomial (whose properties will be further discussed 
in Section \ref{MatrSect}) for a \qMM\ $M$ should be defined. It is worthwhile 
to note that this is a specific case in which the $q$-generalization of some ``ordinary'' object requires a suitable choice among the classically 
equivalent possible definitions. Indeed, the meaningful notion of $q$-characteristic polynomial is referred to  the weighted 
sum of the principal $q$-minors of $M$.

Then we shall show a fundamental property of \qMMs, that is, closedness under matrix ($q$)-inversion. 
We shall then state the analogue of Schur's complement theorem as well as the Sylvester identities. Further, a glance to the 
$q$-Pl\"ucker relations will be addressed (see, for a full treatment in the case of generic non-commutative matrices, 
the paper ~\cite{Lauve04}).

The aim of Section \ref{MatrSect} is twofold. First we shall show, making use of the so-called {\em Pyatov's Lemma}, how 
the theory of \qMMs\ of rank $n$  can be framed within the tensor approach of the Leningrad School, i.e. interpreting $q$-minors and the 
$q$-determinant as suitable elements in the tensor algebra of $\mathbb{C}^n$. This approach will be helpful in reaching the second aim of the section, 
that is, establishing the Cayley-Hamilton theorem and the Newton identities for \qMMs. 

Section \ref{secIntSys} shows, using the formalism introduced in Section \ref{MatrSect}, how the theory of \qMMs\ fits within the scheme of the Quantum Integrable spin  systems of trigonometric type. We shall first show that, if $L(z)$ is a Lax matrix satisfying 
the Yang-Baxter $RLL=LLR$ relations (Eq.~\ref{RLL}) with a trigonometric $R$-matrix $R$, then the matrix
\[
 M=L(z) q^{-2z\frac{\partial}{\partial z}},
\]
is actually a \qMM, and, moreover, that the {\em quantum determinant} of $L(z)$ defined by the Leningrad School 
(see, e.g., \cite{FRT}) coincides with the $q$-determinant of the associated \qMM\ $M$. Then, generalizing the corresponding 
construction of \cite{CF07, CFR08}, we shall show how an alternative set of quantum mutually commuting quantities can be 
obtained considering (ordinary) traces of suitably defined ``quantum'' powers of $M$.

The Appendix contains two subsections. The first is devoted to the detailed proof of the most technical Lemma of Section \ref{MatrSect}, which we extensively use in Section \ref{secIntSys}. In the second, we shall give an alternative proof  of the Dodgson  condensation formula making use of the formalism similar to that of quasi-determinants.

As a general strategy to keep  the present paper within a reasonable size, we shall present detailed proofs only in the case when these proofs are substantially different from the "non-$q$" (that is, $q\neq 1$) case. Otherwise, we shall refer to the corresponding proofs in \cite{CFR08}.

\subsection*{Acknowledgments.} We would like to thank D. Gurevich, P. Pyatov, V. Retakh, Z. Skoda and D. Talalaev for many useful discussions.
V.R. acknowledges financial support from the  project ANR  "DIADEMS", the project RFBR-01-00525a and the project "Cogito" (EGIDE).

The work of G.F. was partially supported by the Italian PRIN2008 project n. 20082K9KXZ. Also, we thank 
INdAM and SISSA for financial support during V.R.'s visits in Italy.

\section{\qMMs }\label{sez2}

\subsection{Definition of \qMMs}

Let $M$ be a $n\times m$ matrix with entries $M_{ij}$ in  an associative algebra 
$\mathfrak R$ over $\CC$ and let $q$ be a non-zero complex number. 
The algebra is not commutative in general, so that the matrix $M$ has non-commutative entries.

\begin{Def}\label{D1a22}
The matrix $M$ is called a $q-${\rm Manin} matrix if the following conditions hold true:
\begin{enumerate}
\item Entries of the same column $q$-commute with each other according to the order of the row indices:
\begin{eqnarray} \label{crossComRel110}
\forall j, i<k: ~~M_{ij}M_{kj} = q^{-1}  M_{kj} M_{ij },
\end{eqnarray}
\item The cross commutation relations:
\begin{eqnarray}
\label{crossComRel11}
\forall i<k, j<l:~~ M_{ij}M_{kl} -  M_{kl} M_{ij}  = q^{-1} M_{kj} M_{il} - q M_{il} M_{kj}
\end{eqnarray}
hold.
\end{enumerate}
\end{Def}

The defining relations for \qMMs{} can be compactly written as follows:
\begin{align}
\label{MMdefSingleFml}
M_{ij}M_{kl}-q^{\s(i-k)}q^{-\s(j-l)}M_{kl}M_{ij}=q^{\s(i-k)} M_{kj}M_{il}- q^{-\s(j-l)}M_{il}M_{kj},
\end{align}
for all $i,k=1,\ldots,n$, $j,l=1,\ldots,m$, where we use the notation
\begin{align}
 \s(k)=\begin{cases}
                  +1, &\text{if $k>0$;} \\
                   0, &\text{if $k=0$;} \\
                  -1, &\text{if $k<0$.}
                \end{cases}
\end{align}
Indeed, for $j= l$ 
and $i\ne k$ one gets the column $q$-commutativity~\eqref{crossComRel110}, while 
for $j\ne l$ and $i\ne k$ one gets the cross commutation relations~\eqref{crossComRel11} and for $i=k$ one gets an identity.

\begin{Rem} \label{remcrossComRel11} \normalfont
Definition~\ref{D1a22} can be reformulated as follows: the matrix $M$ is $q$-Manin matrix if and only if each $2\times2$ submatrix is a $q$-Manin matrix. More explicitly, for any $2\times 2$ submatrix $(M_{(ij)(kl)})$, consisting of rows $i$ and $k$, and columns $j$ and $l$ (where $1 \le i < k \le n$, and $1 \le j < l \le m$)
\begin{eqnarray}
\left(\begin{array}{ccccc}
... & ...& ...&...&...\\
... & M_{ij} &... & M_{il} & ... \\
... & ...& ...&...&...\\
... & M_{kj} &... & M_{kl}& ... \\
... & ...& ...&...&...
\end{array}\right)
\equiv
\left(\begin{array}{ccccc}
... & ...& ...&...&...\\
... & a &... & b& ... \\
... & ...& ...&...&...\\
... & c &... & d& ... \\
... & ...& ...&...&...
\end{array}\right),
\end{eqnarray}
the following commutation relations hold:
\begin{align}
\label{eq.commutation1}
ca &=  qac  &&\text{\small($q$-commutation of the entries in the $j$-th column),} \\
\label{eq.commutation2}
 db &= qbd &&\text{\small($q$-commutation of the entries in the $l$-th  column),} \\
\label{eq.commutation3}
ad - da &= q^{-1}cb-q bc &&\text{\small(cross commutation relation).}
\end{align}
\end{Rem}



{\bf Examples:}
\begin{itemize}
\item Let $\mathfrak R$ be the algebra generated by $a,b,c,d$ over 
$\mathbb C$ with the relations~\eqref{eq.commutation1}, \eqref{eq.commutation2}, \eqref{eq.commutation3} 
(or be an algebra containing some elements $a,b,c,d$ satisfying these relations). Then the matrix
\begin{align} \label{Mabcd}
 M=\begin{pmatrix}a&b\\c&d
   \end{pmatrix}
\end{align}
is a $q$-Manin matrix (over $\mathfrak R$). 
 
\item Let us consider $n$ elements $x_i\in\mathfrak R$, $i=1,\ldots,n$ that 
$q$-commute: $x_i x_j = q^{-1} x_j x_i $ for $i<j$ 
(e.g. we can consider the algebra $\mathfrak R$ generated by $x_i$ with these relations). 
Then the column-matrix
\begin{align} \label{ExManMat11x}
M=
\left(\begin{array}{cccc}
x_1 \\
x_2 \\
...\\
x_n \\
\end{array}\right)
\end{align}
can be considered a $q$-Manin matrix.

\item The elements from the same row of a $q$-Manin matrix are not required to satisfy any relations. So an arbitrary $1\times n$ matrix $M=(r_1,\ldots.,r_n)$ can be considered as a \qMM.

\item If some elements $m_i$ $q$-commute, i.e. $m_im_j=q^{-1}m_jm_i$, $i<j$, then the matrix
\begin{eqnarray}
\left(\begin{array}{ccccc}
m_1 & m_1 \\
m_2 & m_2 \\
... & ...\\
m_n & m_n
\end{array}\right),
\end{eqnarray}
is a \qMM. 
The cross-commutation relation~\eqref{crossComRel11} follows from $q$-commutativity in this case. Moreover, if $q\ne-1$ then the cross-commutation relations~\eqref{crossComRel11} for this matrix implies $q$-commutativity of $m_i$.

\item Consider elements $x_i,r_i\in\mathfrak R$, $i=1,\ldots,n$, 
such that $x_i x_j = q^{-1} x_j x_i $ for $i<j$ and $[x_i,r_j]=0$ $\forall i,j $. Then the matrix
\begin{eqnarray} \label{ExManMat11}
M=
\left(\begin{array}{cccc}
x_1 r_1 &  x_1 r_2 & ...& x_1 r_m \\
x_2 r_1 & x_2 r_2 & ...& x_2 r_m \\
... &...& ... & ...  \\
x_n r_1 & x_n r_2 & ...& x_n r_m \\
\end{array}\right)
=
\left(\begin{array}{cccc}
x_1 \\
x_2 \\
...\\
x_n \\
\end{array}\right)
\left(\begin{array}{cccc}
r_1 & r_2 & ...& r_m \\
\end{array}\right).
\end{eqnarray}
is a \qMM{}. This fact  can be easily checked by direct calculation. 

\item We refer to~\cite{CF07}, \cite{CFR08}, \cite{CM}, \cite{RST} for examples of $q$-Manin matrices for $q=1$ related to integrable systems, Lie algebras, etc. In the $q=1$ case the 
$q$-Manin matrices are simply called {\it Manin matrices}. 

\item Let $\mathbb C[x,y]$ be the algebra generated by $x,y$ with the relation $yx=qxy$. Consider the operators $\partial_x,\partial_y\colon\mathbb C[x,y]\to\mathbb C[x,y]$ of the corresponding differentiations: $\partial_x(x^ny^m)=nx^{n-1}y^m$, $\partial_y(y^mx^n)=my^{m-1}x^n$. Note that we have the following relations in $\End(\mathbb C[x,y])$: $\partial_y\partial_x=q\partial_x\partial_y$, $\partial_x y=qy\partial_x$, $\partial_y x=q^{-1}x\partial_y$, $[\partial_x,x]=[\partial_y,y]=1$, $q^{2y\partial_y}\partial_yq^{-2y\partial_y}=q^{-2}\partial_y$, where $q^{\pm2y\partial_y}=e^{\pm2\log(q)y\partial_y}$ are operators acting as $q^{\pm2y\partial_y}(x^ny^m)=q^{\pm2m}x^ny^m$. Then the matrix
\begin{eqnarray}
M=\left(\begin{array}{cc}
 x & q^{-1}q^{-2y\partial_y}\partial_y \\
 y & q^{-2y\partial_y}\partial_x
        \end{array}
  \right)
\end{eqnarray}
is a $q$-Manin matrix.

\item Examples of $q$-Manin matrices related to the quantum groups $Fun_q(GL_n)$ and $U_q(\wh{\mathfrak{gl}_n})$ will be considered in Subsections~\ref{QG} and~\ref{secIntSys1} respectively.

\end{itemize}

\begin{Rem} \label{tiledM} \normalfont If an $n\times m$  matrix $M$ is a \qMM, then the $n\times m$ matrix $\wt M$ with entries $\wt M_{ij}=M_{n-i+1,m-j+1}$ is a $q^{-1}$-Manin matrix. For example, if $M$ is the matrix~\eqref{Mabcd} then we obtain the following $q^{-1}$-Manin matrix:
\begin{eqnarray}
\wt M
=\left(\begin{array}{cc}
 d & c \\
b & a
\end{array}\right).
\end{eqnarray}
\end{Rem}

Let us note that matrices obtained by permutations
of rows or columns from \qMMs{} are {\em not} \qMMs{} in general.


\subsection{Linear transformations of $q$-(anti)-commuting variables}

The original definition of a Manin matrix \cite{Manin}, for the case $q=1$, was that of a matrix defining a linear transformation of commuting variables -- generators of the polynomial algebra -- as well as for a linear transformation of anti-commuting variables -- generators of the Grassmann algebra. This was largely used in~\cite{CF07,CFR08}.
Let us consider here an analogous interpretation of the \qMMs.

Let us first introduce a $q$-deformation of the polynomial algebra and a $q$-deformation of the Grassmann algebra. 
Define the {\it $q$-polynomial algebra} $\CC[x_1,\ldots,x_m]$ as the algebra with  generators $x_i$, $i=1,\ldots,m$, 
and relations $$x_i x_j=q^{-1} x_j x_i$$ for $i<j$. 
These relations can be rewritten in the form
\begin{align}
x_i x_j=q^{\s(i-j)} x_j x_i,&&& i,j=1,\ldots,m.
\end{align}
The elements $x_i$ are called {\it $q$-polynomial variables}. 
Similarly define the {$q$-Grassmann algebra} $\CC[\psi_1,\ldots,\psi_n]$ as the algebra generated by $\psi_i$, $i=1,\ldots,n$, with relations $\psi_i^2=0$ for $i=1,\ldots,n$ and $\psi_i \psi_j= -q\psi_j \psi_i$, for $i<j$; that is
\begin{align}
\psi_i \psi_j=-q^{-\s(i-j)} \psi_j \psi_i,&&& i,j=1,\ldots,n.
\end{align}
The elements $\psi_i$ are called {\it $q$-Grassmann variables}. 

Let $M$ be a rectangular $n\times m$-matrix over $\mathfrak R$. We can always suppose that $\mathfrak R$ contains the $q$-polynomial algebra $\CC[x_1,\ldots,x_m]$ and the $q$-Grassmann algebra $\CC[\psi_1,\ldots,\psi_n]$ as subalgebras such that the elements of these subalgebras commute with the entries of $M$: $[x_j,M_{pq}]=[\psi_i,M_{kl}]=0$ for all $i,k=1,\ldots,n$ and $j,l=1,\ldots,m$. Indeed, if $\mathfrak R$ does not contain one of these algebras we can regard the matrix $M$ as a matrix over the algebra $\mathfrak R\otimes\mathbb C[x_1,\ldots,x_m]\otimes\CC[\psi_1,\ldots,\psi_n]$, where the elements of different tensor factors pairwise commute.

\begin{Prop} \label{proposizione} 
Let the entries of a rectangular $n\times m$ matrix $M$ commute with the variables $x_1,\ldots,x_m$ and $\psi_1,\ldots,\psi_n$. 
Consider new $q$-polynomial variables $\tilde x_1,\ldots,\tilde x_n\in\mathfrak R$ and new $q$-Grassmann variables $\tilde\psi_1,\ldots,\tilde\psi_m\in\mathfrak R$ obtained by left (right) `action' of $M$ on the old variables:
\begin{align}
\left(\begin{array}{c}
\tilde x_1 \\
... \\
\tilde x_n
\end{array}\right)
&=
\left(\begin{array}{ccc}
 M_{11} & ... &  M_{1m} \\
... \\
 M_{n1} & ... &  M_{nm}
\end{array}\right)
\left(\begin{array}{c}
 x_1 \\
... \\
 x_m
\end{array}\right) \label{xtilde} \\
(\tilde \psi_1, \ldots ,  \tilde \psi_m) &=
(\psi_1, \ldots ,   \psi_n)
\left(\begin{array}{ccc}
 M_{11} & ... &  M_{1m} \\
... \\
 M_{n1} & ... &  M_{nm}
\end{array}\right), \label{psitilde}
\end{align}
Then the following three conditions are equivalent:
\begin{itemize}
\item The matrix $M$ is a \qMM.
\item The variables $\tilde x_i$ $q$-commute: $\tilde x_i \tilde x_j=q^{\s(i-j)} \tilde x_j \tilde x_i$ for all $i,j=1,\ldots,n$.
\item The variables $\tilde \psi_i$ $q$-anticommute: $\tpsi_i \tpsi_j=-q^{-\s(i-j)} \tpsi_j \tpsi_i$ for all $i,j=1,\ldots,m$.
\end{itemize}
\end{Prop}

\begin{Ex} \normalfont Let $x_1,x_2,\psi_1,\psi_2\in\mathfrak R$ be elements such that $x_1x_2=q^{-1}x_2x_1$, $\psi_1\psi_2=-q\psi_2\psi_1$, $\psi_1^2=\psi_2^2=0$. In the case $n=2$ the formulae~\eqref{xtilde}, \eqref{psitilde} take the form
\begin{align}
 \tilde x_1&=a x_1 + b x_2, & \tilde \psi_1&=a \psi_1 + c \psi_2, \\
 \tilde x_2&=c x_1 + d x_2, & \tilde \psi_2&=b \psi_1 + d \psi_2,
\end{align}
where $a,b,c,d\in\mathfrak R$ are elements commuting with $x_1,x_2$ and with $\psi_1,\psi_2$.
The matrix $M=\big(\begin{smallmatrix} a & b \\ c & d \end{smallmatrix}\big)$ is $q$-Manin if and only if $\tilde x_1\tilde x_2=q^{-1}\tilde x_2\tilde x_1$, or if and only if $\tilde\psi_1\tilde\psi_2=-q\psi_2\tilde\psi_1$ and $\tilde\psi_1^2=\tilde\psi_2^2=0$.
\end{Ex}
Both factors of the algebra $\mathfrak{R}\otimes \CC[x_1,\ldots,x_n]$ have a natural grading, and 
so \qMMs{} can be interpreted as matrices of grading-preserving homomorphisms $\CC[x_1,\ldots,x_n]\to\mathfrak R\otimes\CC[x_1,\ldots,x_n]$  
with respect to the variables $x_i$, and/or matrices of grading-preserving homomorphisms $\CC[\psi_1,\ldots,\psi_n]\to\mathfrak R\otimes\CC[\psi_1,\ldots,\psi_n]$ with respect to the variables $\psi_i$.

\begin{Rem} \label{rem3} \normalfont The conditions $\tpsi_i^2=0$ are equivalent to the relations~\eqref{crossComRel110} and the conditions $\tpsi_i \tpsi_j=-q\tpsi_j \tpsi_i$, $i<j$, are equivalent to the relations~\eqref{crossComRel11}.
\end{Rem}

\subsection{Elementary properties}

\begin{Prop} The following properties hold: \label{ElemPropPrpn}
\begin{enumerate} \slshape 
 \item \label{ElemPropPrpnItem1} Any submatrix of a \qMM{} is also a \qMM.
 \item \label{ElemPropPrpnItem2} Any diagonal matrix with commuting entries is a \qMM.
 \item \label{ElemPropPrpnItem3} If $M$ and $N$ are \qMMs{}, then $M+N$ is a \qMM{} if and only if
 \begin{align} \label{MplusN}
  &M_{ij}N_{kl}-q^{\s(i-k)}q^{-\s(j-l)}M_{kl}N_{ij}-q^{\s(i-k)} M_{kj}N_{il}+ q^{-\s(j-l)}M_{il}N_{kj}+ \\
    &+ N_{ij}M_{kl}-q^{\s(i-k)}q^{-\s(j-l)}N_{kl}M_{ij}-q^{\s(i-k)} N_{kj}M_{il}+ q^{-\s(j-l)}N_{il}M_{kj}=0, \notag
 \end{align}
 for all $i,j,k,l=1,\ldots,n$; in particular, if
 \begin{align} \label{MplusN2}
  M_{ij}N_{kl}=q^{\s(i-k)}q^{-\s(j-l)}N_{kl}M_{ij},
\end{align}
for all $i,j,k,l=1,\ldots,n$ then $M+N$ is a \qMM{};
 \item \label{ElemPropPrpnItem4} The product $cM$ of \qMM{} $M$ and a complex constant $c\in\mathbb C$ is also a \qMM{}.
 \item \label{ElemPropPrpnItem5} The product of \qMM{} $M$ and a diagonal complex matrix $D$ from the left $DM$ as well as from the right $MD$ is also a \qMM{}.

\item \label{ElemPropPrpnItem6}
\label{lemProdItem}
Let $M$ and $N$ be $n\times m$ and $m\times r$ \qMMs{} over an algebra $\mathfrak R$ such that their elements commute, i.e. $[M_{ij}, N_{kl}]=0$, $i=1,\ldots,n$, $j,k=1,\ldots,m$, $l=1,\ldots,r$, then the product $MN$ is a \qMM.
\end{enumerate}
\end{Prop}

\noindent{\bf Proof.} The first two properties (\ref{ElemPropPrpnItem1} and \ref{ElemPropPrpnItem2}) are obvious. To prove  property~\ref{ElemPropPrpnItem3} one should write the relations~\eqref{MMdefSingleFml} for $M+N$. Taking into account the relations~\eqref{MMdefSingleFml} for $M$ and $N$ one arrives to the relation~\eqref{MplusN}. The condition~\eqref{MplusN2} implies the condition~\eqref{MplusN}.

 The properties~\ref{ElemPropPrpnItem4} and \ref{ElemPropPrpnItem5} are particular cases of the property~\ref{ElemPropPrpnItem6}. To prove the last one we consider $q$-commuting variables~\footnote{One can just as well use $q$-Grassmann variables $\psi_i$.} $x_i$, $i=1,\ldots,r$, commuting with $M$ and $N$: $x_ix_j=q^{sgn(i-j)}x_jx_i$, $[M_{ij},x_l]=[N_{ij},x_l]=0$. Due to Proposition~\ref{proposizione} new variables $x^N_k=\sum_{l=1}^rN_{kl}x_l$, $k=1,\ldots,m$, $q$-commute: $x^N_ix^N_j=q^{sgn(i-j)}x^N_jx^N_i$. Since $[M_{ij},x^N_k]=0$ we can apply Proposition~\ref{proposizione} to the matrix $M$ and the variables $x^N_k$ -- the variables $x^{MN}_i=\sum_{j=1}^mM_{ij}x^N_j$, $i=1,\ldots,n$, also $q$-commute: $x^{MN}_ix^{MN}_j=q^{sgn(i-j)}x^{MN}_jx^{MN}_i$. Then, the formula $x^{MN}_i=\sum_{l=1}^r(MN)_{il}x_l$ and Proposition~\ref{proposizione} imply that $MN$ is a \qMM{}. \qed

\begin{Rem} \normalfont
 If $M$ and $N$ are matrices over the algebras $\mathfrak R_M$ and $\mathfrak R_N$ respectively we can consider $M$ and $N$ as matrices over the same algebra $\mathfrak R_M\otimes\mathfrak R_N$ and we then have the condition $[M_{ij},N_{kl}]=0$ for all $i,j,k,l$.
\end{Rem}


\subsection{Relations with quantum groups 
 \label{QG} }
One can also define $q$-analogues of \MMs{} characterizing the connections to quantum group theory. Actually $q$-\MMs{} are defined by half of the relations of the corresponding quantum group $Fun_q(GL_n)$\footnote{More  precisely we should write $Fun_q(\Mat_n)$, since we do not localize the $q$-determinant.} (\cite{FRT}). The remaining half consists of relations insuring that $M^\top$ is also a $q$-Manin matrix, where $M^\top$ is the transpose of $M$. 

\begin{Def} \label{QGdef} An $n \times n$ matrix $T$ belongs to the
quantum group $Fun_q(GL_n)$ if the following conditions hold true.
For any $2\times 2$  submatrix $(T_{(ij)(kl)})$, consisting of rows $i$
and $k$, and columns $j$ and $l$ (where $1 \leq i < k \leq n$, and
$1 \leq j < l \leq n$): \begin{eqnarray} \left(\begin{array}{ccccc}
... & ...& ...&...&...\\
... & T_{ij} &... & T_{il} & ... \\
... & ...& ...&...&...\\
... & T_{kj} &... & T_{kl}& ... \\
... & ...& ...&...&...
\end{array}\right)
\equiv \left(\begin{array}{ccccc}
... & ...& ...&...&...\\
... & a &... & b& ... \\
... & ...& ...&...&...\\
... & c &... & d& ... \\
... & ...& ...&...&...
\end{array}\right)
\end{eqnarray}
the following commutation relations hold:
\begin{eqnarray}
ca &= & qac, \quad \text{($q$-commutation of the entries in a column)} \\
 db &=& qbd, \quad \text{($q$-commutation of the entries in a column)} \\
ba &= & qab, \quad \text{($q$-commutation of the entries in a row)} \\
 dc &=& qcd, \quad \text{($q$-commutation of the entries in a row)} \\
ad - da &=& +q^{-1}cb-q bc, \qquad \text{(cross commutation relation 1)} \\
bc  &=& cb, \qquad \text{(cross commutation relation 2)}.
\end{eqnarray}
\end{Def}

As quantum groups are usually defined within the so-called matrix
(Leningrad) formalism, let us briefly recall it. (We will further
discuss this issue in Section \ref{MatrSect}).

\begin{Lem}
 The commutation relations for quantum group matrices can be described in matrix (Leningrad) notations as follows: \begin{eqnarray} R (T\otimes 1) (1\otimes T) =(1\otimes T) (T\otimes 1) R, \end{eqnarray} The $R$-matrix is given, in the case we are considering, by the formula: \begin{eqnarray} R= q^{-1} \sum_{i=1,..,n} E_{ii}\otimes E_{ii} +\sum_{i, j=1,..,n; i\ne j} E_{ii}\otimes E_{jj} +(q^{-1}-q) \sum_{i,j=1,..,n;i>j} E_{ij}\otimes E_{ji}, \end{eqnarray} where $E_{ij}$ are the standard basis of $\End(\mathbb C^n)$, i.e. $(E_{ij})_{kl}=\delta_{ik}\delta_{jl}$ -- zeroes everywhere except $1$ in the intersection of the $i$-th row with the $j$-th column.
\end{Lem}

For example in the $2 \times 2$ case the $R$-matrix is:
\begin{eqnarray} R= \left(\begin{array}{ccccc}
q^{-1}  & 0        & 0 & 0 \\
0       & 1        & 0 & 0 \\
0       & q^{-1}-q & 1 & 0 \\
0       & 0        & 0 & q^{-1}
\end{array}\right).
\end{eqnarray}

\begin{Rem} \normalfont
 This $R$-matrix differs by the change $q\to q^{-1}$ from
that one in~\cite{FRT}, page 185.
\end{Rem}


The relation between $q$-\MMs{} and quantum groups consists in the
following simple  proposition:

\begin{Prop}
 A matrix $T$ is a matrix in the quantum group $Fun_q(GL_n)$ if and only if both $T$ and $T^\top$ are $q$-\MMs.
\end{Prop}

So one gets that \qMMs\  are characterized by a "half" of
the conditions of the corresponding  quantum matrix
group.


\subsection{Hopf structure } \label{HopfSec} 

Let us consider the algebra $\mathbb C[M_{ij}]$ generated (over $\CC$) by $M_{ij}$, $1\le i,j \le n$, with $q$-Manin relations~\eqref{MMdefSingleFml}. One
can see that it can be equipped by a structure of bialgebra with the coproduct
$\Delta(M_{ij})=\sum_k M_{ik}\otimes M_{kj}$. This is usually
denoted as follows:
\begin{eqnarray}
\Delta (M) = M\stackrel{.}{\otimes}  M.
\end{eqnarray}
 It is easy to see that this coproduct is coassociative
(i.e. $ (\Delta \otimes 1)\otimes \Delta =  (1 \otimes
\Delta)\otimes \Delta$).

The natural antipode for this bialgebra should be $S(M)=M^{-1}$. So
it exists only in some extensions of the algebra $\mathbb C[M_{ij}]$.

The "coaction"-proposition~\pref{proposizione} implies that there exist
morphisms of algebras:
\begin{eqnarray} \phi_1: \CC[x_1,\ldots,x_m]\to \CC[M_{ij}]
\otimes \CC[x_1,\ldots,x_m], ~~~
\phi_1(x_i) = \sum_k M_{ik} x_k,\\
\phi_2: \CC[\psi_1,\ldots,\psi_n] \to \CC[M_{ij}]\otimes \CC[\psi_1,\ldots,\psi_n], ~~~
\phi_2(\psi_i) = \sum_k M_{ki} \psi_k
.
\end{eqnarray}
One can check that both maps satisfy
the condition: $(\Delta \otimes 1)(\phi_i) =  (1 \otimes \phi_i )(\phi_i)$, $i=1,2$.

So one can consider the maps $\phi_i$ as "coactions" of \qMMs{} on
the space $\CC^n_q$ and its  Grassmannian version.


\section{The $q$-determinant and Cramer's formula 
}
\label{sec3}
It was shown in~\cite{CF07,CFR08} that the natural generalization of the usual determinant for Manin matrices (i.e., the $q=1$ case) is the column determinant. This column determinant satisfies all 
the properties of the determinant of square matrices over a commuting field, and is defined as in the usuals case, with the proviso in mind that the orders of the column index in the $n!$ summands of the determinant should be always the same.
The role of column determinant for \qMMs{} is played by its $q$-analogue called {\it $q$-determinant}. 
Most of the properties of column determinant of Manin matrices presented in~\cite{CFR08} can be generalized to general $q$.

In this section we recall the definition of the $q$-determinant.
We will see that in the case of $q-$Manin matrices it generalizes  the notion of the usual determinants for the matrices over commutative rings.
We shall start by considering the $q$-determinant for an arbitrary (i.e., not necessarily $q$--Manin) matrix with elements in a noncommutative ring $\kol$.

\subsection{The $q$-determinant}

We define the $q$-determinant of an arbitrary matrix $M$ with non-commutative entries and relate it with coaction of $M$ on the $q$-Grassmann algebra. 
We prove here some formulae used below for $q$-determinants and $q$-minors of \qMMs{}.

\subsubsection{Definition of the $q$-determinant}

\begin{Def} \label{Def3}
The  {\em $q$-determinant},
of  a 
$n\times n$ matrix $M=(M_{ij})$ is  defined by the formula
\begin{eqnarray}
\det_q M:=
\sum_{\tau \in \sg_n} (-q)^{-\inv(\tau)} M_{\tau(1)1} M_{\tau(2)2}
\cdots M_{\tau(n)n} , \label{Def3eq}
\end{eqnarray}
where the sum ranges over the set $\sg_n$ of all permutations of
$\{1, \dots , n\}$.
Recall that $\inv(\tau)$ is number of inversions --
 the number of pairs $1 \leq i<j \leq n$
for which $\tau(i)>\tau (j)$. In other words, $\inv(\tau)$ is the length of $\tau$
with respect to adjacent transpositions $\sigma_k=\sigma_{(k,k+1)}$.
\end{Def}

\begin{Ex} \normalfont
\begin{eqnarray}\label{qdetmat2}
\det_q
\left(\begin{array}{ccc}
 a & b \\
c & d
\end{array}\right) \stackrel{def}{=} ad - q^{-1} c b
.
\end{eqnarray}
\end{Ex}
\begin{Ex} \normalfont
\begin{gather}
\det_q
\left(\begin{array}{ccc}
 M_{11} & M_{12} &  M_{13}  \\
 M_{21} & M_{22} &  M_{23}  \\
 M_{31} & M_{32} &  M_{33}  \\
\end{array}\right) \stackrel{def}{=}
M_{11} M_{22} M_{33}+   q^{-2} M_{21}M_{32}M_{13}
+ q^{-2} M_{31}M_{12}M_{23}  - \nn \\
 \nn \\
-   q^{-1} M_{11}M_{32}M_{23} - q^{-1} M_{21}M_{12}M_{33}
 -q^{-3} M_{31} M_{22} M_{13}
= \\
\nn \\
 = \nn
M_{11} \det_q
\left(\begin{array}{ccc}
 M_{11} & M_{12}   \\
 M_{21} & M_{22}   \\
\end{array}\right)
+
(-q)^{-1} M_{21} \det_q
\left(\begin{array}{ccc}
 M_{12} & M_{13}   \\
 M_{32} & M_{33}   \\
\end{array}\right)
+ \\
\nn \\
(-q)^{-2} M_{31} \det_q
\left(\begin{array}{ccc}
 M_{12} & M_{13}   \\
 M_{22} & M_{23}   \\
\end{array}\right)
 \nn
= (-q)^{-2}
\det_q
\left(\begin{array}{ccc}
 M_{21} & M_{22}   \\
 M_{31} & M_{32}   \\
\end{array}\right) M_{13}
+ \\
\nn \\
(-q)^{-1}
\det_q
\left(\begin{array}{ccc}
 M_{11} & M_{12}   \\
 M_{31} & M_{32}   \\
\end{array}\right) M_{23}
+
\det_q
\left(\begin{array}{ccc}
 M_{11} & M_{12}   \\
 M_{21} & M_{22}   \\
\end{array}\right) M_{33}
 . \label{Ex33detExp1}
\end{gather}
\end{Ex}

It is easy to see from the definition that, if $L$ is an arbitrary lower-triangular matrix with $d_i$ on the diagonal,
and $R$ an arbitrary upper-triangular matrix with $d_i$
on the diagonal, then 
\begin{eqnarray}
\det_q L= \det_q R= \prod_i d_i.
\end{eqnarray}
It is also easy to see that, if, $S$ is the permutation matrix corresponding to the permutation
$\sigma\in \sg_n$, then $det_q S=(-q)^{-\inv(\sigma)}$.
\begin{Rem}
 \normalfont For $q=1$ the present definition of $q$-determinant  coincides with that of  column-determinant i.e. the determinant defined by the column expansion, first taking elements from the first column, then from the second and so on and so forth.
For $q=-1$  the $q$-determinant yields  the (column) permanent of the matrix $M$.  
Just like the usual determinant, the $q$-determinant is linear over $\mathbb C$ both with respect to  its columns and rows.
%
%
\end{Rem}

\subsubsection{The $q$-Grassmann algebra }

Let us 
collect 
some useful relations in the $q$-Grassmann algebra $\mathbb C[\psi_1,\ldots,\psi_n]$.
Some of them hold without assuming $\psi_i^2=0$, while some of them require this property.

\begin{Lem} \label{lempsiProduct}
Let $\psi_i$ satisfy the relations $\psi_i \psi_j=-q\psi_j \psi_i$ for all $1\le i<j\le n$; then their monomials of $n$-th order are related as follows
\begin{align}
 \psi_{\tau(1)}\cdots\psi_{\tau(n)}&=
 (-q)^{-\inv(\tau)} 
 \psi_1\cdots\psi_n, & &\forall\tau\in \sg_n, \label{psiProduct}
\end{align}
or equivalently as
\begin{align}
 \psi_{\sigma\tau(1)}\cdots\psi_{\sigma\tau(n)}&=
 (-q)^{-\inv(\sigma\tau)+\inv(\sigma)}
 \psi_{\sigma(1)}\cdots\psi_{\sigma(n)}, &  & \forall\sigma,\tau\in \sg_n. \label{lemDetPsiR1}
\end{align}
\end{Lem}

\noindent{\bf Proof.} Note that if the relation~\eqref{lemDetPsiR1} is valid for some $\tau_1,\tau_2\in \sg_n$ (for all $\sigma\in \sg_n$) then it is so for $\tau=\tau_1\tau_2$. Indeed
\begin{multline}
\psi_{\sigma\tau_1\tau_2(1)}\cdots\psi_{\sigma\tau_1\tau_2(n)}=
 (-q)^{-\inv(\sigma\tau_1\tau_2)+\inv(\sigma\tau_1)}
 \psi_{\sigma\tau_1(1)}\cdots\psi_{\sigma\tau_1(n)}= \\
 =(-q)^{-\inv(\sigma\tau_1\tau_2)+\inv(\sigma)}
 \psi_{\sigma(1)}\cdots\psi_{\sigma(n)}. \label{lemDetPsiR2}
\end{multline}
Since each $\tau$ can be presented as a product of adjacent transpositions $\sigma_k$, $k=1,\ldots,n-1$, it is sufficient to prove~\eqref{lemDetPsiR1} for $\tau=\sigma_k$. In this case we can write
\begin{multline}
 \psi_{\sigma\sigma_k(1)}\cdots\psi_{\sigma\sigma_k(n)}=
 \psi_{\sigma(1)}\cdots\psi_{\sigma(k+1)} \psi_{\sigma(k)}\cdots\psi_{\sigma(n)}= \\
 =(-q)^{-\s(\sigma(k+1)-\sigma(k))}
 \psi_{\sigma(1)}\cdots\psi_{\sigma(k)} \psi_{\sigma(k+1)}\cdots\psi_{\sigma(n)}. \label{lemDetPsiR3}
\end{multline}
Thus,  formula~\eqref{lemDetPsiR1} for $\tau=\sigma_k$ follows from the equality
\begin{align}
 \inv(\sigma\sigma_k)=\inv(\sigma)+\s(\sigma(k+1)-\sigma(k)). \label{lemDetPsiR4}
\end{align}
\qed

\begin{Cor} \label{corpsiProduct}
Under the condition of the Lemma~\ref{lempsiProduct}, one can relate the monomials of $m$-th order as follows:
\begin{align}
 \psi_{j_{\tau(1)}}\cdots\psi_{j_{\tau(m)}}&=
 (-q)^{-\inv(\tau)}
 \psi_{j_1}\cdots\psi_{j_m}, \label{psiProduct_m} \\
 \psi_{j_{\tau(1)}}\cdots\psi_{j_{\tau(m)}}&=
 (-q)^{-\inv(\tau)+\inv(\sigma)}
 \psi_{j_{\sigma(1)}}\cdots\psi_{j_{\sigma(m)}}, \label{lemDetPsiR1_m}
\end{align}
where $1\le j_1<\ldots<j_m\le n$, $\tau,\sigma\in \sg_m$.
\end{Cor}

\noindent{\bf Proof.} Consider the elements $\tilde\psi_k=\psi_{j_k}$. They satisfy the conditions of Lemma~\ref{lempsiProduct}: $\tpsi_i \tpsi_j=-q\tpsi_j \tpsi_i$ for all $1\le i<j\le m$. Writing the formulae~\eqref{psiProduct}, \eqref{lemDetPsiR1} for $\tilde\psi_k$ we obtain~\eqref{psiProduct_m}, \eqref{lemDetPsiR1_m}. \qed

\begin{Cor}
 Assuming additionally $\psi_i^2=0$, i.e. $\psi_i\psi_j=-q^{-\s(i-j)}\psi_j\psi_i$, 
it is convenient to write the more general formula
\begin{align}
 \psi_{i_1}\cdots\psi_{i_n}&=\varepsilon^q_{i_1,\ldots,i_n}\psi_{1}\cdots\psi_{n}, \label{psiviapsi}
\end{align}
where 
 $1\le i_l\le n$ 
 and the $q$-epsilon-symbol is defined by the formula
\begin{align}
 \varepsilon^q_{\ldots,i,\ldots,i,\ldots}&=0, & \varepsilon^q_{\sigma(1),\ldots,\sigma(n)}&=(-q)^{-\inv(\sigma)}, & \sigma\in \sg_n. \label{vareps}
\end{align}
\end{Cor}

Let us denote by $I\oplus J=(i_1,\ldots,i_m,j_1,\ldots,j_k)$ the contraction of two multi-indices $I=(i_1,\ldots,i_m)$ and $J=(j_1,\ldots,j_k)$. For a multi-index $K=(k_1,k_2,\ldots,k_r)$ denote by $\backslash K$ the multi-index $(1,\ldots,\hat k_1,\ldots,\hat k_m,\ldots,n)$, that is $\backslash K$ is obtained from $(1,2,\ldots,n)$ by deleting $k_i$ for all $i=1,\ldots,r$.

Let $(i_1,\ldots,i_n)$ be a permutation of  $(1,\ldots,n)$ such that $i_1<i_2<\ldots<i_m$ and $i_{m+1}<i_{m+2}<\ldots<i_{n}$ for some $m\le n$. Let $I=(i_1,i_2,\ldots,i_m)$ and $\backslash I= (i_{m+1},i_{m+2},\ldots,i_{n})$. It is quite easy to see that
\begin{align}
\label{EpsIcI1}
\varepsilon^q_{(I \oplus \backslash I)}
&= (-q)^{-\sum_{l=1}^m (i_l-l)}
=(-q)^{+\sum_{l=m+1}^n (i_l-l)},\\
\label{EpsIcI2}
\varepsilon^q_{(\backslash I \oplus  I)}
&= (-q)^{\sum_{l=1}^m i_l-\sum_{l=n-m+1}^n l}
=(-q)^{\sum_{l=1}^{n-m} l-\sum_{l=m+1}^n i_l},
\end{align}
where we used $\sum_{l=1,\ldots,n} i_l = \sum_{l=1,\ldots,n}l$.

Let us also mention that for $I=(n,n-1,n-2,\ldots,1)$ we have
\begin{align}
\label{EpsOmega}
\varepsilon^q_{I}
=
(-q)^{-n(n-1)/2}.
\end{align}

\subsubsection{The $q$-determinant, $q$-minors and the $q$-Grassmann algebra \label{qdetGrasSS} }

Let us give a more conceptual approach to $q$-determinants and $q$-minors -- $q$-determinants of submatrices -- via the $q$-Grassmann algebra. 
We consider an arbitrary matrix $M$, but the fact that 
the notion of the $q$-determinant can reformulated in terms of the 
$q$-Grassmann algebra is a clear hint that it is related with the notion of \qMM{}.

Let $M$ be an arbitrary $n\times m$ 
matrix and $I=(i_1,i_2,\ldots,i_k)$ and $I=(i_1,i_2,\ldots,i_l)$ be two arbitrary multi-indices, where $k\le n$ and $l\le m$. Denote by $M_{IJ}$ the $k\times k$ matrix defined as $(M_{IJ})_{ab}=M_{i_a j_b}$, where $a,b=1,\ldots,k$. Then the $q-$determinant $\det_q(M_{IJ})$ can be considered as 
the $q$-analogue of the minor ({\it $q$-minor}) of the matrix $M$.

\begin{Prop}
\label{lemDetPsi}
Let $M$ be an  arbitrary
  (not-necessarily $q$-Manin)
$n\times n$   matrix,
$\psi_i$ -- $q$-Grassmann variables, i.e. $\psi_i \psi_j = -q^{-\s(i-j)} \psi_j \psi_i$, which commute with the entries of $M$,
and $\psi_i^M=\sum_k \psi_k M_{ki}$. Then
\begin{align}
\psi_1^M\psi_2^M\ldots\psi_n^M = \det_q(M)\psi_1\psi_2\ldots\psi_n. \label{lemDetPsi1}
\end{align}
More generally, for a rectangular $n\times m$ matrix $M$
 and for an arbitrary\footnote{Here it is {\em neither} required that $j_l<j_p$, nor that $j_l\ne j_p$.} multi-index
$J=(j_1,j_2,\ldots,j_k)$
\begin{align}
\psi_{j_1}^M\psi_{j_2}^M\ldots\psi_{j_k}^M
 = \sum_{I=(i_1<i_2<\ldots<i_k)}
  \det_q(M_{IJ})\psi_{i_1}\psi_{i_2}\ldots\psi_{i_k}, \label{lemDetPsi2}
\end{align}
where the sum is taken over all multi-indices $I=(i_1,i_2,\ldots,i_k)$ such that $1\le i_1<i_2<\ldots<i_k\le n$.
\end{Prop}

\noindent{\bf Proof.} The formula~\eqref{lemDetPsi1} follows from \eqref{lemDetPsi2} in the case $m=n=k$. To prove the last one we write left hand side as
\begin{align}
 \psi^M_{j_1}\cdots\psi^M_{j_k}&=\sum_{l_1,\ldots,l_k=1}^n M_{l_1j_1}\cdots M_{l_kj_k}\psi_{l_1}\cdots\psi_{l_k}.
\end{align}
The product $\psi_{l_1}\cdots\psi_{l_k}$ does not vanishes only if $l_1,\ldots,l_k$ is a permutation of some numbers $i_1,\ldots,i_k$ such that $1\le i_1<i_2<\ldots<i_k\le n$. Then using the formula~\eqref{psiProduct_m} one derives:
\begin{multline}
 \psi^M_{j_1}\cdots\psi^M_{j_k}=\sum_{I=(i_1<i_2<\ldots<i_k)}\sum_{\tau\in \sg_k} M_{i_{\tau(1)},j_1}\cdots M_{i_{\tau(k)},j_k}\psi_{i_{\tau(1)}}\cdots\psi_{i_{\tau(k)}}= \\
=\sum_{I=(i_1<i_2<\ldots<i_k)}\sum_{\tau \in \sg_k}(-q)^{-\inv(\tau)} M_{i_{\tau(1)},j_1}\cdots M_{i_{\tau(n)},j_k}\psi_{i_1}\cdots\psi_{i_k}.
\label{lemDetPsiR6}
\end{multline}
 On the other hand the $q$-determinant of $M_{IJ}$ can be written as
\begin{align}
\det_q(M_{IJ})=\sum_{\tau \in \sg_k}(-q)^{-\inv(\tau)} M_{i_{\tau(1)},j_1}\cdots M_{i_{\tau(n)},j_k}.\label{lemDetPsiR7}
\end{align}
Comparing~\eqref{lemDetPsiR6} with \eqref{lemDetPsiR7} one gets~\eqref{lemDetPsi2}.
\qed

\begin{Ex} \normalfont
 For $\psi_1 \psi_2 =- q \psi_2 \psi_1, \psi_i^2=0$ and    $M= \left(\begin{array}{ccc}
 a & b \\
c & d
\end{array}\right)$, one has
\begin{eqnarray}
\psi_1^M \psi_2^M=
(\psi_1 a+\psi_2 c) (\psi_1 b + \psi_2 d) =
\psi_1 \psi_2  a d + \psi_2 \psi_1 cb = \nn\\
=  \psi_1 \psi_2  a d - q^{-1} \psi_1 \psi_2  cb
= det_q(M) \psi_1 \psi_2.
\end{eqnarray}
\end{Ex}

\noindent A tautological corollary of the proposition
above is the following.
\begin{Cor}
Let $M^\sigma$ be the matrix obtained from $M$ by the permutation
$\sigma \in \sg_n$
of columns: $M^\sigma_{ij}=M_{i\sigma(j)}$.
Then
\begin{align}
\psi_{\sigma(1)}^M\psi_{\sigma(2)}^M\cdots\psi_{\sigma(n)}^M = \det_q(M^{\sigma})\psi_1\psi_2\cdots\psi_n. \label{CorR321}
\end{align}
\end{Cor}

In Proposition~\ref{lemDetPsi} it is crucial
that the elements of $M$ commute with $\psi_i$;
however $q$-Grassmann algebra is also helpful
without this condition. This can be seen from the
following formula which is proved in the same way
as Proposition~\ref{lemDetPsi}
\begin{align} 
\sum_{I=(i_1,\ldots,i_r)\atop 1\le i_a \le n} \psi_{i_1} \psi_{i_2}\cdots \psi_{i_r} M_{i_1j_1} M_{i_2j_2}\cdots M_{i_rj_r}
=\sum_{L=(l_1<l_2<\ldots<l_r)} \psi_{l_1} \psi_{l_2}\cdots
\psi_{l_r} \det_q(M_{LJ}). \notag
\end{align}
The Laplace expansion can also be written as follows: 
\begin{multline}
\label{GrVsMinFml3}
\sum_{L_1=(l_1<l_2<\ldots<l_a)}  \psi_{l_1}\cdots\psi_{l_a}
\quad \sum_{L_2=(l_{a+1}<l_{a+2}<\ldots<l_r)}  \psi_{l_{a+1}}\cdots
\psi_{l_r} \det_q A_{L_1 (j_1,\ldots,j_a)}\times \\
\times \det_q (A_{L_2 (j_{a+1},\ldots,j_r)})
=
\sum_{L=(l_1<l_2<\ldots<l_r)} \psi_{l_1} \psi_{l_2}\cdots
\psi_{l_r} \det_q(A_{LJ}).
\end{multline}

\subsection{Properties of the $q$-determinants of \qMMs}
\label{sec32}

The results of the previous subsection lead 
to the properties  of $q$-determinants listed below. Proposition~\ref{lemDetPsi} allows us to derive some properties for $q$-determinants of \qMM{}. Let us remark that some these properties are valid for matrices with entries satisfying a part of the $q$-Manin relations or even for arbitrary matrices.


\begin{Prop} \label{lem.3a} The following properties hold:

\begin{enumerate}  \slshape 

\item\label{qd1} {\rm Linearity in columns and rows.} \\
If some column (row) of a square non-commutative (not necessarily $q$-Manin) matrix $M$ is presented as a sum of column-matrices (rows-matrices) then its $q$-determinant $\det_q M$ is equal to the sum of $q$-determinants of the matrices $M$ with the considered column (row) replaced by the corresponding column-matrices (rows-matrices).

\item\label{qd2} {\rm Permutation of columns.}\\
If $M$ is a square matrix satisfying relations~\eqref{crossComRel11} (in particular, if $M$ is a \qMM) 
and $M^\sigma$ denotes the 
matrix obtained from $M$ by a permutation of columns: 
$M^\sigma_{ij}=M_{i\sigma(j)}$, where $\sigma \in \sg_n $, then
\begin{align}
 \det_q(M^\sigma)=(-q)^{-\inv(\sigma)} \det_q M. \label{qd2R}
\end{align}

Note that  $M^\sigma$ is not a \qMM{} in general (even if $M$ is a $q$-Manin matrix). 
Let us also remark that permutation of rows affects the \qDet in a highly non-trivial 
way since the shuffling of the rows destroys the property of being $q$-Manin. 

\item\label{qd3} {\rm Matrices with coincident columns.}\\
Let $M$ be a square \qMM.
If two columns of $M$ coincide, then
\begin{align}
\label{lemCoinCol}
\det_q M=0.
\end{align}
If $q\ne -1$ the same holds for any square $M$ satisfying relations~\eqref{crossComRel11}.
  
Furthermore, if $M$ is a square \qMM{} and the matrix $\wt M$ is obtained from $M$ by substituting the $r$-th column to the $s$-th column, where $s\ne r$, then
\begin{align}
\label{lemCoinColTilde}
\det_q \wt M =0,
\end{align}
note that if $|s-r|\ge2$ then $\wt M$ is not a \qMM{} in general.
\footnote{Actually this property holds under weaker conditions: relations~\eqref{crossComRel11} and $q$-commutativity of elements in the $r$-th column.}

In particular, the coincidence of {\it rows} do not imply the vanishing of the $q$-determinant.

\item\label{qd4} {\rm Determinant multiplicativity and Cauchy-Binet formula.} \label{CBItem} \\
Let $M$ and $N$ be two matrices such that $[M_{ij}, N_{kl}]=0$ for all possible indices $i,j,k,l$ and let $M$ be a $q$-Manin matrix. If these are $n\times n$ matrices, then
\begin{align}
\label{lemProd}
\det_q(M N)= \det_q(M) \det_q(N).
\end{align}
More generally, if $M$ and $N$ are rectangular matrices and $i_1<i_2<\ldots<i_r$ then  the {\it Cauchy-Binet formula} holds:
\begin{align}
\label{CBfml}
\det_q\big((M N)_{IJ}\big)= \sum_{L=(l_1<l_2<\ldots<l_r)}
 \det_q(M_{IL})\det_q(N_{LJ}),
\end{align}
where $I=(i_1,i_2,\ldots,i_r)$ and $J=(j_1,j_2,\ldots,j_r)$.

Recall that in this case $MN$ is a \qMM,
if $N$ is also a \qMM{} (see the item~\ref{lemProdItem} of Proposition~\pref{ElemPropPrpn}).  
\item\label{qd5} {\rm Column expansion.}\\
For an arbitrary $n\times n$ matrix $M$ (not necessarily a \qMM) the following  expansions with respect to the first and  last columns hold:

\begin{eqnarray}
\det_q(M)=
\sum_{r=1}^{n} (-q)^{1-r} M_{r1}
 \det_q{(M_{\backslash r \backslash 1})}
 =
\sum_{r=1}^{n} (-q)^{r-n}
 \det_q{(M_{\backslash r \backslash n})} M_{rn}, \label{qd5R1}
\end{eqnarray}
where  $M_{\backslash r \backslash s}$ is the $(n-1)\times(n-1)$ matrix obtained by deleting the $r$-th row and the $s$-th column.
(As an example see formulae \eqref{Ex33detExp1}).

For an $n\times n$ matrix $M$ satisfying relations~\eqref{crossComRel11} (in particular, if $M$ is a \qMM) its $q$-determinant 
can be decomposed (for any $s=1,\ldots, n$) along the $s-$ column as follows:
\begin{eqnarray}
\label{ColExpFml}
\det_q(M)=\sum_{r=1}^{n} (-q)^{s-r} M_{rs}
 \det_q{(M_{\backslash r \backslash s})}
=\sum_{r=1}^{n} (-q)^{r-s}
 \det_q{(M_{\backslash r \backslash s})} M_{rs} .
\end{eqnarray}
To the best of our knowledge, there are no analogous formulae for the row expansion.

\item\label{qd6} {\rm Laplace expansion.} \label{LaplItem} \\
For an arbitrary $n\times n$ matrix $M$ (not necessarily \qMM)
the following is true:
\begin{align}
\det_q M=\sum_{K=(k_1,\ldots,k_m)\atop1\le k_1<\ldots <k_m \le n}
(-q)^{-\sum_{l=1}^m(k_l-l)}
\det_q M_{ K,(1\ldots m)}
\det_q M_{ \backslash K,(m+1,\ldots, n) }.
\label{LaplFml1}
\end{align}

For an $n\times n$  \qMM{} $M$,
and  arbitrary pair of multi-indices $I_1=(i_1,\ldots,i_m)$ and $I_2= (i_{m+1},\ldots,i_n)$
one can also write:  
\begin{align}
 \varepsilon^q_{i_1,\ldots,i_n}
\det_q M=\sum_{K=(k_1<\ldots<k_m)}
(-q)^{-\sum_{l=1}^m(k_l-l)}
\det_q M_{K,I_1}
\det_q M_{\backslash K,I_2 },
\label{LaplFml2}
\end{align}
where the $q$-epsilon-symbol is defined in ~\eqref{vareps}.

One can also write similar formulae with products of more determinants in the right hand side. 
Consider $r$ non-negative numbers $l_1,\ldots,l_r$ such that $l_1+l_2+\ldots+l_r=n$ and arbitrary multi-indices $I_j=(i_{j1},i_{j2},\ldots,i_{j,l_j})$. Let $(i_1,\ldots,i_n)=I_1\oplus I_2\oplus\ldots\oplus I_r$, then
\begin{align} \label{LaplFml3}
 \varepsilon^q_{i_1,\ldots,i_n}
\det_q M=
\sum_{K_{1},\ldots, K_r}
 \varepsilon^q_{k_1,\ldots,k_n}
\prod_{j=1}^r \det_q M_{ K_j I_j},
\end{align}
where the sum is taken over all multi-indices $K_j=(k_{j1},k_{j2},\ldots,k_{j,l_j})$ such that $k_{j1}<k_{j2}<\ldots<k_{j,l_j}$ and $(k_1,\ldots,k_n)=K_1\oplus K_2\oplus\ldots\oplus K_r$.


If $(i_1,\ldots,i_n)=(1,2,3,\ldots,n)$ then the formula~\eqref{LaplFml3} holds for arbitrary square matrix. For $l_1=l_2=\ldots=l_n=1$ this formula becomes the definition of the $q$-determinant.

\item\label{qd7} {\rm Non centrality of \qDet.}\\
We remark  that $\det_q(M)$ (contrary to the case of Quantum Group theory)  is not a central element in the algebra 
$\mathbb C[M_{ij}]$ defined in the subsection~\ref{HopfSec};
moreover, in  general $[\det_q(M), \tr(M)]\ne 0$.
This difference between  \qMMs\ and quantum matrices can be traced back to the fact that the commutation relations between elements of 
a \qMM\ are more general than those that the elements of  a "quantum matrix" do satisfy.
\end{enumerate}
\end{Prop}

{\sf Proof of Property \sl\ref{qd1}}. Let us show the linearity in the first column and in the first row. The first one is a direct sequence of Definition~\ref{Def3}: let $M_{ij}=N_{ij}$ for $j\ne1$ then
\begin{align}
\det_q(M+N)=
\sum_{\tau \in \sg_n} (-q)^{-\inv(\tau)} (M_{\tau(1)1}+N_{\tau(1)1}) M_{\tau(2)2}
\cdots M_{\tau(n)n}=\det_q(M)+\det_q(N). \notag
\end{align}
Analogously, if $M_{ij}=N_{ij}$ for $i\ne1$, 
\begin{multline}
\det_q(M+N)=
\sum_{\tau \in \sg_n} (-q)^{-\inv(\tau)} M_{\tau(1)1} M_{\tau(2)2}\cdots(M_{1\tau^{-1}(1)}+N_{1\tau^{-1}(1)})
\cdots M_{\tau(n)n}= \\ =\det_q(M)+\det_q(N). \notag
\end{multline}
 \qed

{\sf Proof of Property \sl\ref{qd2}.} Let $\psi_i$ be $q$-Grassmann variables commuting with $M$, that is, $\psi_i\psi_j\hm=-q^{-\s(i-j)}\psi_j\psi_i$, $[\psi_i,M_{jk}]=0$, and $\psi^M_i=\sum_{k=1}^n \psi_k M_{ki}$. From Remark~\ref{rem3} we conclude that $\psi^M_i$ satisfy the condition of the lemma~\ref{lempsiProduct}. Applying the formula~\eqref{psiProduct} to the left hand side of~\eqref{CorR321} we obtain~\eqref{qd2R}. \qed

{\sf Proof of Property \sl\ref{qd3}.} Let the matrix $M$ satisfy~\eqref{crossComRel11}. Let suppose also that its $k$-th and $l$-th columns coincide. If $\psi^M_j=\sum_{i=1}\psi_i M_{ij}$ then $\psi^M_k=\psi^M_l$ and $M=M^{\sigma_{kl}}$, where $\sigma_{kl}$ is a permutation interchanging $k$ and $l$.
If $M$ is $q$-Manin the variables $\psi^M_j$ are $q$-Grassmann and hence the equality~\eqref{lemDetPsi1} implies $\det_q(M)=0$.
If $q\ne-1$ the formula $\det_q(M)=0$ is a consequence of the lemma~\ref{lem.3a} for $\sigma=\sigma_{kl}$.
Let $\psi^{\wt M}_j=\sum_{i=1}\psi_i M_{ij}$, then $\psi^{\wt M}_i\psi^{\wt M}_r=-q^{-\s(i-r)}\psi^{\wt M}_r\psi^{\wt M}_i$ and $\psi^{\wt M}_r=\psi^{\wt M}_s$. It is sufficient to obtain the formula $\det_q(\wt M)=0$ using the equality~\eqref{lemDetPsi1}. \qed

{\sf Proof of Property \sl\ref{qd4}.} Let $\psi_i$, $i=1,\ldots,n$, be $q$-Grassmann variables, which commute with the entries of $M$ and $N$. Since the matrix $M$ is $q$-Manin the variables $\psi^M_l=\sum_{i=1}^n \psi_i M_{il}$ are also $q$-Grassmann. Let $\psi^{MN}_j=\sum_{l=1}^m \psi^M_l N_{lj}=\sum_{i=1}^n \psi_i (MN)_{ij}$. Note that $\psi_i$, $\psi^M_l$ and $\psi_i^{MN}$ commute with $M$, $N$, $MN$ respectively. So that we can write the formula~\eqref{lemDetPsi2} for these matrices:
\begin{align}
\psi^M_{l_1}\ldots\psi^M_{l_r}&=\sum_{I=(i_1<\ldots<i_r)}\psi_{i_1}\ldots\psi_{i_r}\det_q(M_{IL}), \label{lemProdDetR1} \\
\psi^{MN}_{j_1}\ldots\psi^{MN}_{j_r}&=\sum_{L=(l_1<\ldots<l_r)}\psi^M_{l_1}\ldots\psi^M_{l_r}\det_q(N_{LJ}),\label{lemProdDetR2} \\
\psi^{MN}_{j_1}\ldots\psi^{MN}_{j_r}&=\sum_{I=(i_1<\ldots<i_r)}\psi_{i_1}\ldots\psi_{i_r}\det_q\big((MN)_{IJ}\big). \label{lemProdDetR3}
\end{align}
where $L=(l_1<\ldots<l_r)$ and $J=(j_1<\ldots<j_r)$. Substituting~\eqref{lemProdDetR1} to the right hand side of~\eqref{lemProdDetR2} and comparing the result with~\eqref{lemProdDetR3} we derive~\eqref{CBfml}. \qed

{\sf Proof of Property \sl\ref{qd5}.} The formula~\eqref{qd5R1} follows immediately from the definition of the $q$-determinant. For example,
\begin{multline}
 \det_q(M)=\sum_{\tau\in \sg_n}(-q)^{-\inv(\tau)}M_{\tau(1),1}\cdots M_{\tau(n),n}= \\
 =\sum_{r=1}^n\sum_{\tau\in \sg_n\atop\tau(n)=r}(-q)^{-\inv(\tau)}M_{\tau(1),1}\cdots M_{\tau(n-1),n-1}\;M_{rn}=\sum_{r=1}^n(-q)^{r-n}\det_q(M_{\backslash r\backslash n})M_{rn}, \notag
\end{multline}
where we used the fact that if $\tau(n)=r$ then $\big(\tau(1),\ldots,\tau(n-1)\big)$ is a permutation of $(1,\ldots,r-1,r+1,\ldots,n)$ of length $\inv(\tau)-(n-r)$. The expansion~\eqref{qd5R1} with respect to the first column can be shown in the same way.

To prove the formula~\eqref{ColExpFml} we consider the permutation $\sigma=\left(\begin{smallmatrix}
1,\ldots,s-1,s\phantom{+1},s+1,\ldots,n-1,n \\
1,\ldots,s-1,s+1,s+2,\ldots,n\phantom{-1},s \end{smallmatrix}\right)$.
Taking into account $(M^\sigma)_{\backslash r\backslash n}=M_{\backslash r\backslash s}$, $\inv(\sigma)=n-s$ and the formula~\eqref{qd2R} one gets
\begin{multline}
\det_q(M)=(-q)^{\inv(\sigma)}\det_q(M^\sigma)= (-q)^{\inv(\sigma)}\sum_{r=1}^{n} (-q)^{r-n} \det_q\big((M^\sigma)_{\backslash r\backslash n}\big) M_{rs}= \\
 =\sum_{r=1}^{n} (-q)^{r-s} \det_q\big(M_{\backslash r\backslash s}\big)M_{rs}.
\end{multline}
The first expansion~\eqref{ColExpFml} can be proved similarly. \qed

{\sf Proof of Property \sl\ref{qd6}.} Let $\psi_i$ be $q$-Grassmann variables commuting with $M$, $\psi^M_i=\sum_{k=1}^n \psi_k M_{ki}$ and $I_1=(i_1,\ldots,i_m)$, $I_2=(i_{m+1},\ldots,i_n)$ such that $1\le i_l\le n$. Using~\eqref{lemDetPsi1} we obtain
\begin{multline} \label{jrt_p3}
 \psi^M_{i_1}\cdots\psi^M_{i_n}= (\psi^M_{i_1}\cdots\psi^M_{i_m})(\psi^M_{i_{m+1}}\cdots\psi^M_{i_n})= \\
 =\sum_{K_1=(k_1<\ldots<k_m)\atop K_2=(k_{m+1}<\ldots<k_n)} \psi_{k_1}\cdots\psi_{k_m}\psi_{k_{m+1}}\cdots\psi_{k_n} \det_q(M_{K_1I_1})\det_q(M_{K_2I_2}).
\end{multline}
where $1\le k_1<\ldots<k_m\le n$, $1\le k_{m+1}<\ldots<k_n\le n$.
Note that $\psi_{k_1}\cdots\psi_{k_m}\psi_{k_{m+1}}\cdots\psi_{k_n}$ does not vanish only if $K_2=\backslash K_1$. In this case $(k_1,\ldots,k_m,k_{m+1},\ldots,k_n)$ is a permutation of $(1,\ldots,m,m+1,\ldots,n)$ of length $\sum_{l=1}^m(k_l-l)$:
\begin{align} \label{jrt_p31}
 \psi^M_{i_1}\cdots\psi^M_{i_n} =\sum_{K=(k_1<\ldots<k_m)}(-q)^{-\sum_{l=1}^m(k_l-l)} \psi_{1}\cdots\psi_{n} \det_q(M_{KI_1})\det_q(M_{\backslash K I_2}).
\end{align}
Substituting $i_l=l$ for $l=1,\ldots,n$ we arrive at~\eqref{LaplFml1}. If the matrix $M$ is $q$-Manin the variables $\psi^M_l$ are $q$-Grassmann. Substituting the formula~\eqref{psiviapsi} in the left hand side of~\eqref{jrt_p31} we obtain ~\eqref{LaplFml2}. Formula~\eqref{LaplFml3} can be proved in the same way. \qed

{\sf Proof of Property \sl\ref{qd7}.}  We need to provide a counterexample, which is readily found in the $2\times 2$ case. Indeed, from formula~(\ref{qdetmat2}) and the commutation relations of Remark (\ref{remcrossComRel11}) we have, for the \qMM\ 
$
M
=\left(\begin{array}{cc}
 a & b \\
c & d
\end{array}\right),
$ 
\[
a(\det_q M)=a(ad-q^{-1}cb)=ada-qabc\neq\text{ (in general) } ada-q^{-1}cba=(\det_q M)a.
\]
\qed
\begin{Rem} \label{tildeMqDet} \normalfont 
It is easy to see that if  $M$ is a \qMM, then $\wt M_{ij}=M_{n-i+1,n-j+1}$ is
a $q^{-1}$-\MM{} and  one can also prove, that
$\det_q M= \det_{q^{-1}} \wt M$.
For example
\begin{eqnarray}
det_q
\left(\begin{array}{cc}
 a & b   \\
 c & d
\end{array}\right)
=ad-q^{-1}cb=da-qbc=
det_{q^{-1}}
\left(\begin{array}{cc}
 d & c   \\
 b & a
\end{array}\right).
\end{eqnarray}
\end{Rem}

\noindent{\bf Proof.}The fact  that if $M$ is a $q$-Manin matrix then $\wt M$ is a $q^{-1}$-Manin matrix  easy follows from the formulae~\eqref{crossComRel110}, \eqref{crossComRel11}(see also Remark~\ref{tiledM}). Let us apply Proposition~\ref{lemDetPsi} to the matrix $\wt M$ and $q^{-1}$-Grassmann variables $\wt\psi_i=\psi_{n+1-i}$. The formula~\eqref{lemDetPsi1} gives us $\wt\psi^{\wt M}_1\cdots\wt\psi^{\wt M}_n=\wt\psi_1\cdots\wt\psi_n\det_{q^{-1}}(\wt M)$, where $\wt\psi^{\wt M}_j=\wt\psi_i\wt M_{ij}$. Note that $\wt\psi^{\wt M}_j=\sum_i\psi_i M_{i,n+1-j}=\psi^M_{n+1-j}$, where $\psi^M_j=\sum_i\psi_iM_{ij}$ satisfy~\eqref{lemDetPsi1} with $\psi_i$ and $\det_q(M)$. This yields
\begin{multline}
\wt\psi^{\wt M}_1\cdots\wt\psi^{\wt M}_n=\psi^M_n\cdots\psi^M_1=(-q)^{-\frac{n(n-1)}{2}}\psi^M_1\cdots\psi^M_n= \\
=(-q)^{-\frac{n(n-1)}{2}}\det_q(M)\psi_1\cdots\psi_n=\det_q(M)\psi_n\cdots\psi_1=\det_q(M)\wt\psi_1\cdots\wt\psi_n.
\end{multline}
So we obtain $\det_{q^{-1}}(\wt M)=\det_q(M)$. \qed

\subsection{The $q$-Characteristic polynomial  \label{qCharSS}}

Let us discuss now some subtleties  in the definition of the characteristic polynomial for $q-$Manin matrices. 
In the classical case the  characteristic polynomial of $M$ was defined as $\det(\lambda -M)$, and the same definition (provided the column-determinant is used) hold for the $q=1$ case of the "ordinary" Manin matrices discussed in \cite{CFR08}. 
However, the fact that, for $\lambda \in {\mathbb{C}}$ the matrix 
$(\lambda-M)$ is {\em not} {a \qMM} clearly signals that the na\"ive generalization -- (i.e., the $q$-determinant of $(\lambda-M)$) -- cannot be
the correct one. To get a "good" definition of $q$-characteristic polynomial, we must rely on another property-definition of the  characteristic 
polynomial of a matrix.
\begin{Def}\label{Def4} The $q$-characteristic polynomial for a matrix $M$
is defined as follows:
\begin{align}
\Char_q(\lambda, M)&=\sum_{k=0}^n (-1)^k\lambda^{n-k} \sum_{I=(i_1<i_2<\ldots<i_k)}
\det_q( M_{II})= \\
& =  \lambda^n - \lambda^{n-1} \tr M+ \ldots. + (-1)^n \det_q M,
\end{align}
that is, it is the weighted sum of principal $q$-minors.
Here $(M_{II})_{ab}=M_{i_a i_b}$ are principal submatrices of $M$ of the size equal to the cardinality of $I$.
\end{Def}

Clearly enough, for $q=1$ one gets the usual definition: $\Char_{q=1}(\lambda,M)=\det_{q=1}(\lambda-M)$.

An equivalent definition could be given as follows. Let $\Lambda^k[\psi_1,\ldots,\psi_n]$ be the subspace of the $q$-Grassmann algebra $\mathbb C[\psi_1,\ldots,\psi_n]$ consisting of elements of degree $k$. For a $n\times n$ matrix $M$ over $\mathfrak R$ we define the (left) action of $M$ on $\mathfrak R\otimes\mathbb C[\psi_1,\ldots,\psi_n]$ as $M(r_0\psi_{i_1}\psi_{i_2}\cdots\psi_{i_k})=M(\psi_{i_1}) M(\psi_{i_2})\cdots M(\psi_{i_k})r_0$, where $r_0\in\mathfrak R$ and $M(\psi_j)=\sum_i \psi_i M_{ij}$. The subspace $\mathfrak R\otimes\Lambda^k[\psi_1,\ldots,\psi_n]$ is invariant under this action and then the $q$-characteristic polynomial can be written in the following form
\begin{align}
\Char_q(\lambda, M)=\sum_{k=0}^n
\lambda^{n-k} (-1)^k \tr_{\Lambda^k[\psi_1,\ldots,\psi_n]} M.
\end{align}

We will show in this paper that the coefficients of the 
$q$-characteristic polynomial satisfy many important properties. Indeed, 
they enter $q$-analogues of  the Newton, Cayley-Hamilton, and MacMahon-Wronski
identities and have important applications in quantum integrable systems.

\subsection{A $q$-generalization of the Cramer formula and  quasi-de\-ter\-mi\-nants}

Here we consider a square $n\times n$ \qMM{} and present some relation between its $q$-determinant and the inverse matrix. In the commutative case they are reduced to the Cramer's formula. We formulate these relations in terms of Gelfand-Retakh-Wilson's quasi-determinants and apply the theory of quasi-determinants to the Gauss decomposition.

\subsubsection{Left adjoint matrix}
\label{sec341}

First we introduce the left adjoint matrix in terms of $(n-1)\times(n-1)$ $q$-minors. The main aim of this subsubsection is to express the (left) inverse matrix through the adjoint matrix, or equivalently, through these $q$-minors.

\begin{Prop} \label{lem.3c}
Let $M$ be a \qMM{} and $M^{adj}$ be the matrix with the entries
\begin{align}
 M^{adj}_{sr}=(-q)^{r-s}\det_q M_{\backslash r \backslash s}, \label{M_left_ad_entr}
\end{align}
where $M_{\backslash r \backslash s}$ is the $(n-1)\times(n-1)$ submatrix
of $M$ defined by deleting $r$-th row and $s$-th column.
Then
\begin{align} \label{CramerFormula1}
M^{adj}M=\det_q M\cdot 1_{n\times n},
\end{align}
 i.e. $M^{adj}$ is a left adjoint matrix for the matrix $M$.
\end{Prop}

\noindent{\bf Proof.} Using the
formula \eqref{ColExpFml} for the column expansion
and vanishing of \qDet for coincident columns
(formula \eqref{lemCoinCol})
one gets
\begin{align}
 \sum_{r=1}^n M^{adj}_{sr} M_{rk}=\sum_{r=1}^n(-q)^{r-s} \det_q ( M_{\backslash r \backslash s}) M_{rk}
 =\det_q\wt M=\delta_{sk}\det_q M,
\end{align}
where $\wt M$ is the matrix obtained from $M$ by the replacement the $s$-th column with  the $k$-th one. \qed

\begin{Ex} \normalfont
 In the case $n=2$ the formula~\eqref{CramerFormula1} is deduced as follows (in the notations~\eqref{Mabcd})
\begin{eqnarray}
\left(\begin{array}{ccc}
 d & -qb \\
 -q^{-1}c  & a
\end{array}\right)
\left(\begin{array}{ccc}
 a & b \\
 c & d
\end{array}\right)
=
\left(\begin{array}{ccc}
 da-qbc & db -qbd \\
 -q^{-1}ca +ac & -q^{-1}cb+ad
\end{array}\right)
=
\left(\begin{array}{ccc}
 \det_q M & 0 \\
 0 & \det_q M
\end{array}\right). \notag
\end{eqnarray}
\end{Ex}

\begin{Cor} \label{lem.3d}
Let $M$ be a \qMM{} and suppose that its $q$-determinant $\det_q M$ is invertible from the left, so that there exists an element (a left inverse of the $q$-determinant) $(\det_q M)^{-1}$ such that $(\det_q M)^{-1}\det_q M=1$. Then the product $M^{-1}=(\det_q M)^{-1} M^{adj}$ is a left inverse of $M$, that is $M^{-1}M=1$. The matrix $M^{-1}$ consists of the following entries
\begin{align}
 M^{-1}_{sr}=(\det_q M)^{-1} M^{adj}_{sr}=(-q)^{r-s}(\det_q M)^{-1} \det_q M_{\backslash r \backslash s}. \label{M_left_inv_entr}
\end{align}
In particular, the existence of a left inverse of the $q$-determinant $\det_q(M)$ of a \qMM{} $M$ implies the existence of a left inverse of $M$.
\end{Cor}
Let us remark that the left invertibility of a \qMM{} $M$ does not imply the left 
invertibility of $\det_qM$. On the other hand, neither left nor right invertibility of $\det_q M$ nor 
both of them implies the right invertibility of $M$. 
Let us also remind that the left (right) invertibility of an element of some non-commutative algebra does not guarantee that the left (right) inverse is unique. In Corollary~\ref{lem.3d} we claim that there exists at least one inverse $M^{-1}=(\det_q M)^{-1} M^{adj}$.

\subsubsection{Relation with the quasi-determinants 
 \label{SectQuasiDet} }

We will herewith recall a few constructions from the theory of
quasi-determinants  of I.\,Gelfand and V.\,Retakh and
discuss their counterparts in the case of \MMs. It is fair to say
that the general theoretical set-up of quasi-determinants
developed by Gelfand, Retakh and collaborators (see, e.g.,
\cite{GGRW02}, \cite{GR91}, \cite{GR97})  can be
briefly presented as follows: {\em the basic  facts of linear algebra can
be reformulated with the only use of an inverse matrix}. Thus it can
be extended to the non-commutative set-up and can be  applied, for
example, to some questions considered here. We must stress the
difference between our set-up and the much more general one of \cite{GGRW02}: 
we herewith consider
{\em a special class} of  matrices with non-commutative entries (the
\qMMs), and for this class we can extend theorems of linear algebra
basically in the same form as in the commutative case, (in
particular, as we have seen,  there exists a well-defined notion of
determinant). On the other hand, in \cite{GGRW02} {\em generic
matrices} are considered; thus there is no natural notion of the
determinant, and facts of linear algebra are not exactly presents in
the same form as in the commutative case.

Let us recall (\cite{GGRW02} definition 1.2.2) that the $(p,q)$-th quasi-determinant $|A|_{pq}$ of an invertible matrix $A$
is defined as $|A|_{pq}= (A^{-1}_{qp})^{-1}$, i.e. as the inverse to
the $(q,p)$-element of the matrix inverse to $A$. It is also denoted
by:
\begin{eqnarray}
|A|_{pq}
= \left| \begin{array}{cccccc}
A_{11} & A_{12} &
\ldots &  \ldots    & \ldots & A_{1n} \cr
 ...  & ...  & ... & ... & ...  & ... \cr
...        & ... & ... & \bo {A_{pq} } &... &... \cr
 ...  & ...  & ... & ... & ...  & ... \cr
\end{array} \right| , 
\end{eqnarray}

Let $M$ be a $n\times n$ \qMM{}. Let $M$ and $\det_q M$ be two-sided invertible. Then their inverse are unique and the entries of $M^{-1}$ are expressed by the formula~\eqref{M_left_inv_entr}. Moreover, in the case when $\det_q M$ is two-sided invertible the left (right) invertibility of the entry $M^{-1}_{sr}$ is equivalent to the left (right) invertibility of the determinant $\det_q(M_{\backslash r\backslash s})$. Thus the following lemma holds.

\begin{Lem} If the matrix $M$ and the elements $\det_q M$ and an entry $M^{-1}_{sr}$ are two-sided invertible (for some $1\le s,r\le n$) then the $(r,s)$-th quasi-determinant $|M|_{rs}=(M^{-1}_{sr})^{-1}$ is expressed by the formula
\begin{eqnarray}
|M|_{rs}= (-q)^{s-r}
 \big(\det_q(M_{\backslash r\backslash s})\big)^{-1} \det_q(M). \label{quasidetDef}
\end{eqnarray}
\end{Lem}

Using the formula~\eqref{quasidetDef} we can generalize the notion of quasi-determinant to the case when $\det_q M$ is not two-sided invertible.

\begin{Def}
 Let the determinant $\det_q(M_{\backslash r\backslash s})$ be two-sided invertible. Then the element $|M|_{rs}$ defined by the formula~\eqref{quasidetDef} is called the $(r,s)$-th quasi-determinant $|M|_{rs}$ for a \qMM{} $M$.
\end{Def}

The following lemma is often useful in applications of
quasi-determinants to determinants. It holds thanks to
the Cramer  rule for \qMMs.
\begin{Prop} (c.f.~\cite{GR91,GR92}). \label{DetQuasDetLem}
Assume $M$ is a \qMM, then
\begin{eqnarray} && \det_q \left(
\begin{array}{cccc}
M_{11}& M_{12} &\ldots & M_{1n}\cr
M_{21}& M_{22} &\ldots & M_{2n}\cr
\vdots &\vdots &\ddots & \vdots\cr
M_{n1}& M_{n2} &\ldots & M_{nn}\cr
\end{array}
\right)
=  \label{DetQuasDetLemR} \\
 = &&
M_{nn}
\left|\begin{array}{cccc}
\bo{M_{n-1\,n-1}} & M_{n-1\,n}\cr
M_{n\,n-1} & M_{nn}\cr
\end{array}
\right|
\ \ldots \
\left|\begin{array}{cccc}
\bo{M_{22}}& \ldots & M_{2n}\cr
\vdots &\ddots & \vdots\cr
M_{n2}& \ldots & M_{nn}\cr
\end{array}
\right|
\left|\begin{array}{cccc}
\bo{M_{11}}& M_{12} &\ldots & M_{1n}\cr
M_{21}& M_{22} &\ldots & M_{2n}\cr
\vdots &\vdots &\ddots & \vdots\cr
M_{n1}& M_{n2} &\ldots & M_{nn}\cr
\end{array}
\right|
= \notag \\ = &&
M_{11}
\left|\begin{array}{cccc}
{M_{1\,1}} & M_{1\,2}\cr
M_{2\,1} & \bo{M_{22}}\cr
\end{array}
\right|
\ \ldots \
\left|\begin{array}{cccc}
M_{11}& \ldots & M_{1 n-1}\cr
\vdots &\ddots & \vdots\cr
M_{n-11}& \ldots &\bo{ M_{n-1n-1} }\cr
\end{array}
\right|
\left|\begin{array}{cccc}
{M_{11}}& M_{12} &\ldots & M_{1n}\cr
M_{21}& M_{22} &\ldots & M_{2n}\cr
\vdots &\vdots &\ddots & \vdots\cr
M_{n1}& M_{n2} &\ldots & \bo{ M_{nn}}\cr
\end{array}
\right| \ , \notag
\end{eqnarray}
if the corresponding quasi-determinants are defined.
\end{Prop}

\noindent{\bf Proof.} By the formula~\eqref{quasidetDef} we obtain
\begin{align}
&\left|\begin{array}{cccc}
\bo{M_{kk}}& \ldots & M_{kn}\cr
\vdots &\ddots & \vdots\cr
M_{nk}& \ldots & M_{nn}\cr
\end{array}
\right|=
\left[\det_q\left(\begin{array}{cccc}
M_{k+1,k+1}& \ldots & M_{k+1,n}\cr
\vdots &\ddots & \vdots\cr
M_{n,k+1}& \ldots & M_{nn}\cr
\end{array}
\right)\right]^{-1}
\det_q\left(\begin{array}{cccc}
M_{kk}& \ldots & M_{kn}\cr
\vdots &\ddots & \vdots\cr
M_{nk}& \ldots & M_{nn}\cr
\end{array}
\right). \notag
\end{align}
Substituting it to~\eqref{DetQuasDetLemR} we obtain the identity for the first formula. The second formula is proved in the same way. \qed

\subsubsection{Gauss decomposition and \qDet}

Here we show that the $q$-determinant of a \qMM{} can be expressed via the
diagonal part of the Gauss decomposition exactly in the same way as
in the commutative case.

\begin{Prop} \label{Gauss-pr}
 Let $M$ be a \qMM. If it is factorized into Gauss form
\begin{eqnarray} \label{GaussM11}
M=\left(\begin{array}{ccc}
 1 &{}&x_{\alpha \beta}\\
{}&\ddots&{}\\
0&{}&1
\end{array}\right)
\left(\begin{array}{ccc}
y_1&{}&0\\
{}&\ddots&{}\\
0&{}&y_n
\end{array}\right)
\left(\begin{array}{ccc}
1&{}&0\\
{}&\ddots&{}\\
z_{\beta \alpha}&{}&1
\end{array}\right)
\end{eqnarray}
where $y_1,\ldots,y_n$ are two-sided invertible then
\begin{eqnarray} \label{GaussMyn1}
det_q M=  y_n \ldots y_1.
\end{eqnarray}
Analogously if
\begin{eqnarray}
M=\left(\begin{array}{ccc}
 1 &{}&0\\
{}&\ddots&{}\\
x_{\alpha \beta}' &{}&1
\end{array}\right)
\left(\begin{array}{ccc}
y_1'&{}&0\\
{}&\ddots&{}\\
0&{}&y_n'
\end{array}\right)
\left(\begin{array}{ccc}
1&{}&z_{\beta \alpha}'\\
{}&\ddots&{}\\
0&{}&1
\end{array}\right),
\end{eqnarray}
with two-sided invertible $y'_1,\ldots,y'_n$ then
\begin{eqnarray} \label{GaussMy1n}
det_q M= y_1' \ldots y_n'.
\end{eqnarray}
\end{Prop}

\begin{Ex} \normalfont In the case $n=2$ we have
\begin{eqnarray}
&& M= \left(\begin{array}{ccc}
a  & b \\
c  & d \\
\end{array}\right)
=
\left(\begin{array}{ccc}
1  & bd^{-1} \\
0  & 1 \\
\end{array}\right)
\left(\begin{array}{ccc}
a-bd^{-1}c  & 0 \\
0  & d \\
\end{array}\right)
\left(\begin{array}{ccc}
1  & 0 \\
 d^{-1} c   & 1 \\
\end{array}\right),  \\
&& det_q(M)=ad-q^{-1}cb= da-q bc=d(a-bd^{-1} c)
,\\
&& M= \left(\begin{array}{ccc}
a  & b \\
c  & d \\
\end{array}\right)
=
\left(\begin{array}{ccc}
1  & 0 \\
ca^{-1}  & 1 \\
\end{array}\right)
\left(\begin{array}{ccc}
a  & 0 \\
0  & d- c a^{-1} b \\
\end{array}\right)
\left(\begin{array}{ccc}
1  & a^{-1} b \\
 0   & 1 \\
\end{array}\right),  \\
&& det_q(M)=ad-q^{-1}cb=a(d-ca^{-1} b).
\end{eqnarray}
\end{Ex}

Remark that in the second of these example one only needs the
first column $q$-commutativity relation.\\
\noindent{\bf Proof.} First we note that the decomposition~\eqref{GaussM11} for the following submatrices holds with the same $y_i$:
\begin{eqnarray} \notag
M_{(k)}=\left(\begin{array}{cccc}
M_{kk}& \ldots & M_{kn}\cr
\vdots &\ddots & \vdots\cr
M_{nk}& \ldots & M_{nn}\cr
\end{array}
\right)
=\left(\begin{array}{ccc}
 1 &{}&x_{\alpha \beta}\\
{}&\ddots&{}\\
0&{}&1
\end{array}\right)
\left(\begin{array}{ccc}
y_k&{}&0\\
{}&\ddots&{}\\
0&{}&y_n
\end{array}\right)
\left(\begin{array}{ccc}
1&{}&0\\
{}&\ddots&{}\\
z_{\beta \alpha}&{}&1
\end{array}\right).
\end{eqnarray}
Then note that a left (right) triangular matrix with non-commutative entries and with unity on the diagonal is two-sided invertible and its inverse is a left (right) triangular matrix with non-commutative entries and with unity on the diagonal. Hence the inverse of $M_{(k)}$ has the form
\begin{eqnarray}
M_{(k)}^{-1}=\left(\begin{array}{ccc}
1&{}&0\\
{}&\ddots&{}\\
\tilde z_{\beta \alpha}&{}&1
\end{array}\right)
\left(\begin{array}{ccc}
y_k^{-1}&{}&0\\
{}&\ddots&{}\\
0&{}&y_n^{-1}
\end{array}\right)
\left(\begin{array}{ccc}
 1 &{}&\tilde x_{\alpha \beta}\\
{}&\ddots&{}\\
0&{}&1
\end{array}\right).
\end{eqnarray}
Taking into account $(M_{(k)}^{-1})_{11}=y_k^{-1}$ we obtain the corresponding quasi-determinant
\begin{align} \label{ykMkkqd}
&\left|\begin{array}{cccc}
\bo{M_{kk}}& \ldots & M_{kn}\cr
\vdots &\ddots & \vdots\cr
M_{nk}& \ldots & M_{nn}\cr
\end{array}
\right|=y_k.
\end{align}
Thus Proposition~\ref{DetQuasDetLem} and $M_{nn}=y_n$ imply~\eqref{GaussMyn1}. The formula~\eqref{GaussMy1n} is proved similarly. \qed

\begin{Rem} \normalfont
 The Gauss components $x_{\alpha\beta}$, $y_k$, $z_{\beta\alpha}$, $x_{\alpha\beta}$, $y_k$, $z_{\beta\alpha}$ are unique and can be expressed in terms of quasideterminants such as~\eqref{ykMkkqd} (see~\cite{GR92,GGRW02} for details). Thus they can be expressed in terms of $q$-minors of $M$.
\end{Rem}

\section{$q$-Minors of a $q$-Manin matrix and its inverse }
\label{sec4}
In this section we discuss several formulae holding (in a manner substantially similar to that of the classical case)
for \qMMs. In particular, in the first Subsection we prove the $q$-analogue of the Jacobi ratio theorem (which expresses  minors
of inverse matrix in terms of the minors of the original matrix). 
Building on this, in the second Subsection, we recover 
($q$-analogues of several statements of linear algebra, such as the so-called Dodgson (or Lagrange -- Desnanot-- Jacobi-- Lewis Carroll) condensation formula, the Schur decomposition 
formula, the Sylvester identities, and (although limiting ourselves to the simplest non-trivial case), the Pl\"ucker 
relations. We also show that the inverse of a \qMM is a $q^{-1}$-Manin matrix.

\subsection{Jacobi ratio theorem}
This subsection is devoted to the proof of the $q$ analogue of the Jacobi ratio theorem. We deem useful to preliminarily prove 
a few propositions and lemmas, to be used in the main proof.. 
They also might be of independent interest. 
In particular, we obtain a formula for the right inverse of a two-sided invertible $q$-Manin matrix.

\subsubsection[Preliminary propositions]{Preliminary propositions. Right inverse for the $q$-determinant of a \qMM}




\begin{Lem} \label{ComInvLem} Assume that an element $\psi$ commute with all the entries of
a two-sided invertible matrix $M$; then
$\psi$ commute with all the entries of its inverse $M^{-1}$.
\end{Lem}

\noindent{\bf Proof.} Multiplying the relation $\psi M_{ij}=  M_{ij} \psi$ by $M^{-1}_{jk}$ on the right,  $M^{-1}_{li}$ on the left and taking summation over $j$ and $i$ we obtain $M^{-1}_{lk}\psi = \psi M^{-1}_{lk}$. \qed

\begin{Lem} \label{lemDetPsiR1mrev}
Let $\psi_1,\ldots,\psi_n$ be $q$-Grassmann variables (or at least non-commuting elements satisfying $\psi_i\psi_j=-q\psi_j\psi_i$ for $i<j$). Then
\begin{align} \label{lemDetPsiR1mrevR1}
 \psi_{k_{\tau(m)}}\cdots\psi_{k_{\tau(1)}}=(-q)^{\inv(\tau)}\psi_{k_m}\cdots\psi_{k_1},
\end{align}
where $k_1<\ldots<k_m$ and $\tau\in\sg_m$.
\end{Lem}

\noindent{\bf Proof.} Let $\sigma$ be the longest element of $\sg_m$, i.e. $\sigma(i)=m+1-i$. Then the use of the formula~\eqref{lemDetPsiR1_m} and $\inv(\tau\sigma)=\inv(\sigma)-\inv(\tau)$ yields~\eqref{lemDetPsiR1mrevR1}. 

\begin{Prop} \label{PropGrasDetInv}
Consider a two-sided  invertible \qMM{} $M$. Let $\psi_i$ be $q$-Grassmann variables commuting with entries of $M$ and let $\psi^M_j= \sum_i \psi_i M_{ij}$. Then for a general multi-index $J=(j_1,j_2,\ldots,j_m)$ we have
\begin{align}
 \psi_{j_m}\cdots\psi_{j_1}&=
\sum_{L=(l_1<\ldots<l_m)}
 \psi^M_{l_m}\cdots\psi^M_{l_1}\det_{q^{-1}} (M^{-1}_{LJ}), \label{qinv-minG}
\end{align}
where $M^{-1}_{L J}=(M^{-1})_{L J}$ is the corresponding $m\times m$ submatrix of $M^{-1}$.
\end{Prop}

%
%
%

\noindent{\bf Proof.} 
 Using $\psi_j=\sum_l\psi_l^M M^{-1}_{lj}$ and taking into account that $\psi_j$ commute with the entries of $M^{-1}$ (Lemma~\ref{ComInvLem}) we obtain, 
 at first
\begin{multline}
\psi_{j_m}\cdots\psi_{j_1}= (\sum_{l_m} \psi_{l_m}^M M^{-1}_{l_mj_m}) \psi_{j_{m-1}} \cdots\psi_{j_1}
 \stackrel{Lemma:~[M^{-1}_{ij},\psi_k]=0}{=}
\sum_{l_m} \psi_{l_m}^M  \psi_{j_{m-1}} \cdots\psi_{j_1} M^{-1}_{l_mj_m}.
\end{multline}
Then, the second step and further iterations yield
\begin{multline}
\sum_{l_m} \psi_{l_m}^M  (\sum_{l_{m-1}} \psi_{l_{m-1}}^M  M^{-1}_{l_{m-1}j_{m-1}}) \psi_{l_{m-2}}
\cdots\psi_{j_1} M^{-1}_{l_mj_m}
\stackrel{Lemma:~[M^{-1}_{ij},\psi_k]=0}{=} \\ =
\sum_{l_m}\sum_{l_{m-1}} \psi_{l_1}^M   \psi_{l_2}^M  \psi_{j_{m-3}}
\cdots\psi_{j_1}  M^{-1}_{l_{m-1}j_{m-1}} M^{-1}_{l_mj_m}.
\end{multline}
We obtain making this transformation $m$ times :
\begin{align}
\sum_{l_1,\ldots,l_m}
 \psi_{l_m}^M \cdots \psi_{l_1}^M
 M^{-1}_{l_1j_1}\cdots M^{-1}_{l_mj_m}.
\end{align}
But using the Manin's property (Proposition \pref{proposizione}) one has that the elements $\psi_{l_1}^M,\ldots,\psi_{l_m}^M$ are $q$-Grassmann variables, so
we obtain finally \begin{align}
\sum_{l_1<l_2<\ldots<l_m}
 \sum_{\tau \in \sg_m}
 \psi_{l_{\tau(m)}}^M \cdots \psi_{l_{\tau(1)}}^M
 M^{-1}_{l_{\tau(1)} j_1}\cdots M^{-1}_{l_{\tau(m)} j_m}.
\end{align}
Using Lemma~\ref{lemDetPsiR1mrev} we can rewrite it as
\begin{multline*}
\sum_{\tau \in \sg_m }
\sum_{l_1<l_2<\ldots<l_m }
(-q^{-1})^{-\inv(\tau)} \psi_{l_m}^M \cdots \psi_{l_1}^M
 M^{-1}_{\tau(l_1) j_1}\cdots M^{-1}_{\tau(l_m)j_m} = \\
=
\sum_{L=(l_1<\ldots<l_m)}
 \psi^M_{l_m}\cdots\psi^M_{l_1}\det_{q^{-1}} (M^{-1}_{LJ}) .
\end{multline*}
Thus the formula~\eqref{qinv-minG} is proved. %
%
 \qed

\begin{Rem} \normalfont
Let us note that two-sided invertibility of $M$ is crucial since we need it to use Lemma~\ref{ComInvLem}.
\end{Rem}

\begin{Prop} \label{PropQCBinet}
Let $M$ be a two-sided invertible \qMM{} and consider a multi-index $K=(k_1< k_2<\ldots<k_m)$.
 Then for a general multiindex of a cardinality $m$, say, $J=(j_1 , j_2, \ldots ,j_m)$ we have 
\begin{align}
 \label{CHq222G}
\sum_{L=(l_1<l_2<\ldots<l_m)} \det_q(M_{KL})
\det_{q^{-1}}(M^{-1}_{LJ})=\sum_{\tau\in \sg_m}(-q)^{\inv(\tau)}\delta^{j_1}_{k_{\tau(1)}}\cdots\delta^{j_m}_{k_{\tau(m)}},
\end{align}
where $M^{-1}_{L J}=(M^{-1})_{L J}$ is a submatrix of $M^{-1}$.
\end{Prop}

%
%
%

\noindent{\bf Proof.} Multiplying the relation~\eqref{lemDetPsi2} by  $(-q)^{-m(m-1)/2}$ and re-ordering $q$-Grassmann variables we obtain
\begin{align}
 \psi^M_{l_m}\cdots \psi^M_{l_1}
 & =\sum_{K=(k_1<\ldots<k_m)}\psi_{k_m}\cdots\psi_{k_1}
 \det_q (M_{KL}). \label{q-minrev}
\end{align}
Substituting the formula~\eqref{q-minrev} in~\eqref{qinv-minG} we obtain
\begin{align}
 \psi_{j_m}\cdots\psi_{j_1}=\sum_{L=(l_1<l_2<\ldots<l_m)\atop K=(k_1<\ldots<k_m)} \psi_{k_m}\cdots\psi_{k_1} \det_q(M_{KL})
 \det_{q^{-1}}(M^{-1}_{LJ}). \label{qminrev1}
\end{align}
Comparing the coefficients at $\psi_{k_m}\cdots\psi_{k_1}$ with $k_1<\ldots<k_m$ in sides of~\eqref{qminrev1} we see that the expression in the left hand side of~\eqref{CHq222G} does not vanish only if $K$ is a permutation of $J$, that is if $(j_1,\ldots,j_m)=(k_{\tau(1)},\ldots,k_{\tau(m)})$ for some $\tau\in\sg_m$.
The use of  to Lemma~\ref{lemDetPsiR1mrev}  yields in this case
\begin{align}
 (-q)^{\inv(\tau)}\psi_{k_m}\cdots\psi_{k_1}=\sum_{L=(l_1<l_2<\ldots<l_m)\atop K=(k_1<\ldots<k_m)} \psi_{k_m}\cdots\psi_{k_1} \det_q(M_{KL})
 \det_{q^{-1}}(M^{-1}_{LJ}).
\end{align}
Hence the expression in the left hand side of~\eqref{CHq222G} in this case is equal to $(-q)^{\inv(\tau)}$. This gives us the formula~\eqref{CHq222G}. \qed

In the $q=1$ commutative case Proposition~\ref{PropQCBinet} is a corollary of the Cauchy-Binet formula~\eqref{CBfml} applied to the identity $M M^{-1} = 1$. In spite of the fact that the entries of $M$ and $M^{-1}$ do not commute in general the formula~\eqref{CBfml} works also for \qMMs with $\det_{q^{-1}}(M^{-1})$ and $\det_q(M)$.

Considering the formula~\eqref{CHq222G} in the case $K=J=(1,2,\ldots,n)$ we obtain the following corollary.

\begin{Cor} \label{cor431}
 If $M$ is two-sided invertible \qMM{} then the $q^{-1}$-determinant of its inverse $M^{-1}$ is right inverse of $\det_q(M)$:
\begin{align} \label{Invfml11}
\det_q(M)\det_{q^{-1}}(M^{-1})=1.
\end{align}
In particular, if the $q$-determinant of a two-sided invertible $q$-Manin matrix is right invertible.
\end{Cor}

\begin{Rem} \normalfont Let $M$ be an $n\times n$ $q$-Manin matrix and $\wt M$ be a $n\times n$ matrix satisfying $\sum_{j=1}^n\wt M_{jk}M_{ij}=\sum_{j=1}^nM_{jk}\wt M_{ij}=\delta_{ik}$. Then one can derive the following formulae in similar way:
\begin{align}
\sum_{L=(l_1<l_2<\ldots<l_m)} \det_{q^{-1}}(\wt M_{LJ}) \det_q(M_{KL})
&=\sum_{\tau\in \sg_m}(-q)^{\inv(\tau)}\delta^{j_1}_{k_{\tau(1)}}\cdots\delta^{j_m}_{k_{\tau(m)}}, \label{CHq222Gwt} \\
\det_{q^{-1}}\wt M\det_q M&=1. \label{Invfml11wt}
\end{align}
\end{Rem}

\begin{Rem}\normalfont  Is it also true that $\det_{q^{-1}} (M^{-1}) \det_q M =1$ for a two-sided invertible \qMM{} $M$? In other words, does the two-sided invertibility of a $q$-Manin matrix implies the left invertibility of its $q$-determinant? We don't know an answer to this question.
\end{Rem}

\subsubsection{Jacobi ratio theorem for \qMMs}

Here we will formulate and prove a  $q$-analogue of the Jacobi ratio theorem.

Let $J=(j_1,\ldots,j_m)$ where $1\le j_l\le n$. We define the symbol $\epsilon^{q^{-1}}_J=\epsilon^{q^{-1}}_{j_1,\ldots,j_m}$ by the formulae $\varepsilon^{q^{-1}}_{\ldots,j,\ldots,j,\ldots}=0$ and
\begin{align}
 &\varepsilon^{q^{-1}}_{k_{\tau(1)},\ldots,k_{\tau(m)}}=(-q)^{\inv(\tau)} &  \text{for $1\le k_1<\ldots<k_m\le n$ and $\tau\in \sg_m$}.
\end{align}
It generalizes the symbol defined in ~\eqref{vareps} for $q^{-1}$ and $m\le n$.

\begin{Th} \label{qjrtThQMM}
Let $M$ be a two-sided  invertible \qMM{}. Consider multi-indices $I=(i_1< i_2<\ldots<i_m)$ and $J=(j_1,\ldots,j_m)$. Then the minors of $M^{-1}$ can be expressed through the complementary minors of $M$ by the formula
\begin{align}
 \label{qJRTfml0}
  \det_q (M) \det_{q^{-1}} (M^{-1}_{IJ})
 =(-q)^{\sum\limits_{l=1}^m (j_l-i_l)}
 \varepsilon^{q^{-1}}_J\det_q(M_{\backslash J\backslash I}),
 \end{align}
where $M^{-1}_{L J}=(M^{-1})_{L J}$ and, as usual,  we denote by $\backslash I$ the multi-index obtained by deleting $i_1,\ldots,i_m$ from the sequence $(1,2,3,\ldots,n)$.~\footnote{If $J=(\ldots,j,\ldots,j,\ldots)$ the expression $\det_q(M_{\backslash J\backslash I})$ is not defined, but $\varepsilon^{q^{-1}}_J=0$ in this case and then we can formally consider $\varepsilon^{q^{-1}}_J\det_q(M_{\backslash J\backslash I})=0$. In the case $I=J=\emptyset$ it is natural to define $\det_q\big(M_{IJ}\big)=1$ for an arbitrary matrix $M$.}
In particular
\begin{align}
  \det_q (M) \det_{q^{-1}} (M^{-1}_{IJ})
 &=(-q)^{\sum\limits_{l=1}^m (j_l-i_l)}
 \det_q (M_{\backslash J\backslash I}), && \text{if $j_1<j_2<\ldots<j_m$,} \label{qJRTfml1} \\
  \det_q (M) \det_{q^{-1}} (M^{-1}_{IJ}) &=0 && \text{if $j_a=j_b$ for some $a\ne b$.} \label{qJRTfml2}
 \end{align}
\end{Th}

Let us remark that in the case $I=J=(1,2,\ldots,n)$ the Jacobi ratio formula~\eqref{qJRTfml1} coincides with~\eqref{Invfml11}.


\noindent{\bf Proof.}
Let us rewrite the $q$-Laplace identity~\eqref{LaplFml2} in the form
\begin{align}
\label{qLap111}
\varepsilon^q_{ (\backslash I  \oplus L) }
\det_q(M) = \sum_{K=(k_1<k_2<\ldots<k_{m})}
(-q)^{-\sum\limits_{l=1}^{n-m} (k_{m+l}-l)}
\det_q(M_{\backslash K \backslash I }) \det_q(M_{K L }),
\end{align}
where $L=(l_1<l_2<\ldots<l_m)$ is an arbitrary multi-index, $\varepsilon^q$ is the usual $q$-epsilon-symbol defined by~\eqref{vareps};
multi-index $(\backslash I \oplus L)=(i_{m+1},i_{m+2},\ldots,i_{n},l_1,l_2,\ldots,l_m)$ where $i_{m+1},i_{m+2},\ldots,i_{n}$ are integers such that $\backslash I=(i_{m+1},i_{m+2},\ldots,i_{n})$. Notice that
\begin{align}
 \sum\limits_{l=1}^{n-m} (k_{m+l}-l)=\sum\limits_{l=m+1}^nk_l-\sum\limits_{l=1}^{n-m}l
=\sum\limits_{l=1}^n l-\sum\limits_{l=1}^m k_l-\sum\limits_{l=1}^{n-m}l
=\sum\limits_{l=n-m+1}^n l-\sum\limits_{l=1}^m k_l.
\label{qLap111sum}
\end{align}

Let us multiply the identity \eqref{qLap111} by
$ \det_{q^{-1}} (M^{-1}_{L J})$ from the right, take a summation over $L=(l_1<l_2<\ldots<l_m)$ and transform the right hand side using~\eqref{qLap111sum} and Proposition~\ref{PropQCBinet}:
\begin{multline} \label{qLap222}
\sum_{L=(l_1<l_2<\ldots<l_m)}
\varepsilon^q_{\backslash I \oplus L }
\det_{q}(M) \det_{q^{-1}} (M^{-1}_{L J}) =\\
 = \sum_{L=(l_1<l_2<\ldots<l_m)\atop K=(k_1<k_2<\ldots<k_m)}
(-q)^{\sum\limits_{l=1}^m k_l-\sum\limits_{l=n-m+1}^nl}
\det_{q}(M_{\backslash K \backslash I}) \det_{q}(M_{K L})
\det_{q^{-1}} (M^{-1}_{L J})= \\
=\sum_{K=(k_1<k_2<\ldots<k_m)\atop \tau\in \sg_m}
(-q)^{\sum\limits_{l=1}^m k_l-\sum\limits_{l=n-m+1}^nl}
\det_{q}(M_{\backslash K \backslash I})(-q)^{\inv(\tau)}\delta^{j_1}_{k_{\tau(1)}}\cdots\delta^{j_m}_{k_{\tau(m)}}.
\end{multline}
The summation in the right hand side gives us $(-q)^{\sum\limits_{l=1}^m j_l-\sum\limits_{l=n-m+1}^nl} \varepsilon^{q^{-1}}_J\det_{q}(M_{\backslash J \backslash I})$.
 At the left hand side we have $\varepsilon^q_{\backslash I \oplus L } = 0$, unless $L=I$.
So the formula~\eqref{qLap222} transforms to
\begin{align}
\label{qLap222222}
\varepsilon^q_{\backslash I \oplus I }
\det_{q}(M) \det_{q^{-1}} (M^{-1}_{I J})
 =
(-q)^{\sum\limits_{l=1}^m j_l-\sum\limits_{l=n-m+1}^nl} \varepsilon^{q^{-1}}_J\det_{q}(M_{\backslash J \backslash I}),
\end{align}
Recalling that  (see \eqref{EpsIcI2} page \pageref{EpsIcI2})
\begin{align}
\label{VE222}
\varepsilon^q_{\backslash I \oplus I } =
(-q)^{\sum\limits_{l=1}^m i_l-\sum\limits_{l=n-m+1}^nl},
\end{align}
we come to the formula~\eqref{qJRTfml0}. \qed

\begin{Rem} \normalfont
The factor in the theorem  can be rewritten via the complementary indices:
\begin{align}
 \sum\limits_{l=1}^m (j_l-i_l)=
 -\sum\limits_{l=m+1}^n (j_l-i_l)
\end{align}
since $\sum_{l=1}^n i_l=\sum_{l=1}^n j_l=1+2+\ldots+n$.
\end{Rem}


\subsection{Corollaries }

Theorem~\ref{qjrtThQMM} has a number of important consequences. In particular, the following one is the key to prove the fact that the inverses of some $q$-Manin matrices are $q^{-1}$-Manin matrices.

\subsubsection{Lagrange-Desnanot-Jacobi-Lewis Caroll formula 
}

Let us discuss the special case of the $q$-Jacobi ratio theorem~\ref{qjrtThQMM} for indices of length two. It is interesting by many reasons. Besides the mentioned application (see the next subsubsection) it plays role in wide range of 
questions~\cite{DI}.

Recall that the Lagrange-Desnanot-Jacobi-Lewis Carroll formula reads:
\begin{align}
 \det(M_{\backslash j\backslash i})\det(M_{\backslash l\backslash k})-\det(M_{\backslash j\backslash k})\det(M_{\backslash l\backslash i})&=\det M \det(M_{\backslash (jl)\backslash (ik)}), \label{LDJLC}
\end{align}
where $M$ is a matrix over $\mathbb C$, $j<l$ and $i<k$. According to  (\cite{Bres99} page 111), 
Lagrange found this identity for $n=3$, Desnanot proved it for $n\le
6$, Jacobi proved his general theorem,
C. L. Dodgson -- better known as Lewis Carroll
-- used it to derive an algorithm for calculating determinants that
required only $2\times 2$ determinants
(\href{http://www.jstor.org/pss/112607}{"Dodgson's condensation"}
method \cite{D1866}).

Considering the case $m=2$ of Theorem~\ref{qjrtThQMM} we obtain the following non-commutative $q$-generalization of this formula.

\begin{Prop} ~\label{MainLemInv}
Consider a two-sided  invertible \qMM{} $M$
(i.e. $\exists M^{-1}:~~M^{-1} M=MM^{-1}=1$).
Then for $1\le i_1<i_2\le n$,~~~$1\le j\le n$ and $1\le j_1<j_2\le n$ 
\begin{align} \label{MainLemInv-fml1}
 \det_q M \,\big(M^{-1}_{i_1j}M^{-1}_{i_2j}
- q M^{-1}_{i_2j}M^{-1}_{i_1j}\big)&=0, \\
\det_q M\,\big(M^{-1}_{i_1j_1}M^{-1}_{i_2j_2}
- q M^{-1}_{i_2j_1}M^{-1}_{i_1j_2}\big)
&= (-q)^{j_1+j_2-i_1-i_2} \label{MainLemInv-fml2}
\det_q \big(M_{\backslash (j_1j_2)\backslash (i_1i_2)}\big), \\
\det_q M\,\big(M^{-1}_{i_1j_2}M^{-1}_{i_2j_1}
- q M^{-1}_{i_2j_2}M^{-1}_{i_1j_1}\big)
&= (-q)^{j_1+j_2-i_1-i_2+1}
\det_q \big(M_{\backslash (j_1j_2)\backslash (i_1i_2)}\big), \label{MainLemInv-fml3}
\end{align}
where $M^{-1}_{ij}$ are entries of the inverse matrix $M^{-1}$ and $M_{\backslash (j_1j_2)\backslash (i_1i_2)}$ as usual
is the $(n-2)\times (n-2)$ matrix obtained from $M$
deleting rows $j_1,j_2$ and columns $i_1,i_2$.
\end{Prop}

\begin{Rem} \label{RemLDJLC}\normalfont
In the commutative case, i.e. $q=1$ and $[M_{ij},M_{ij}]=0$, one has $\det(M_{\backslash j\backslash i})=(-1)^{i+j}\det M \cdot M^{-1}_{ij}$ and then both formulae~\eqref{MainLemInv-fml2} and \eqref{MainLemInv-fml3} imply the formula~\eqref{LDJLC} of Lagrange, Desnanot, Jacobi and Lewis Caroll.
\end{Rem} 

\begin{Rem} \label{RemLDJLC2} \normalfont
Using the formula for the left adjoint matrix $M^{adj}=\det_q M \cdot M^{-1}$ one can rewrite the relations~\eqref{MainLemInv-fml1}--\eqref{MainLemInv-fml3} in the form
\begin{align}
 M^{adj}_{i_1j}M^{-1}_{i_2j}-qM^{adj}_{i_2j}M^{-1}_{i_1j}&=0, \label{remLDJLC1} \\
 M^{adj}_{i_1j_1}M^{-1}_{i_2l_2}-qM^{adj}_{i_2j_1}M^{-1}_{i_1l_2}&=(-q)^{j_1+l_2-i_1-i_2} \det_q\big(M_{\backslash (j_1l_2)\backslash (i_1i_2)}\big), \label{remLDJLC2} \\
 M^{adj}_{i_1l_2}M^{-1}_{i_2j_1}-qM^{adj}_{i_2l_2}M^{-1}_{i_1j_1}&=(-q)^{j_1+l_2-i_1-i_2+1} \det_q\big(M_{\backslash(j_1l_2)\backslash(i_1i_2)}\big), \label{remLDJLC3}
\end{align}
where $1\le i_1<i_2\le n$,~~~$1\le j\le n$ and $1\le j_1<j_2\le n$.
The formulae~\eqref{remLDJLC1}--\eqref{remLDJLC3} are still valid if the matrix $M$ is invertible {\itshape only} from the right, where $M^{-1}$ is some right inverse: $MM^{-1}=1$, 
(see Appendix~\ref{AppPropLDJLC}).
\end{Rem}



\subsubsection{The inverse of a \qMM{} is a $q^{-1}$-Manin matrix}

In~\cite{CF07},\cite{CFR08} it was proved that the inverse of a two-sided invertible Manin matrix (the case $q=1$) 
is again a Manin matrix. There, it  was also shown that this fact has a series of applications. Here we present its analogue for \qMMs.

Consider a right invertible \qMM{} $M$ with left invertible $q$-determinant $\det_q M$. The left invertibility of $\det_q M$ implies the left invertibility of $M$ (Corollary~\ref{lem.3d}). Then according to the Corollary~\ref{cor431} the two-sided invertibility of $M$ implies the right invertibility of $\det_q M$. In this way that the matrix $M$ and its $q$-determinant $\det_q M$ are two-sided invertible; i.e., there exist a matrix $M^{-1}$ and an element $(\det_q M)^{-1}$ such that
\begin{gather}
 M^{-1}M=M^{-1}M=1_{n\times n} \\
 (\det_q M)^{-1}\det_q(M)=\det_q(M)(\det_q M)^{-1}=1, \label{corLDJLC_}
\end{gather}
Notice also that due to the corollaries~\ref{lem.3d}, \ref{cor431} they are related as
\begin{align}
 M^{-1}&=(\det_q M)^{-1} M^{adj}, \\
 (\det_q M)^{-1}&=\det_{q^{-1}}(M^{-1}). \label{corLDJLC__}
\end{align}

\begin{Th} \label{corLDJLC}
If a \qMM{} $M$ is invertible from the right and its $q$-determinant $\det_q M$ is invertible from the left then the inverse matrix $M^{-1}$ is a $q^{-1}$-Manin matrix.
\end{Th}

\noindent{\bf Proof.} Due to the left invertibility of the $q$-determinant of $M$ one can multiplying the equations~\eqref{MainLemInv-fml1}, \eqref{MainLemInv-fml2}, \eqref{MainLemInv-fml3} by $(\det_q M)^{-1}$ from the left. In the first equation we obtain
\begin{align}
 M^{-1}_{i_1j}M^{-1}_{i_2j}-qM^{-1}_{i_2j}M^{-1}_{i_1j}&=0,  \label{MainLemInv-fml11}
\end{align}
where $i_1<i_2$. Comparing the last two equations one gets
\begin{align}
M^{-1}_{i_1j_1}M^{-1}_{i_2j_2}
- q M^{-1}_{i_2j_1}M^{-1}_{i_1j_2}
=-q^{-1}M^{-1}_{i_1j_2}M^{-1}_{i_2j_1}
+ M^{-1}_{i_2j_2}M^{-1}_{i_1j_1}.
\end{align}
where $i_1<i_2$, $j_1<j_2$. So that we derive the commutation relations for the entries of a $q^{-1}$-Manin matrix (see the formulae~\eqref{crossComRel110}, \eqref{crossComRel11}). \qed

\subsubsection{Schur complements}
Here we apply Theorem~\ref{qjrtThQMM} to the Schur complements.  
Let $M$ be an arbitrary $n\times n$ matrix over a non-commutative ring, which we consider partitioned  in  four  blocks:
\begin{align} \label{M_ABCD}
 M=\begin{pmatrix}
     A & B\\
         C & D
     \end{pmatrix},
\end{align}
where $A$, $B$, $C$ and $D$ are $m\times m$, $m\times(n-m)$, $(n-m)\times m$ and $(n-m)\times(n-m)$ matrices respectively.

 Suppose first that the matrix $M$ is left-invertible and let $M^{-1}$ be its left inverse, to be analogously partitioned into
\begin{align} \label{Minv_ABCD}
 M^{-1}=\begin{pmatrix}
     \wt A & \wt B\\
     \wt C & \wt D
     \end{pmatrix}.
\end{align}
Suppose also that the matrices $A$ and $D$ are right-invertible. Then $\wt A(A-BD^{-1}C)\hm=1$ and $\wt D(D-CA^{-1}B)\hm=1$, where $A^{-1}$ and $D^{-1}$  and let $A^{-1}$ and $D^{-1}$ are right inverses of $A$ and $D$ respectively.

If the matrix $M$ is right-invertible with right inverse of the form~\eqref{Minv_ABCD} and the matrices $A$ and $D$ are left-invertible then $(A-BD^{-1}C)\wt A\hm=1$ and $(D-CA^{-1}B)\wt D\hm=1$, where $A^{-1}$ and $D^{-1}$ are left inverse of $A$ and $D$ respectively. 
So finally we have the following lemma.

\begin{Lem} \label{lemsch_comps}
If $M$, $A$ and $D$ are two-sided invertible then the matrices $\wt A$ and $\wt D$ defined from~\eqref{Minv_ABCD} are also two-sided invertible and
\begin{align}
\wt A^{-1}&=A-BD^{-1}C, & \wt D^{-1}&=D-CA^{-1}B. \label{sch_comps}
\end{align}
\end{Lem}

\begin{Def} \normalfont
 The matrices $A-BD^{-1}C$ and $D-CA^{-1}B$ are called {\it Schur complements} of the blocks $D$ and $A$ respectively.
\end{Def}

In terms of multi-indices we have
\begin{align} \label{ADMij}
 A&=M_{II}, & D&=M_{\backslash I\backslash I}, & \wt A&=M^{-1}_{II}, & \wt D&=M^{-1}_{\backslash I\backslash I},
\end{align}
where $I=(1,\ldots,m)$, $\backslash I=(m+1,\ldots,n)$.
Using~\eqref{sch_comps} and \eqref{ADMij} we can apply Theorem~\ref{qjrtThQMM} \big(the case $I=J=(1,\ldots,m)$\big) to the Schur complements.

\begin{Prop} \label{propdet_Sch_comps}
 Let $M$ be a \qMM{} and let $A$ and $D$ be its submatrices defined by the formula~\eqref{M_ABCD} (which are also $q$-Manin). Suppose $M$, $A$ $D$ are right-invertible and their $q$-determinants $\det_q M$, $\det_q A$ and $\det_q D$ are left-invertible (hence they all are two-sided invertible). Then the Schur complements~\eqref{sch_comps} are \qMMs{} as well and their $q$-determinants satisfy the relations
\begin{align} \label{det_Sch_comps}
  \det_q M&=\det_q D\det_q(A-BD^{-1}C), & \det_q M&=\det_q A \det_q(D-CA^{-1}B).
\end{align}
\end{Prop}
\noindent{\bf Proof.} In order to prove that Schur complements are $q$-Manin we show that the matrices $\wt A$ and $\wt D$ defined by the formula~\eqref{Minv_ABCD} satisfy the condition of Theorem~\ref{corLDJLC} (replacing $q$ by  $q^{-1}$). First let us note that the matrix $M$ satisfies the conditions of this theorem and hence the inverse matrix $M^{-1}$ is $q^{-1}$-Manin. The matrices $\wt A$ and $\wt D$ are also $q^{-1}$-Manin as submatrices of $M^{-1}$. Lemma~\ref{lemsch_comps} implies that the matrices $\wt A$ and $\wt D$ are two-sided invertible. Then, writing down the formula~\eqref{qJRTfml1} for $I=J=(1,\ldots,m)$ and $I=J=(m+1,\ldots,n)$ and taking into account~\eqref{ADMij} we obtain
\begin{align} \label{jrt_AD}
 \det_{q^{-1}}\wt A&=(\det_q M)^{-1}\det_q D, & \det_{q^{-1}}\wt D&=(\det_q M)^{-1}\det_q A.
\end{align}
This implies that the $q^{-1}$-determinants $\det_{q^{-1}}\wt A$ and $\det_{q^{-1}}\wt D$ are also two-sided invertible. Applying Theorem~\ref{corLDJLC} to $q^{-1}$-Manin matrices $\wt A$ and $\wt D$ we conclude that their inverse $\wt A^{-1}$ and $\wt D^{-1}$, i.e. the Schur complements~\eqref{sch_comps}, are \qMMs. The formulae~\eqref{det_Sch_comps} follow from~\eqref{jrt_AD} if one takes into account the equality
\begin{align}
 \big(\det_{q^{-1}}\wt A\big)^{-1}=\det_q(\wt A^{-1})=\det_q(A-BD^{-1}C), \\
 \big(\det_{q^{-1}}\wt D\big)^{-1}=\det_q(\wt D^{-1})=\det_q(D-CA^{-1}B),
\end{align}
\big(see the formula~\eqref{corLDJLC__}\big). \qed

\subsubsection{Sylvester's theorem} 
Sylvester's identity is a classical determinantal identity (see, e.g.
\cite{GGRW02} theorem 1.5.3). 
The Sylvester identities for the non commutative case were discussed, in the framework of the theory of quasi-determinants,
in \cite{KL94}, and, 
using combinatorial methods for the $q$-analogues of matrices of the form
$1+M$ in 
\cite{Konvalinka07-2}. Here
we show that for \qMMs\  the identity easily follows from the Schur Theorem \ref{propdet_Sch_comps}.
For the sake of simplicity, we shall suppose that all matrices are two-sided invertible as well as their determinants.

Let us first recall the commutative case: 
 Let $A$ be a matrix $(a_{ij})_{m \times m}$; take $n < i,j \leq m$; denote:
\bea
A_0 = \begin{pmatrix} a_{11} & a_{12} & \cdots & a_{1n} \\
               a_{21} & a_{22} & \cdots & a_{2n} \\
               \vdots & \vdots & \ddots & \vdots \\
               a_{n1} & a_{n2} & \cdots & a_{nn} \end{pmatrix}, \quad
 a_{i*} = \begin{pmatrix} a_{i1} & a_{i2} & \cdots & a_{in} \end{pmatrix}, \quad a_{*j} = \begin{pmatrix} a_{1j} \\ a_{2j} \\ \vdots \\ a_{nj}
\end{pmatrix}.
\eea Define the $(m-n)\times (m-n)$ matrix $B$ as follows:
 \bea
B_{ij} = {\det}\begin{pmatrix} A_0 & a_{*j}
\\ a_{i*} & a_{ij} \end{pmatrix}, \quad B = (B_{ij})_{n+1 \leq i,j \leq m}.
\eea
 Then
 \bea\label{csyl}
{\det}B = {\det}A \cdot ({\det}A_0)^{m-n-1} . \eea 
{\Th (Sylvester's
identity for \qMMs.)
 Let $M$ be $m\times m$ a \qMM~ with right and left inverse; take $n < i,j \leq m$ and denote:
\bea
M_0 = \begin{pmatrix} M_{11} & M_{12} & \cdots & M_{1n} \\
               M_{21} & M_{22} & \cdots & M_{2n} \\
               \vdots & \vdots & \ddots & \vdots \\
               M_{n1} & M_{n2} & \cdots & M_{nn} \end{pmatrix}, \quad
 M_{i*} = \begin{pmatrix} M_{i1} & M_{i2} & \cdots & M_{in} \end{pmatrix},
\quad M_{*j} = \begin{pmatrix} M_{1j} \\ M_{2j} \\ \vdots \\ M_{nj}
\end{pmatrix}.
\eea Define the $(m-n)\times (m-n)$ matrix $B$ as follows:
 \bea\label{bsylv}
B_{ij} = ({\det}_q(M_0))^{-1} \cdot
{\det}_q\begin{pmatrix} M_0 & M_{*j}
\\ M_{i*} & M_{ij} \end{pmatrix}, \quad B = (B_{ij})_{n+1 \leq i,j \leq m}.
\eea
 Then the matrix $B$ is a \qMM~ and
 \bea\label{ncsyl}
{\det}_qB = ({\det}_qM_0)^{-1} \cdot {\det}_qM.
\eea
}
%

\PRF Once chosen $M_0$, we consider the resulting block
decomposition of  $M$, \bea M= \left(
\begin{array}{cc}
M_0 & M_1 \\
M_2 & M_3
\end{array}\right).
\eea

The key observation is that the matrix $B$ defined by (\ref{bsylv}) equals to the Schur
complement matrix: $M_0-M_2 (M_3)^{-1} M_1$.  To see this, we use Schur complement Theorem  (Theorem \ref{propdet_Sch_comps}) again:
 \bea
B_{ij} = ({\det}_q(M_0))^{-1} \cdot
{\det}_q\begin{pmatrix} M_0 & M_{*j}
\\ M_{i*} & M_{ij} \end{pmatrix}
= ({\det}_q(M_0))^{-1} \cdot ( ({\det}_q(M_0)) ( M_{ij} - M_{i*} M_{0}^{-1} M_{*j}) )= \\
( M_{ij} - M_{i*} M_{0}^{-1} M_{*j}) =(M_0-M_2 (M_3)^{-1} M_1)_{ij}.
\eea In particular, we used the Schur formula
${\det}_q(M)={\det}_q(A){\det}_q(D-CA^{-1}B)$ for blocks $M_0=A$, $ M_{*j} =B$,
$M_{i*}=C$, $M_{ij}=D$, the last being a $1\times 1$ matrix. 
The theorem now follows from the Schur property of the matrix $B$ immediately; indeed, $B=(M_0-M_2
(M_3)^{-1} M_1)$ is a \qMM, since a Schur complement is a \qMM~ by
Theorem  \ref{propdet_Sch_comps}. Then ${\det}_q B = ({\det}_qM_0)^{-1} \cdot {\det}_qM$ follows from the formula \ref{det_Sch_comps} 
for the determinant of Schur complements.
\BX

\subsection{Pl\"ucker relations: an example}
In this subsection we shall briefly address the problem of the existence and of the form of Pl\"ucker relations for $q$-determinants of \qMMs  (giving an example). On general grounds (see \cite{Lauve04}) we know that quasi-Pl\"ucker identities exists within the theory of quasi-determinants of Gel'fand-Retakh-Wilson. 
For \qMMs, these identities acquire a form that (up to the appearance of suitable powers of $q$), is the same as that of the commutative case. We shall only deal with the case of the Pl\"ucker identities in the case of $q-Gr(2,4)$, which is sufficiently enlightening. The analysis of the generic case can be easily performed with some combinatorics.
\begin{Prop}
 \label{Plucker} Consider a $4\times 2$  \qMM\ $A$, and
 let $\pi_{ij}$ be the $q$-minors made from the $i$-th and $j$-th rows of $A$.
Then it holds:
\begin{equation}
\label{Pluck} (\pi_{12} \pi_{34} +  q^{-4}\pi_{34} \pi_{12} )-
(q^{-1}\pi_{13} \pi_{24} +  q^{-3}\pi_{24} \pi_{13} ) +q^{-2}(\pi_{14} \pi_{23} +
\pi_{23} \pi_{14} )=0 \end{equation}
\end{Prop}
\PRF The proof is basically the same as in the commutative case. Consider the
$q$-Grassmann algebra $\CC[\psi_1,....,\psi_4]$, and the  variables
$\tpsi_1,\tpsi_2$ defined as: \bea 
(\tpsi_1, \tpsi_2)
 = (\psi_1,  \psi_2 , \psi_3, \psi_4) \cdot A. 
\eea It is clear that $ \tpsi_1 \wedge \tpsi_2= \sum_{i\neq j} \pi_{ij}
\psi_i \wedge \psi_j$.  By the defining property of \qMMs\ in terms of coaction on $q$-Grassmann variables (Proposition \ref{proposizione}),
$\tpsi_1$ and $\tpsi_2$ are again Grassmann variables. Hence
$ (\tpsi_1\wedge \tpsi_2)^2$ must vanish. If we
write explicitly this vanishing relation, taking into account the $q$-commutations between the $\psi_j$'s, we get the Pl\"{u}cker
relations (\ref{Pluck}).\qed
\\
{\bf Remark.} Formula (\ref{Pluck}) can be compactly written as
\[
\sum_{\sigma\in S_4^\prime} q^{-\inv(\sigma)} \pi_{\sigma(1)\sigma(2)}\pi_{\sigma(3)\sigma(4)}=0,\]
that is the formula proven in (\cite{Lauve04}) for the {\em quantum} matrix algebra, where $S_4^\prime$ is the subset of the group of permutations of four objects satisfying $\sigma(1)<\sigma(2), \sigma(3)<\sigma(4)$.

\section{Tensor approach to \qMMs}\label{MatrSect}

In previous sections we mainly considered the \qMMs{} over a non-commutative algebra $\mathfrak R$ as homomorphisms from the corresponding Grassmann algebra to the tensor product of this algebra with $\mathfrak R$. Indeed, in Sections~\ref{sec3}, \ref{sec4} we have used  formula $\psi^M_j=\sum_{i=1}^n\psi_i M_{ij}$. Another option is to work with them as with linear maps from the vector space $\mathbb C^n$ to $\mathfrak R\otimes\mathbb C^n$, or in other words as with matrices over the non-commutative algebra $\mathfrak R$. In this case we shall interpret the $q$-determinant and $q$-minors in terms of certain higher tensors. Therefore one can call that {\it tensor approach}.

The approach used above is more natural to consider algebraic properties of \qMMs, while the tensor approach is more useful in area of application to the quantum integrable systems. In this section we are going to reformulate some of the notions presented in previous Sections  of the paper as well as to derive some new properties directly applicable to quantum integrable systems described by Lax matrices.


\subsection{Leningrad tensor notations}

Let us remind tensor notations known as Leningrad notations. First we shall identify $ab$ with $a\otimes b$ for $a\in\mathfrak R$ and $b\in\End(\mathbb C^n)^{\otimes N}$, where $N$ is a number of tensor factors.~\footnote{Here we consider $\mathfrak R\otimes\End(\mathbb C^n)^{\otimes N}$ as a right $\mathfrak{R}$-module and a left $\End(\mathbb C^n)^{\otimes N}$-module.}
Let $\{X_\ell\}$ be a basis in the space $\End(\mathbb C^n)$ and $C\in\mathfrak R\otimes\End(\mathbb C^n)^{\otimes N}$ be a tensor over the algebra $\mathfrak R$. Then $C$ can be written in the form
\begin{align}
 C=\sum_{\ell_1,\ldots,\ell_N} C_{\ell_1,\ldots,\ell_N}\cdot \big(X_{\ell_1}\otimes X_{\ell_2}\otimes\cdots\otimes X_{\ell_N}\big),
\end{align}
where $C_{\ell_1,\ldots,\ell_N}\in \mathfrak{R}$ and $X_{\ell_i}\in \End(\mathbb C^n), \> i=1,\ldots, N.$
Introduce the notation $C^{(k_1,\ldots,k_N)}$ (Leningrad notation) for an element of $\mathfrak R\otimes\End(\mathbb C^n)^{\otimes N'}$, where $N'\ge N$. It is defined as
\begin{align}
 C^{(k_1,\ldots,k_N)}&=\sum_{\ell_1,\ldots,\ell_N}C_{\ell_1,\ldots,\ell_N}\cdot \\
   &\cdot\big(1\otimes\cdots\otimes1\otimes X_{\ell_1}\otimes1\otimes\cdots\otimes1\otimes X_{\ell_2}\otimes1\otimes\cdots\cdots\otimes1\otimes X_{\ell_N}\otimes1\otimes\cdots\otimes1\big), \notag
\end{align}
where each $X_{\ell_i}$ is placed in the $k_i$-th tensor factor (obviously $1\le k_i\le N'$). For example let  $N=N'=2$,  and 
$C=\sum_{\ell_1,\ell_2} C_{\ell_1,\ell_2}X_{\ell_1}\otimes X_{\ell_2}$; then
\begin{align}
  C^{(12)}&=\sum_{\ell} C_{\ell_1,\ell_2}\cdot(X_{\ell_1}\otimes1)(1\otimes X_{\ell_2}), &
  C^{(21)}&=\sum_{\ell} C_{\ell_1,\ell_2}\cdot(1\otimes X_{\ell_1})(X_{\ell_2}\otimes1).
\end{align}
Their product reads as
\begin{align}
  C^{(12)}C^{(21)}=\sum_{\ell_1,\ell_2,\ell_1'\ell_2'} C_{\ell_1,\ell_2}\cdot C_{\ell_1',\ell_2'}\cdot(X_{\ell_1}\otimes X_{\ell_2'})\otimes (X_{\ell_1'}\otimes X_{\ell_2}) . 
\end{align}

Let $E_{ij}$ be the standard matrices with entries $(E_{ij})_{kl}=\delta_{ik}\delta_{jl}$. Then the set $\{E_{ij}\mid i,j=1,\ldots,n\}$ is a basis in $\End(\mathbb C^n)$ and each matrix $M\in\mathfrak R\otimes\End(\mathbb C^n)$ is decomposed as $M=\sum_{i,j=1}^n M_{ij}E_{ij}$, where $M_{ij}\in\mathfrak R$ are entries of $M$. Let $\{e_1,\ldots,e_n\}$ be the standard basis in $\mathbb C^n$: $(e_i)^j=\delta^j_i$, so that $E_{ij}e_k=\delta_{kj}e_i$. Then in this basis the action of the matrix $M$ reads $Me_j=\sum_{i=1}^n M_{ij}e_i$. According to the Leningrad notation the tensor $M^{(1)}M^{(2)}\cdots M^{(m)}$ can be written as
\begin{align}
 M^{(1)}M^{(2)}\cdots M^{(m)}=\sum_{i_1,\ldots,i_m\atop j_1,\ldots,j_m}M_{i_1j_1}M_{i_2j_2}\cdots M_{i_mj_m}\cdot\big(E_{i_1j_1}\otimes E_{i_2j_2}\otimes\cdots\otimes E_{i_mj_m}\big).
\end{align}

\subsection{Tensor relations for \qMMs{}}

In this section we present the defining relations for Manin matrices in the tensor notations. We also present important higher order relations following from these defining quadratic relations.

\subsubsection{Pyatov's Lemma}

Let $P^q\in \End(\CC^n\otimes\CC^n)$ be the $q$-permutation operator: $P^q (e_i\otimes e_j)
=  q^{-\s(i-j)} e_j \otimes e_i$. The matrix of $P^q$ (which we denote by the same symbol) can be written as
\begin{align}
 P^q=\sum_{i,j=1}^n q^{\s(i-j)} E_{ij} \otimes E_{ji}. \label{PqDef}
\end{align}
As the usual permutation matrix $P^1$ (when $q=1$) it satisfies $(P^q)^2=1$. Introduce the following tensors called $q$-anti-symmetrizer and $q$-symmetrizer:
\begin{align}
 &A^q=\frac{1-P^q}2, & &S^q=\frac{1+P^q}2. \label{AqSqDef}
\end{align}
Notice that these are orthogonal idempotents: $(A^q)^2= A^q$, $(S^q)^2= S^q$, $A^q S^q\hm=S^q A^q\hm=0$.

The following lemma was suggested to us by P.\,Pyatov.
Let $M$ be a $n\times n$ matrix with entries belonging to a associative algebra $\mathfrak R$ over $\CC$.
\begin{Lem} \label{P-lem} (Pyatov's lemma).
Matrix $M$ is a \qMM{} if and only if any of the following formul\ae\ holds:
\begin{align}
 M^{(1)} M^{(2)} - P^q M^{(1)} M^{(2)} P^q &= P^q M^{(1)} M^{(2)}-  M^{(1)} M^{(2)} P^q , \label{lemPyat1} \\
A^q M^{(1)} M^{(2)} A^q&= A^q M^{(1)} M^{(2)}, \label{lemPyat2} \\
 S^q M^{(1)} M^{(2)}S^q&=  M^{(1)} M^{(2)}S^q, \label{lemPyat3} \\
(1-P^q) M^{(1)} &M^{(2)} (1+P^q)= 0. \label{lemPyat4}
\end{align}
\end{Lem}

\noindent{\bf Proof.} It is easy to see that the matrix equations~\eqref{lemPyat1}--\eqref{lemPyat4} are equivalent to each other. Rewriting, for instance, the equation~\eqref{lemPyat1} by entries one yields exactly the relations~\eqref{MMdefSingleFml}. \qed

Lemma~\ref{P-lem} can be considered as a definition of \qMMs{} in the tensor form. For instance one can construct $q$-Manin matrices from $RLL$-relations (see Section~\ref{secIntSys}) using the relation~\eqref{lemPyat2} in the form given by the  following 

\begin{Lem} \label{lemma2}
 If the matrix $M$ satisfy the equation
 \begin{align}
 A^qM^{(1)}M^{(2)}=T A^q \label{AMM_TA}
 \end{align}
for some matrix $T\in\mathfrak R\otimes\End(\mathbb C^n\otimes\mathbb C^n)$ then $M$ is a $q$-Manin matrix.
\end{Lem}

\noindent{\bf Proof.} Multiplying~\eqref{AMM_TA} by $A^q$ from the right and taking into account the equality $A^qA^q=A^q$ we obtain
 \begin{align}
 A^qM^{(1)}M^{(2)}A^q=T A^q. \label{AMMA_TA}
 \end{align}
Using again the equation~\eqref{AMM_TA} we can substitute $TA^q$ in the right hand side by $A^qM^{(1)}M^{(2)}$. So we derive~\eqref{lemPyat2} and hence $M$ is a $q$-Manin matrix due to Lemma~\ref{P-lem}. \qed

\begin{Lem} \label{P-lemQInv} 
Since $(P^q)^{(-1)}=P^{q^{-1}}$ an $n\times n$ matrix $M$ is a \qMM{} if and only if any of these formul\ae\ holds:
\begin{align}
  M^{(2)} M^{(1)}  - P^{q^{-1}} M^{(2)} M^{(1)}  P^{q^{-1}}
& =  P^{q^{-1}} M^{(2)} M^{(1)}
- M^{(2)} M^{(1)} P^{q^{-1}},\\
A^{q^{-1}} M^{(2)} M^{(1)} A^{q^{-1}} &= A^{q^{-1}} M^{(2)} M^{(1)}, \label{AMMAAMM21}\\
S^{q^{-1}}  M^{(2)} M^{(1)} S^{q^{-1}} &= M^{(2)} M^{(1)} S^{q^{-1}},\\
A^{q^{-1}} M^{(2)} M^{(1)} S^{q^{-1}} &= 0.
\end{align}
\end{Lem}

Since $(P^q)^\top=P^{q^{-1}}$ the transpose matrix $M^\top$ is a  \qMM{}  if and only if any of the formulae holds:
\begin{align}
  M^{(1)} M^{(2)}  - P^{q^{-1}} M^{(1)} M^{(2)}   P^{q^{-1}}
& =   M^{(1)} M^{(2)} P^{q^{-1}}
- P^{q^{-1}} M^{(1)} M^{(2)},\\
A^{q^{-1}} M^{(1)} M^{(2)}  A^{q^{-1}} &=  M^{(1)} M^{(2)} A^{q^{-1}}, \label{AMMAMMAqinv} \\
S^{q^{-1}}  M^{(1)}  M^{(2)} S^{q^{-1}} &=  S^{q^{-1}}  M^{(1)} M^{(2)},\\
 S^{q^{-1}} M^{(1)} M^{(2)} A^{q^{-1}}  &= 0, \\
  M^{(2)} M^{(1)}  - P^{q} M^{(2)} M^{(1)}   P^{q}
& =   M^{(2)} M^{(1)} P^{q}
- P^{q} M^{(2)} M^{(1)},\\
A^{q} M^{(2)} M^{(1)}  A^{q} &=  M^{(2)} M^{(1)} A^{q}, \label{AMMAMMA}\\
S^{q}  M^{(2)}  M^{(1)} S^{q} &=  S^{q}  M^{(2)} M^{(1)},\\
 S^{q} M^{(2)} M^{(1)} A^{q}  &= 0.
\end{align}

\subsubsection{Higher $q$-(anti)-symmetrizers and \qMMs}

Consider the group homomorphism $\pi_q\colon \sg_m\to\Aut\big((\mathbb C^n)^{\otimes m}\big)$ defined on the generators of $\sg_m$ by formula
\begin{align}
 \pi_q(\sigma_k)=(P^q)^{(k,k+1)}=\sum_{i,j=1}^n q^{\s(i-j)}\big(1^{\otimes(k-1)}\otimes E_{ij} \otimes E_{ji}\otimes1^{\otimes(m-k-1)}\big), \label{piqDef}
\end{align}
where $\sigma_k=\sigma_{k,k+1}$ are adjacent permutations. This is a $q$-deformation of the standard representation of the symmetric group $\sg_m$ on the space $(\mathbb C^n)^{\otimes m}$.
The $q$-anti-symmetrizer and $q$-symmetrizer acting in the space $(\mathbb C^n)^{\otimes m}$ reads
\begin{align}
 A^q_m&=\frac1{m!}\sum_{\tau\in \sg_m}(-1)^{\sigma}\pi_q(\tau), & S^q_m&=\frac1{m!}\sum_{\tau\in \sg_m}\pi_q(\tau),
\end{align}
where $(-1)^{\sigma}=(-1)^{\inv(\sigma)}$ is the sign of the permutation $\sigma\in \sg_m$. These formulae generalize the definitions~\eqref{AqSqDef} as $A^q_2=A^q$, $S^q_2=S^q$.

We can also generalize the relations~\eqref{lemPyat2} and \eqref{lemPyat3} for the higher $q$-anti-symmetrizers and $q$-symmetrizers.

\begin{Th} \label{prop_AMMM_AMMMA}
Let $M$ be an $n\times n$ \qMM{}. The tensor $A^q_m M^{(1)}\cdots M^{(m)}$ is invariant under multiplication by the $q$-anti-symmetrizer from the right:
\begin{align}
  A^q_m M^{(1)}\cdots M^{(m)}=A^q_m M^{(1)}\cdots M^{(m)}A^q_m. \label{AMMM_AMMMA}
\end{align}
The tensor $M^{(1)}\cdots M^{(m)}S^q_m$ is invariant under multiplication by the $q$-symmetrizer from the left:
\begin{align}
  M^{(1)}\cdots M^{(m)}S^q_m=S^q_m M^{(1)}\cdots M^{(m)}S^q_m. \label{MMMS_SMMMS}
\end{align}
\end{Th}

\noindent{\bf Proof.} Let us note that the formula~\eqref{AMMM_AMMMA} is equivalent to the following proposition: the tensor $A^q_m M^{(1)}\cdots M^{(m)}$ is anti-invariant under multiplication by any element $\tau\in\sg_m$ in the representation $\pi_q$ from the right: 
\begin{align} \label{AMMM_AMMMtau}
  A^q_m M^{(1)}\cdots M^{(m)}=(-1)^\tau A^q_m M^{(1)}\cdots M^{(m)}\pi_q(\tau).
\end{align}
Indeed, since $(-1)^\tau A^q_m\pi_q(\tau)=A^q_m$ the right hand side of~\eqref{AMMM_AMMMA} is invariant under multiplication by $(-1)^\tau\pi_q(\tau)$ from the right and hence~\eqref{AMMM_AMMMtau} holds. Conversely, summing the formula~\eqref{AMMM_AMMMtau} over $\tau\in\sg_m$ and dividing by $m!$  yields~\eqref{AMMM_AMMMA}.

Since the formula~\eqref{AMMM_AMMMA} have been already proved for $m=2$ (Lemma~\ref{P-lem}) we have the relation $A^qM^{(1)}M^{(2)}(-P^q)=A^qM^{(1)}M^{(2)}$ or equivalently
\begin{align}
 (A^q)^{(k,k+1)}M^{(k)}M^{(k+1)}(-P^q)^{(k,k+1)}=(A^q)^{(k,k+1)}M^{(k)}M^{(k+1)}.
\end{align}
To prove the formula~\eqref{AMMM_AMMMtau} for general $m$ it is sufficient to prove it for the generators $\tau=\sigma_k$. In this case $\pi_q(\tau)=(P^q)^{(k,k+1)}$. Applying the formula $A^q_m=A^q_m(A^q)^{(k,k+1)}$ in the right hand side of~\eqref{AMMM_AMMMtau} and taking into account $[(A^q)^{(k,k+1)},M^{(l)}]=0$, $1\le l<k$, and $[M^{(l)},(P^q)^{(k,k+1)}]=0$, $k+1<l\le m$, we obtain the left hand side of~\eqref{AMMM_AMMMtau}.

   The relation~\eqref{MMMS_SMMMS} is equivalent to the relations
\begin{align}
  M^{(1)}\cdots M^{(m)}S^q_m=\pi_q(\tau) M^{(1)}\cdots M^{(m)}S^q_m, &&&\tau\in\sg_m, \label{MMMS_tauMMMS}
\end{align}
which can be proved in the same way. \qed \\

Analogously one can obtain the relations
\begin{align}
  A^{q^{-1}}_m M^{(m)}\cdots M^{(1)}&=A^{q^{-1}}_m M^{(m)}\cdots M^{(1)}A^{q^{-1}}_m, \label{AMMM_AMMMA21qinv} \\
  M^{(1)}\cdots M^{(m)}S^{q^{-1}}_m&=S^{q^{-1}}_m M^{(1)}\cdots M^{(m)}S^{q^{-1}}_m, \label{MMMS_SMMMS21qinv}
\end{align}
where $M$ is a $q$-Manin matrix. In particular, any $q^{-1}$-Manin matrix $M$ satisfies
\begin{align}
  A^q_m M^{(m)}\cdots M^{(1)}&=A^q_m M^{(m)}\cdots M^{(1)}A^q_m, \label{AMMM_AMMMA21} \\
  M^{(1)}\cdots M^{(m)}S^q_m&=S^q_m M^{(1)}\cdots M^{(m)}S^q_m, \label{MMMS_SMMMS21}
\end{align}

\begin{Rem} \normalfont
Let us remark that for $m\ge3$ the relations~\eqref{AMMM_AMMMA} and \eqref{MMMS_SMMMS} are not equivalent and do not imply that $M$ is a \qMM.
\end{Rem}

\subsection{The $q$-determinant and the $q$-minors as tensor components}

The vectors $e_{j_1,\ldots,j_m}=e_{j_1}\otimes\ldots\otimes e_{j_m}$ form a basis of the space $(\mathbb C^n)^{\otimes m}$. Each tensor (over $\mathbb C$ or $\mathfrak R$) can be decomposed by entries with respect to this basis (and its dual). Here we obtain the relation between components of the tensor~\eqref{AMMM_AMMMA} and $m\times m$ $q$-minors the matrix $M$.

\subsubsection{Action of the $q$-anti-symmetrizer in the basis $\{e_{j_1,\ldots,j_m}\}$}

Let us consider the representation $\sigma\mapsto(-1)^\sigma\pi_q(\sigma)$ of the permutation group $\sg_m$ in terms of the basis $\{e_{j_1,\ldots,j_m}\}$. For the adjacent transpositions $\sigma_k$ its action has the form
\begin{align}
 -\pi_k(\sigma_k)e_{j_1,\ldots,j_k,j_{k+1},\ldots,j_m}=-q^{\s(j_{k+1}-j_k)}e_{j_1,\ldots,j_{k+1},j_k,\ldots,j_m},
\end{align}
(not $(-q)^{\s(j_{k+1}-j_k)}$). Then one can obtain the following formula
\begin{align}
 (-1)^{\tau^{-1}}\pi_q(\tau^{-1})e_{k_{\sigma(1)},\ldots,k_{\sigma(m)}}=(-q)^{\inv(\sigma\tau)-\inv(\sigma)} e_{k_{\sigma\tau(1)},\ldots,k_{\sigma\tau(m)}}, \label{pi_q_tau_e}
\end{align}
where $k_1<\ldots<k_m$, $\sigma,\tau\in \sg_m$. It is proved in the same way as~\eqref{lemDetPsiR1}. The action of $A^q_m$ on the basis elements $e_{k_{\sigma(1)},\ldots,k_{\sigma(m)}}$ can be obtained by summing~\eqref{pi_q_tau_e} over all $\tau\in \sg_m$, its action on the other elements of the basis is zero:
\begin{align}
 A^q_me_{k_{\sigma(1)},\ldots,k_{\sigma(m)}}=\frac1{m!}\sum_{\tau\in\sg_m}(-q)^{\inv(\sigma\tau)-\inv(\sigma)} e_{k_{\sigma\tau(1)},\ldots,k_{\sigma\tau(m)}},&&&A^q_m e_{\ldots,i,\ldots,i,\ldots}=0. \label{Am_tau_e}
\end{align}

For the multi-indices $I=(i_1,\ldots,i_m)$ and $J=(j_1,\ldots,j_m)$ we will use the notation $M^{i_1,\ldots,i_m}_{j_1,\ldots,j_m}=M_{IJ}$ for an $n\times n$ matrix $M$, that is $M^{i_1,\ldots,i_m}_{j_1,\ldots,j_m}$ is the $m\times m$ matrix with entries $\big(M^{i_1,\ldots,i_m}_{j_1,\ldots,j_m}\big)_{kl}=M_{i_k,j_l}$. This is the matrix consisting of the entries of the matrix $M$ lying in the intersection of $i_k$-th row and $j_l$-th column.
For example, $M^i_j=M_{ij}$ are corresponding entries of $M$.

\subsubsection{Components of the tensor $A^q_mM^{(1)}\cdots M^{(m)}$}

Let us denote by $\{e^{i_1,\ldots,i_n}\}$ the basis of $\big((\mathbb C)^{\otimes m}\big)^*$ dual to the basis $\{e_{j_1,\ldots,j_m}=e_{j_1}\otimes\ldots\otimes e_{j_m}\}$ of $(\mathbb C)^{\otimes m}$, i.e. $\big\langle e^{i_1,\ldots,i_m},e_{j_1,\ldots,j_m}\big\rangle= \delta^{i_1}_{j_1}\cdots\delta^{i_m}_{j_m}$\,.

\begin{Lem} \label{lemAMMM}
Let $M$ be an arbitrary $n\times n$ matrix (over an algebra $\mathfrak R$), then the tensor $A^q_m M^{(1)}\cdots M^{(m)}$ has the following components
\begin{align}
 \big\langle e^{\ldots,i,\ldots,i,\ldots},A^q_m M^{(1)}\cdots M^{(m)}e_{j_1,\ldots,j_m}\big\rangle&=0, \label{AMMM0} \\
 \big\langle e^{k_{\sigma(1)},\ldots,k_{\sigma(m)}},A^q_m M^{(1)}\cdots M^{(m)}e_{j_1,\ldots,j_m}\big\rangle&=\frac1{m!}(-q)^{\inv(\sigma)}\det_q\Big(M^{k_1,\ldots,k_m}_{j_1,\ldots,j_m}\Big), \label{AMMM_det}
\end{align}
where $k_1<\ldots<k_m$, $\sigma\in \sg_m$.
\end{Lem}

\noindent{\bf Proof.} Substituting the action of the matrix $M$ on the basis $\{e_j\}$ in terms of its entries we obtain
\begin{align}
 \big\langle e^{i_1,\ldots,i_m},A^q_m M^{(1)}\cdots M^{(m)}e_{j_1,\ldots,j_m}\big\rangle=
 \sum_{s_1,\ldots,s_m=1}^n M_{s_1,j_1}\cdots M_{s_m,j_m}\big\langle e^{i_1,\ldots,i_m},A^q_m e_{s_1,\ldots,s_m}\big\rangle. \label{eAMMeR1}
\end{align}
If $i_l=i_{l'}=i$ for some $l\ne l'$ then the covector $e^{i_1,\ldots,i_m}=e^{\ldots,i,\ldots,i,\ldots}$ is orthogonal to the image of the operator $A^q_m$ (see the formula~\eqref{Am_tau_e}), hence the corresponding components~\eqref{eAMMeR1} vanish -- yielding the formula~\eqref{AMMM0}. Otherwise the indices $i_1,\ldots,i_m$ are pairwise different and can be represented as $i_l=k_{\sigma(l)}$ for some $\sigma\in \sg_m$ and $k_1,\ldots,k_m$ such that $k_1<\ldots<k_m$. On the other hand the action of $A^q_m$ vanishes on the vectors $e_{s_1,\ldots,s_m}$ if $s_l=s_{l'}$ for some $l\ne l'$. It means that the sum in~\eqref{eAMMeR1} reduced in this case to the sum over the indices $s_1,\ldots,s_m$ that are represented as $s_l=r_{\sigma'(l)}$ for some $\sigma'\in \sg_m$ and $r_1,\ldots,r_m$ such that $r_1<\ldots<r_m$. Substituting $i_l=k_{\sigma(l)}$ and $s_l=r_{\sigma'(l)}$ to~\eqref{eAMMeR1} and calculating the action of $A^q_m$ on the basis vectors by the formula~\eqref{pi_q_tau_e} we obtain
\begin{multline}
 \big\langle e^{k_{\sigma(1)},\ldots,k_{\sigma(m)}},A^q_m M^{(1)}\cdots M^{(m)}e_{j_1,\ldots,j_m}\big\rangle= \\
 =\frac1{m!}\sum_{\sigma',\tau\in \sg_m}\sum_{r_1<\ldots<r_m}M_{r_{\sigma'(1)},j_1}\cdots M_{r_{\sigma'(m)},j_m}\big\langle e^{k_{\sigma(1)},\ldots,k_{\sigma(m)}},(-1)^{\tau^{-1}}\pi_q(\tau^{-1}) e_{r_{\sigma'(1)},\ldots,r_{\sigma'(m)}}\big\rangle= \\
 =\frac1{m!}\sum_{\sigma',\tau\in \sg_m}\sum_{r_1<\ldots<r_m} M_{r_{\sigma'(1)},j_1}\cdots M_{r_{\sigma'(m)},j_m}\times \\
 \times \big\langle e^{k_{\sigma(1)},\ldots,k_{\sigma(m)}},
(-q)^{\inv(\sigma'\tau)-\inv(\sigma')} e_{r_{\sigma'\tau(1)},\ldots,r_{\sigma'\tau(m)}} \big\rangle= \\
 =\frac1{m!}\sum_{\sigma'\in \sg_m}(-q)^{\inv(\sigma)-\inv(\sigma')} M_{k_{\sigma'(1)},j_1}\cdots M_{k_{\sigma'(m)},j_m}=\frac1{m!}(-q)^{\inv(\sigma)}\det_q\Big(M^{k_1,\ldots,k_m}_{j_1,\ldots,j_m}\Big), \notag 
\end{multline}
where we have used formula~\eqref{Def3eq}. \qed

\begin{Cor} \label{CorAMMMm}
 For a $q$-Manin matrix $M$ we obtain
\begin{multline}
A^q_m M^{(1)}\cdots M^{(m)}
=m!\sum_{K=(k_1<\ldots<k_m) \atop J=(j_1<\ldots<j_m)}
\det_q(M_{KJ})\; A^q_m\big(E_{k_1j_1}\otimes\cdots\otimes E_{k_mj_m}\big) A^q_m.
\end{multline}
\end{Cor}

We denote by $\tr_{1,\ldots,m}(\cdot)$ the trace over the space $(\mathbb C^n)^{\otimes m}$. Notice that this trace is the composition of the traces over each tensor factor: $\tr_{1,\ldots,m}(\cdot)=\tr_1\tr_2\ldots\tr_m(\cdot)$\,. The tensor $A^q_m M^{(1)}\cdots M^{(m)}$ can be considered as a matrix (over $\mathfrak R$) acting on the space $(\mathbb C^n)^{\otimes m}$. So we can consider the trace of this matrix:
\begin{align}
 \tr_{1,\ldots,m}\big(A^q_m M^{(1)}\cdots M^{(m)}\big)
 =\sum_{k_1,\ldots,k_m=1}^n\big\langle e^{k_1,\ldots,k_m},A^q_m M^{(1)}\cdots M^{(m)}e_{k_1,\ldots,k_m}\big\rangle. \label{treAMMeR1}
\end{align}

\begin{Cor} \label{CorAMMM1}
For an arbitrary $n\times n$ matrix $M$ the trace~\eqref{treAMMeR1} has the form
\begin{align}
 \tr_{1,\ldots,m}\big(A^q_m M^{(1)}\cdots M^{(m)}\big) 
=\sum_{k_1<\ldots<k_m}\frac1{m!}\sum_{\sigma\in \sg_m}(-q)^{\inv(\sigma)}\det_q \Big(M^{k_1,\ldots,k_m}_{k_1,\ldots,k_m}\Big)^{\sigma}, \label{CorAMMM2eq}
\end{align}
where $\Big(M^{k_1,\ldots,k_m}_{k_1,\ldots,k_m}\Big)^{\sigma}=M^{k_1,\ldots,k_m}_{k_{\sigma(1)},\ldots,k_{\sigma(m)}}$ is the matrix $M^{k_1,\ldots,k_m}_{k_1,\ldots,k_m}$ with columns permuted by $\sigma$ (see the subsection~\ref{sec32}).
\end{Cor}

In some cases we can rewrite the formula~\eqref{CorAMMM2eq} using the Property~\ref{qd2} of Proposition~\ref{lem.3a}.

\begin{Cor} \label{CorAMMM2}
For $n\times n$ matrices $M$ satisfying the relation~\eqref{crossComRel11}(in particular for $q$-Manin matrices ) the trace of $A^q_m M^{(1)}\cdots M^{(m)}$ is the sum of all principal $m\times m$ $q$-minors of this matrix:
\begin{align}
 \tr_{1,\ldots,m}\big(A^q_m M^{(1)}\cdots M^{(m)}\big)
=\sum_{k_1<\ldots<k_m}\det_q\Big(M^{k_1,\ldots,k_m}_{k_1,\ldots,k_m}\Big);
\end{align}
in the notations of the Sections~\ref{sec3}, \ref{sec4} we can rewrite it as
\begin{align}
 \tr_{1,\ldots,m}\big(A^q_m M^{(1)}\cdots M^{(m)}\big)
=\sum_{K=(k_1<\ldots<k_m)}\det_q\big(M_{KK}\big),
\end{align}
In particular, the trace of $A^q_n M^{(1)}\cdots M^{(n)}$ in this case is the $q$-determinant of this matrix:
\begin{align}
 \tr_{1,\ldots,n}\big(A^q_n M^{(1)}\cdots M^{(n)}\big)
=\det_q M.
\end{align}
\end{Cor}

\begin{Cor} \label{AMMdet}
 Any \qMM{} $M$ satisfies the relation
\begin{align}
  A^q_n M^{(1)}\cdots M^{(n)}=A^q_n\det_q M. \label{AMMdet_}
\end{align}
\end{Cor}

\noindent{\bf Proof.} Since the image of the idempotent operator $A^q_n$ is one-dimensional the Proposition~\ref{prop_AMMM_AMMMA} for $m=n$ implies $A^q_n M^{(1)}\cdots M^{(n)}=\mathfrak e_n A^q_n$, where $\mathfrak e_n\in\mathfrak R$. Using the formula $\tr_{1,\ldots,n}(A^q_n)=1$ one gets $\mathfrak e_n=\tr_{1,\ldots,n}\big(A^q_n M^{(1)}\cdots M^{(n)}\big)=\det_q M$. \qed \\


\subsection{$q$-powers of $q$-Manin matrices. Cayley-Hamilton theorem and Newton identities}

We remind the definition of "$q$-corrected" powers $M^{[n]}$ for \qMMs{} and show that the Cayley-Hamilton 
theorem and the Newton identities hold\footnote{Compare for the rational case with \cite{CF07, CFR08}, and for the elliptic case with~\cite{RST}.}

\subsubsection{Cayley-Hamilton theorem}

Remind that (see Definition \ref{Def4}) that
the $q$-characteristic polynomial of an $n\times n$ \qMM{}  $M$ is the following linear combinations of the sums of principal $q$-minors:
\begin{align}
\Char_M(t)=\sum_{m=0}^n(-1)^m \mathfrak e_m t^{n-m}, \label{CharPol}
\end{align}
where $\mathfrak e_0=1$, $\mathfrak e_m=\sum\limits_{i_1<\ldots<i_m}\det_q M^{i_1,\ldots,i_m}_{i_1,\ldots,i_m}$.

\begin{Def}\label{def8}
Let  $M$ be an $n\times n$ \qMM. We call $q$-powers of $M$ the matrices $M^{[m]}$ defined by the formul\ae
\begin{align}
M^{[0]}&=1, & M^{[m]}&=\tr_1\Big(P^q\big(M^{[m-1]}\big)^{(1)}M^{(2)}\Big).
\end{align}
Using the notation $A\qstar B=\tr_1\big(P^qA^{(1)}B^{(2)})\big)=\sum_{ijk}q^{\s(j-i)}A_{ij}B_{jk}E_{ik}$ one can write
\begin{align}
M^{[m]}=M^{[m-1]}\qstar M=(\ldots((M\qstar M)\qstar M)\qstar\ldots\qstar M),
\end{align}
(it is a polynomial of order $m$ in the matrix $M$).
\end{Def}

\begin{Ex} \normalfont The second $q$-power
\begin{math}
M^{[2]}= \tr_{1}\big(P^q M^{(1)} M^{(2)}\big)
\end{math} looks as
\begin{eqnarray}
M^{[2]}=
\left(\begin{array}{cc}
a^2 + qbc & ab + db \\
ac+dc & d^2+ q^{-1} cb
\end{array}\right) =
\left(\begin{array}{cc}
a^2 + qbc & ab + qbd \\
q^{-1} ca + dc & d^2+ q^{-1} cb
\end{array}\right).
\end{eqnarray}
\end{Ex}

\begin{Th} (c.f.~\cite{GIOPS}) \label{propCharPol_qPower}
Any $n\times n$ \qMM{} annihilates its characteristic polynomial by its $q$-powers via right substitution:
\begin{align}
\sum_{m=0}^n(-1)^m \mathfrak e_m M^{[n-m]}=0. \label{CharPol_qPower}
\end{align}
\end{Th}

\begin{Ex} \normalfont The Cayley-Hamilton theorem in the case $n=2$ reads:
\begin{eqnarray}
M^{[2]}- \tr(M) M + det_q (M) 1_{2\times 2} =0.
\end{eqnarray}
\end{Ex}

\noindent{\bf Proof.} Let us start with a  proof of  the following formulae
\begin{align} \label{propCharPol_qPowerR1}
 m\tr_{1,\ldots,m-1}(A^q_m M^{(1)}\cdots M^{(m)})= \\
  = \mathfrak e_{m-1}M-(m-1)\tr_{1,\ldots,m-2}(A^q_{m-1} M^{(1)}\cdots M^{(m-1)})\qstar M.
\end{align}
Using $A^q_m=\frac1m A^q_{m-1}\big(1-(m-1)(P^q)^{(m-1,m)}\big)A^q_{m-1}$ and Corollary~\ref{CorAMMM2}  yields
\begin{multline}\label{propCharPol_qPowerR2}
 m\tr_{1,\ldots,m-1}(A^q_m M^{(1)}\cdots M^{(m)}) = \\
  =\mathfrak e_{m-1}M-(m-1)\tr_{1,\ldots,m-1}\big(A^q_{m-1}(P^q)^{(m-1,m)}A^q_{m-1} M^{(1)}\cdots M^{(m)}\big).
\end{multline}
In the ~\eqref{propCharPol_qPowerR2} we can apply Proposition~\ref{prop_AMMM_AMMMA} for $m-1$ and use periodicity of trace to present it in the form
\begin{multline} \label{propCharPol_qPowerR4}
    -(m-1)\tr_{1,\ldots,m-1}\big((P^q)^{(m-1,m)}A^q_{m-1} M^{(1)}\cdots M^{(m)}\big)= \\
    =-(m-1)\tr_m\big((P^q)^{(m-1,m)}(\tr_{1,\ldots,m-2}A^q_{m-1} M^{(1)}\cdots M^{(m-1)})M^{(m)}\big)= \\
    =-(m-1)(\tr_{1,\ldots,m-2}A^q_{m-1} M^{(1)}\cdots M^{(m-1)})\qstar M.
\end{multline}
So, the formula~\eqref{propCharPol_qPowerR1} is proven. Iterating it one obtains
\begin{align} \label{propCharPol_qPowerR5}
 m(-1)^{m-1}\tr_{1,\ldots,m-1}(A^q_m M^{(1)}\cdots M^{(m)}) =\sum_{k=0}^{m-1}(-1)^k \mathfrak e_k M^{[m-k]}.
\end{align}
The formula~\eqref{CharPol_qPower} is a particular case of~\eqref{propCharPol_qPowerR5} corresponding to $m=n$. Indeed, taking the trace over the first $n-1$ spaces in the formula~\eqref{AMMdet_} and taking into account $\tr_{1,\ldots,n-1}A^q_n=\frac1n$ we have $\tr_{1,\ldots,n-1}A^q_n M^{(1)}\cdots M^{(n)}=\frac1n \mathfrak e_n$. \qed

\subsubsection{Newton theorem}

The classical Newton formulae are relations between elementary symmetric functions 
$\sum_{1\le i_1<\ldots<i_m\le n}\lambda_{i_1}\cdots\lambda_{i_m}$ and sums of powers $\sum_{k=1}^n\lambda_k^m$. 
They allow to inductively express elementary symmetric functions via the sums of powers and vice versa. 
These functions can be presented as sums of principal minors and as traces of powers of the matrix 
$diag(\lambda_1,\ldots,\lambda_n)$. The Newton formulae can be represented in the matrix form, 
which is valid for any square number matrix. In~\cite{CF07,CFR08} they were generalized for ($q=1$)
Manin matrices. Here we present a version of the Newton formulae for \qMMs.

\begin{Th} (c.f.~\cite{GIOPS}) \label{propNewton}
For any $n\times n$ \qMM{} $M$ and any $m\ge0$ we have the relations
\begin{align} \label{Newton_th}
 m \mathfrak e_m=\sum_{k=0}^{m-1}(-1)^{m+k+1} \mathfrak e_k \tr(M^{[m-k]}),
\end{align}
where  $\mathfrak e_m=0$ for $m>n$.
\end{Th}

\noindent{\bf Proof.} Note that the formula~\eqref{propCharPol_qPowerR5} is still valid for $m>n$ if we put $\mathfrak e_k=0$ for $k>n$. Taking the trace over the left space in this formula we obtain~\eqref{Newton_th}. \qed

\begin{Ex} \normalfont
 Let $n=2$, $M=\left(\begin{smallmatrix}a & b\\ c & d\end{smallmatrix}\right)$. Then $\mathfrak e_0=1$, $\mathfrak e_1=\tr(M^{[1]})=\tr M=a+d$, $\mathfrak e_2=\det_q M=ad-q^{-1}cb$, $\mathfrak e_3=\mathfrak e_4=\ldots=0$, $\tr(M^{[2]})=a^2+qbc+d^2+q^{-1}cb$, $\tr(M^{[3]})=a^3+qbca+qabc+qdbc+q^{-1}acb+q^{-1}dcb+d^3+q^{-1}cbd$. One can explicitly check the first tree (non-trivial) relations
\begin{align*}
 \mathfrak e_1&=\tr(M^{[1]}), & 2\mathfrak e_2&=-\tr(M^{[2]})+\mathfrak e_1\tr(M^{[1]}), &  0&=\tr(M^{[3]})-\mathfrak e_1\tr(M^{[2]})+\mathfrak e_2\tr(M^{[1]}).
\end{align*}
Using them one can express $\tr(M^{[3]})$ via $\tr(M^{[1]})$ and $\tr(M^{[2]})$ only:
\begin{align}
 \tr(M^{[3]})=\tr(M^{[1]})\tr(M^{[2]})+\frac12\tr(M^{[2]})\tr(M^{[1]})-\frac12\big(\tr(M^{[1]})\big)^3.
\end{align}
Analogously one can express $\tr(M^{[4]})$, $\tr(M^{[5]})$ etc.
\end{Ex}

Let us now consider the first three relations for general $n$
\begin{align}
 \mathfrak e_1&=\tr(M^{[1]}), \qquad \mathfrak e_2=-\frac12\tr(M^{[2]})+\frac12\mathfrak e_1\tr(M^{[1]}), \\
  \mathfrak e_3&=\frac13\tr(M^{[3]})-\frac13\mathfrak e_1\tr(M^{[2]})+\frac13\mathfrak e_2\tr(M^{[1]}).
\end{align}
Substituting them iteratively we can express $\mathfrak e_1$, $\mathfrak e_2$, $\mathfrak e_3$ via traces of $q$-powers:
\begin{align}
 \mathfrak e_1&=\tr(M^{[1]}), \\ \mathfrak e_2&=-\frac12\tr(M^{[2]})+\frac12\big(\tr(M^{[1]})\big)^2, \\
  \mathfrak e_3&=\frac13\tr(M^{[3]})-\frac13\tr(M^{[1]})\tr(M^{[2]})-\frac16\tr(M^{[2]})\tr(M^{[1]})+\frac16\big(\tr(M^{[1]})\big)^3.
\end{align}
Thus iteratively substituting the first $k$ relations~\eqref{Newton_th} (divided by $m$) we can express $\mathfrak e_k$ via the traces of $q$-powers $\tr(M^{[1]}),\ldots,\tr(M^{[k]})$ for any $k=1,\ldots,n$. Note also that if $m>n$ then the first $m$ relations~\eqref{Newton_th} allow to express $\tr(M^{[m]})$ via $\tr(M^{[1]}),\ldots,\tr(M^{[n]})$ as in the example.

Conversely, one can express the traces of $q$-powers of a \qMM{} via $\mathfrak e_1,\ldots,\mathfrak e_n$. Indeed rewriting~\eqref{Newton_th} as 
\begin{align} \label{Newton_th_trMm}
 \tr(M^{[m]})=(-1)^{m+1}m \mathfrak e_m+\sum_{k=1}^{m-1}(-1)^{k+1} \mathfrak e_k \tr(M^{[m-k]})
\end{align}
and after iterative substitutions one can express $\tr(M^{[m]})$ via $\mathfrak e_1,\ldots,\mathfrak e_m$. For example the first three relations give us
\begin{align}
 \tr(M^{[1]})&=\mathfrak e_1, \\  \tr(M^{[2]})&=-2\mathfrak e_2+\mathfrak e_1^2, \\
 \tr(M^{[3]})&=3\mathfrak e_3-2\mathfrak e_1\mathfrak e_2+\mathfrak e_1^3-\mathfrak e_2\mathfrak e_1.
\end{align}

\subsection{Inverse to a $q$-Manin matrix}

Here we obtain some tensor relations for the inverse to a $q$-Manin matrix and investigate the conditions when these relations holds. First consider a $q$-Manin matrix $M$ over the algebra $\mathfrak R$ invertible from the right with the $q$-determinant invertible from the left (as in Theorem~\ref{corLDJLC}). Then there exists a two-sided inverse $M^{-1}$, which is $q^{-1}$-Manin matrix. Due to Lemma~\ref{P-lemQInv} it satisfies
\begin{align}
A^q (M^{-1})^{(2)} (M^{-1})^{(1)} A^q &= A^q (M^{-1})^{(2)} (M^{-1})^{(1)}. \label{AmiMi}
\end{align}

Suppose that $M^{-1}$ is an $n\times n$ matrix satisfying~\eqref{AmiMi}. According to Lemma~\ref{P-lemQInv} this is a $q^{-1}$-Manin matrix. In particular, it satisfies the relation~\eqref{AMMM_AMMMA21}, which takes the form
\begin{align} \label{AMiMiA}
A^q_m (M^{-1})^{(m)}\cdots (M^{-1})^{(1)} A^q_m &= A^q_m(M^{-1})^{(m)}\cdots (M^{-1})^{(1)}.
\end{align}
Also from Corollary~\ref{CorAMMM2} we obtain
\begin{align}
 \tr_{1,\ldots,m}\Big(A^{q^{-1}}_m (M^{-1})^{(1)}\cdots (M^{-1})^{(m)}\Big)&= 
 \tr_{1,\ldots,m}\Big(A^q_m (M^{-1})^{(m)}\cdots (M^{-1})^{(1)}\Big)= \notag \\
&=\sum_{K=(k_1<\ldots<k_m)}\det_{q^{-1}}\big(M^{-1}_{KK}\big),
\end{align}
In particular one has
\begin{align}
 \tr_{1,\ldots,n}\Big(A^q_n (M^{-1})^{(n)}\cdots (M^{-1})^{(1)}\Big)
=\det_{q^{-1}}(M^{-1}).
\end{align}
Then taking the trace in the both hand sides of the relation
\begin{align}
A^q_n (M^{-1})^{(n)}\cdots (M^{-1})^{(1)} &=\wt{\mathfrak e_n} A^q_n, 
\end{align}
where $\wt{\mathfrak e_n}\in\mathfrak R.$

Using ~\eqref{AMiMiA} we obtain
\begin{align}
A^q_n (M^{-1})^{(n)}\cdots (M^{-1})^{(1)}=A^q_n \det_{q^{-1}}(M^{-1}). \label{AMiMidetMi}
\end{align}

\begin{Prop}
 Let $M$ be $n\times n$ right-invertible matrix over the algebra $\mathfrak R$. Assume that there exists a matrix $B\in\mathfrak R\otimes\End(\mathbb C^n\otimes\mathbb C^n)$ such that
\begin{align}
 BA^qM^{(1)}M^{(2)}=A^q. \label{BAMMA}
\end{align}
Then the right inverse $M^{-1}$ is a $q^{-1}$-Manin matrix and satisfies consequently the relations~\eqref{AmiMi}--\eqref{AMiMidetMi}.

If in addition $M$ is a $q$-Manin matrix then its determinant is invertible from the right:
\begin{align}
 \det_q M\det_{q^{-1}}(M^{-1})=1, \label{detMtetMi}.
\end{align}
Moreover if $M$ is invertible from the left in this case then its determinant is invertible from the left:
\begin{align}
 \det_{q^{-1}}(M^{-1})\det_q M=1. \label{detMitetM}
\end{align}
\end{Prop}

\noindent{\bf Proof.}
Multiplying~\eqref{BAMMA} by $(M^{-1})^{(2)} (M^{-1})^{(1)}$ from the right we obtain $$BA^q=A^q(M^{-1})^{(2)} (M^{-1})^{(1)},$$ which implies~\eqref{AmiMi} (c.f. Lemma~\ref{lemma2}). Hence $M^{-1}$ is a $q^{-1}$-Manin matrix satisfying also~\eqref{AMiMiA}--\eqref{AMiMidetMi}. If $M$ is $q$-Manin then the relations~\eqref{AMMM_AMMMA} and \eqref{AMMdet_} hold. Multiplying~\eqref{AMMdet_} with~\eqref{AMiMidetMi} and using~\eqref{AMMM_AMMMA} we obtain $A^q_n=A^q_n\det_q M\det_{q^{-1}}(M^{-1})$, which implies~\eqref{detMtetMi}. If $M$ is invertible from the left then $M^{-1}M=1$ for any right inverse $M^{-1}$. Hence multiplying~\eqref{AMiMidetMi} with~\eqref{AMMdet_} and using~\eqref{AMiMiA} we obtain $A^q_n=A^q_n\det_{q^{-1}}(M^{-1})\det_q M$, which implies~\eqref{detMitetM}. \qed

Conversely if the determinant of a $q$-Manin matrix $M$ is left invertible then there exists a matrix $B$ satisfying~\eqref{BAMMA}. It can be given in the form $B=\sum_{i<j\atop k<l}B^{ij}_{kl}(E_{ik}\otimes E_{jl}-qE_{jk}\otimes E_{il})$ where
\begin{align}
 B^{ij}_{kl}=\frac12(-q)^{k+l-i-j}\big(\det_q M\big)^{-1}\det_{q^{-1}}(M_{KJ}),
\end{align}
$K=\backslash(kl)$, $J=\backslash(ij)$. One can check it using the Laplace expansion formula~\eqref{LaplFml2} for the case $m=n-2$ in the form
\begin{align}
\sum_{k<l}(-q)^{k+l-i-j}\det_q M_{KJ}
\det_q \big(M^{ij}_{ab}\big)=(\delta^i_a\delta^j_b-\delta^i_b\delta^j_a)\det_q M.
\label{LaplFml22}
\end{align}

Thus for a two-sided invertible $q$-Manin matrix $M$ the left invertibility of $\det_q M$ is equivalent to the existence of the matrix $B$ satisfying~\eqref{BAMMA}.

\begin{Rem} \normalfont
 The condition~\eqref{BAMMA} is satisfied for some $B\in\mathfrak R\otimes\End(\mathbb C^n\otimes\mathbb C^n)$ if and only if there exists $\wt B\in\mathfrak R\otimes\End(\mathbb C^n\otimes\mathbb C^n)$ such that $\wt BA^qM^{(1)}M^{(2)}\in\End(\mathbb C^n\otimes\mathbb C^n)\backslash\{0\}$. Indeed decomposing $\wt BA^qM^{(1)}M^{(2)}$ as $\alpha A^q+\beta(1-A^q)$ where $\alpha,\beta\in\mathbb C$ and multiplying~\eqref{lemPyat2} by $\wt B$ from the left we obtain $\beta=0$. Hence $\alpha\ne0$ and $B=\wt B/\alpha$ satisfies~\eqref{BAMMA}.
\end{Rem}

%
%


\section{Integrable systems. $L$-operators}
\label{secIntSys}

\subsection{$L$-operators and \qMMs}
\label{secIntSys1}

Let us consider the trigonometric $R$-matrix $R(z)\in\End(\mathbb C^n\otimes\mathbb C^n)(z)$
\footnote{Up to inessential changes the $R$-matrix is the same as in \cite{TV07,OPS07}
but slightly different from \cite{ArnaudonCrampeDoikouFrappRagoucy05,  HM06}.
}
:
\begin{align} \label{Rmatry}
 R(z)= \frac{1}{z-1}   \Bigl( (qz-q^{-1}) \sum_{i=1}^n E_{ii}\otimes E_{ii}+ & (z-1) \sum_{i\ne j} E_{ii}\otimes E_{jj}+ \notag \\
  +(q-q^{-1})&\sum_{i<j}\big(z E_{ij}\otimes E_{ji}+ E_{ji}\otimes E_{ij}\big) \Bigr),
\end{align}
where $E_{ij}$ are basis elements of $\End(\mathbb C^n)$: $E_{ij}e_k=\delta_{jk}e_i$ and $z$ is a complex parameter.

It is convenient sometimes to substitute $z=u/v$:
\begin{align} \label{Rmatruv}
 R(u/v)=\frac{qu-q^{-1}v}{u-v}\sum_{i=1}^n E_{ii}\otimes E_{ii}+&\sum_{i\ne j} E_{ii}\otimes E_{jj}+ \notag \\
  +\frac{q-q^{-1}}{u-v}&\sum_{i<j}\big(u E_{ij}\otimes E_{ji}+v E_{ji}\otimes E_{ij}\big).
\end{align}
The $R$-matrix satisfies the Yang-Baxter equation
\begin{align} \label{YBE}
 R^{(12)}(z_1/z_2)R^{(13)}(z_1/z_3)R^{(23)}(z_2/z_3)=R^{(23)}(z_2/z_3)R^{(13)}(z_1/z_3)R^{(12)}(z_1/z_2).
\end{align}
\begin{Lem}
 The $R$-matrix~\eqref{Rmatry} at $z=q^{-2}$ has the form
\begin{align}
 R(q^{-2})=1-P^q=2A^q, \label{R_A}
\end{align}
where $P^q$ and $A^q$ are defined by the formulae~\eqref{PqDef} and \eqref{AqSqDef}.
\end{Lem}

\noindent{\bf Proof.} Substituting $u=q^{-1}$, $v=q$ to~\eqref{Rmatruv} one obtain 
$$R(q^{-2})=\sum_{i\ne j}\big(E_{ii}\otimes E_{jj}-q^{\s(i-j)} E_{ji}\otimes E_{ij}\big)=1-P^q.$$
\qed \\

Consider an $n\times n$ matrix $L(z)$ with elements in a non-commutative algebra $\mathfrak{R}$ depending on the parameter $z$ (in general the entries of $L(z)$ belong to $\mathfrak{R}[[z,z^{-1}]]$)  which satisfies $RLL$-relation
\begin{align}
 R(z/w)L^{(1)}(z)L^{(2)}(w)=L^{(2)}(w)L^{(1)}(z)R(z/w), \label{RLL}
\end{align}
where $R(z/w)$ is the $R$-matrix~\eqref{Rmatruv}. The matrix $L(z)$ is called {\it $L$-operator} or {\it Lax matrix}. Notice that the Hopf algebra $U_q(\widehat{\mathfrak{gl}_n})$ can be described by a pair of $L$-operators such that the relations~\eqref{RLL} for them are the defining commutation relation for this algebra.

The  basic observation of the present Section is the following
\begin{Prop}
Let $L(z)$ be an $L$-operator satisfying~\eqref{RLL}. Then
\begin{align}
 M=L(z)q^{2z\frac{\partial}{\partial z}}
\end{align}
is a $q$-Manin matrix, where $q^{2z\frac{\partial}{\partial z}}$ is the operator acting as $q^{2z\frac{\partial}{\partial z}}f(z)=f(q^2z)$. Also $\wt M=L(z)^{\top}q^{-2z\frac{\partial}{\partial z}}$ is a $q$-Manin matrix, where $L(z)^{\top}$ is the transpose of $L(z)$.
\end{Prop}


\noindent{\bf Proof.} Substituting $w=q^2z$ to~\eqref{RLL} and using~\eqref{R_A} we obtain
\begin{align}
 A^qL^{(1)}(z)L^{(2)}(q^2z)=L^{(2)}(q^2z)L^{(1)}(z)A^q. \label{ALL_LLA}
\end{align}
Multiplying~\eqref{ALL_LLA} by the operator $q^{4z\frac{\partial}{\partial z}}$ from the right
and using $f(z) \lambda^{z \p} = \lambda^{z \p} f(\lambda^{-1}z)$ one gets:
\begin{align}
 A^q L^{(1)}(z)q^{2z\frac{\partial}{\partial z}} L^{(2)}(z)q^{2z\frac{\partial}{\partial z}} =L^{(2)}(q^2z)L^{(1)}(z)q^{4z\frac{\partial}{\partial z}}A^q. \label{ALL_LLA_q}
\end{align}
Then, applying Lemma~\ref{lemma2} we conclude that the operator $M=L(z)q^{2z\frac{\partial}{\partial z}}$ is $q$-Manin's matrix. Replacing $z\to q^2z$ in the formula~\eqref{ALL_LLA} one can analogously deduce that the matrix $\left.\wt M\right.^\top=L(z)q^{-2z\frac{\partial}{\partial z}}=q^{-2z\frac{\partial}{\partial z}}L(q^2z)$ satisfies the relation~\eqref{AMMAMMA}. \qed

\begin{Rem} \normalfont
 Similarly using~\eqref{AMMAAMM21} and \eqref{AMMAMMAqinv} one can show that if an $L$-operator $L(z)$ satisfies the relation
\begin{align}
 R(z/w)L^{(2)}(z)L^{(1)}(w)=L^{(1)}(w)L^{(2)}(z)R(z/w), \label{RLL21}
\end{align}
(which replaces ~\eqref{RLL}) then the matrices $M=L(z)q^{-2z\frac{\partial}{\partial z}}$ and $\wt M=L(z)^\top q^{2z\frac{\partial}{\partial z}}$ are $q$-Manin matrices.
\end{Rem}

\subsection{Quantum determinant of $L$-operator}

Now we consider the quantum determinant of $L$-operators. They were introduced and extensively studied by L. Faddeev's  school (Kulish, Sklyanin et al.)
and play a fundamental role in the theory of affine quantum groups and corresponding quantum integrable systems. 
Remind that the quantum determinants of the $L$-operators of the Hopf algebra $U_q(\wh{\mathfrak{gl}_n})$ 
describes the center of this algebra. 

Consider an $L$-operator $L(z)$ satisfying~\eqref{RLL}. The {\it quantum determinant} of $L(z)$ is defined by the formula
\begin{align}
\qdet L(z):= \sum_{\tau \in \sg_n} (-q)^{-\inv(\tau)} L_{\tau(1)1}(z) L_{\tau(2)2}(q^2z)\cdots L_{\tau(n)n}(q^{2(n-1)}z). \label{quanDetLDef}
\end{align}
This notion is closely related with the notion of the $q$-determinant~\eqref{Def3eq} of the associated \qMM\ $L(z)\;q^{2nz\frac{\partial}{\partial z}}$.
\begin{Prop}\label{qdetdetq}
The $q$-determinant of  $M=L(z)q^{2z\frac{\partial}{\partial z}}$ satisfies
\begin{align}
 \det_q M=\qdet L(z)\;q^{2nz\frac{\partial}{\partial z}}.
\end{align}
\end{Prop}
{\bf Proof.} The proposition is proven simply noticing that, for any permutation $\tau\in\mathfrak{S}_n$, 
\begin{equation}\begin{split}
&L_{\tau(1),1} \;q^{2z\frac{\partial}{\partial z}}\cdot  L_{\tau(2),2} \;q^{2z\frac{\partial}{\partial z}}\cdots L_{\tau(n),n} \;q^{2z\frac{\partial}{\partial z}}=\\
&(-q)^{-\inv(\tau)} L_{\tau(1)1}(z) L_{\tau(2)2}(q^2z)\cdots L_{\tau(n)n}(q^{2(n-1)}z) \;q^{2nz\frac{\partial}{\partial z}}.\end{split}
\end{equation}
\qed\\
{\bf Remark.} The same result can be obtained using Corollary \ref{AMMdet}. Indeed,
substituting $M=L(z)q^{2z\frac{\partial}{\partial z}}$ in the formula~\eqref{AMMdet_} one gets
\begin{align}
  A^q_n L^{(1)}(z)L^{(2)}(q^2z)\cdots L^{(n)}(q^{2(n-1)}z)=A^q_n\qdet L(z). \label{ALLLdet}
\end{align}
More generally, consider multi-indices $I=(i_1<\ldots<i_m)$ and $J=(j_1<\ldots<j_m)$. It follows from the structure of the $R$-matrix~\eqref{Rmatry} that the submatrix $L_{IJ}(z)$ of an $L$-operator $L(z)$ also satisfies~\eqref{RLL} with the $m^2\times m^2$ $R$-matrix~\eqref{Rmatry}. The corresponding $q$-Manin matrix is the submatrix $M_{IJ}=L_{IJ}(z)q^{2z\frac{\partial}{\partial z}}$ of the $q$-Manin matrix $M=L(z)q^{2z\frac{\partial}{\partial z}}$. So, the quantum minors of $L(z)$ and the $q$-minors of $M$ are related via
\begin{multline}
\det_q(M_{IJ})=\sum_{\tau \in \sg_n} (-q)^{-\inv(\tau)} L_{i_{\tau(1)}j_1}(z) L_{i_{\tau(2)}j_2}(q^2z)\cdots L_{i_{\tau(m)}j_m}(q^{2(m-1)}z)q^{2mz\frac{\partial}{\partial z}}= \\
=\qdet L_{IJ}(z)\;q^{2mz\frac{\partial}{\partial z}}.
\end{multline}
By virtue of Corollary~\ref{CorAMMMm} we can derive
\begin{multline}
 A^q_m L^{(1)}(z)L^{(2)}(q^2z)\cdots L^{(m)}(q^{2(m-1)}z)= \\
=m!\sum_{I=(i_1<\ldots<i_m)\atop J=(j_1<\ldots<j_m)}
\qdet L_{IJ}(z)\; A^q_m\big(E_{i_1j_1}\otimes\cdots\otimes E_{i_mj_m}\big)A^q_m.
\end{multline} 
\begin{Ex} \normalfont
As an application of this formul\ae\ we consider the the Gauss decomposition of an $L$-operator:
\begin{align}
L(z)=\left(\begin{array}{ccc}
 1 &{}&F_{\alpha \beta}(z)\\
{}&\ddots&{}\\
0&{}&1
\end{array}\right)
\left(\begin{array}{ccc}
k_1(z)&{}&0\\
{}&\ddots&{}\\
0&{}&k_n(z)
\end{array}\right)
\left(\begin{array}{ccc}
1&{}&0\\
{}&\ddots&{}\\
E_{\beta \alpha}(z)&{}&1
\end{array}\right)= \\
=\left(\begin{array}{ccc}
 1 &{}&0\\
{}&\ddots&{}\\
E'_{\alpha \beta}(z) &{}&1
\end{array}\right)
\left(\begin{array}{ccc}
k'_1(z)&{}&0\\
{}&\ddots&{}\\
0&{}&k'_n(z)
\end{array}\right)
\left(\begin{array}{ccc}
1&{}&F'_{\beta \alpha}(z)\\
{}&\ddots&{}\\
0&{}&1
\end{array}\right).
\end{align}
If $k_i(z)$ and $k'_i(z)$ are invertible Propositions~\ref{Gauss-pr} and \ref{qdetdetq} imply
\begin{align}
\qdet L(z)=k_n(z)k_{n-1}(q^2z)\cdots k_1(q^{2(n-1)}z)\;q^{-2 n z\frac{\partial}{\partial z}}\\
=k'_1(z)k'_2(q^2z)\cdots k'_n(q^{2(n-1)}z)\;q^{-2nz\frac{\partial}{\partial z}}.
\end{align}
\end{Ex}

\subsection{Generating functions for the commuting operators}

As it well known in the theory of quantum Yang-Baxter equations, the $\mathfrak R$-valued function
\begin{align} \label{t1z_def}
 t_1(z)=\tr L(z)
\end{align}
commute with itself for different values of the spectral parameters, i.e., 
\begin{align}
 t_1(z)t_1(w)=t_1(w)t_1(z), \> \forall\, z,w\in \mathbb{C}.
\end{align}
It means that it generates the commuting operators:
\begin{align}
 t_1(z)=\sum_{k\in\mathbb Z}H_k z^k, &&&H_kH_l=H_lH_k. \label{HkHl}
\end{align}
Hence, whenever one of these operators is the Hamiltonian of a quantum system, we see that this Hamiltonian admits 
a family of commuting operators (in the classical limit, they are usually called constants of the motion). 
However, if the size $n$ of the $L$ operator is  greater than two,
then they are generally not enough to ensure simplicity of the joint spectrum (in the quantum case) and/or guarantee the applicability 
of the Arnol'd-Liouville theorem in the classical case.
Here we construct the highest generating functions of commuting operators 
(integrals of motion) in terms of the  q-Manin matrix $M=L(z)q^{2z\frac{\partial}{\partial_z}}$.

Let us introduce the following matrix acting in the space $(\mathbb C^n)^{\otimes m}$:
\begin{multline}
 \mathbb R_m(z_1,\ldots,z_m)=\mathop{\overrightarrow\prod}\limits_{1\le i<j\le m}R^{(i,j)}(z_i/z_j)= \\
 =\big(R^{(1,2)}R^{(1,3)}\cdots R^{(1,m)}\big)\cdots\big(R^{(m-2,m-1)}R^{(m-2,m-1)}\big)\big(R^{(m-1,m)}\big)= \\
 =\big(R^{(1,2)}\big)\big(R^{(1,3)}R^{(2,3)}\big)\cdots\big(R^{(1,m)}\cdots R^{(m-1,m)}\big), \label{RR_def}
\end{multline}
where $R^{(i,j)}=R^{(i,j)}(z_i/z_j)$. Note that one can rearrange it in the form
\begin{multline}
 \mathbb R_m(z_1,\ldots,z_m)=\mathbb R_k(z_1,\ldots,z_k)
\mathop{\overrightarrow\prod}\limits_{1\le i\le k}\;
\mathop{\overrightarrow\prod}\limits_{k+1\le j\le m}R^{(i,j)}(z_i/z_j)\;\;
\mathbb R_{m-k}^{(k+1,\ldots,m)}(z_{k+1},\ldots,z_m), \label{RRkm}
\end{multline}
where $k\le m$. Using the Yang-Baxter equation~\eqref{YBE} we can rewrite~\eqref{RRkm} in another order:
\begin{multline}
 \mathbb R_m(z_1,\ldots,z_m)=\mathbb R_k(z_1,\ldots,z_k)
\mathbb R_{m-k}^{(k+1,\ldots,m)}(z_{k+1},\ldots,z_m)
\mathop{\overrightarrow\prod}\limits_{1\le i\le k}\;
\mathop{\overleftarrow\prod}\limits_{k+1\le j\le m}R^{(i,j)}(z_i/z_j)= \\
=\mathop{\overleftarrow\prod}\limits_{1\le i\le k}\;
\mathop{\overrightarrow\prod}\limits_{k+1\le j\le m}R^{(i,j)}(z_i/z_j)\;\;
\mathbb R_k(z_1,\ldots,z_k) \mathbb R_{m-k}^{(k+1,\ldots,m)}(z_{k+1},\ldots,z_m). \label{RRrevor}
\end{multline}
There holds (see, e.g.~\cite{ArnaudonCrampeDoikouFrappRagoucy05}) the
\begin{Lem} \label{lemRA}  
\begin{align}
 \mathbb R_m(z,q^2z,\ldots,q^{2(m-k-1)}z)=m!A^q_m. \label{lemRAkm}
 \end{align}
\end{Lem}

\noindent This lemma is proved in the Appendix~\ref{AppLemRA}.

Substituting $m=k+l$, $z_i=q^{2(i-1)}z$ for $i=1,\ldots,k$ and $z_{k+i}=q^{2(i-1)}w$ for $i=1,\ldots,l$ in~\eqref{RRrevor} and taking into account~\eqref{lemRAkm} we obtain
\begin{align} \label{AklRkl}
 A^q_{k,l}\mathbb R_{k,l}(z,w)=A^q_{k,l}\mathbb R_{k,l}(z,w)A^q_{k,l},
\end{align}
where $A^q_{k,l}=\big(A^q_k\big)^{(1,\ldots,k)}\big(A^q_l\big)^{(k+1,\ldots,k+l)}$ and
\begin{align}
 \mathbb R_{k,l}(z,w)&=\mathop{\overrightarrow\prod}\nolimits_{1\le i\le k} \mathop{\overleftarrow\prod}\nolimits_{k+1\le j\le k+l}
 R^{(i,j)}(q^{2(i-j+k)}z/w).
\end{align}
Similarly one gets
\begin{align} \label{AklRklinv}
 A^q_{k,l}\mathbb R_{k,l}(z,w)^{-1}=A^q_{k,l}\mathbb R_{k,l}(z,w)^{-1}A^q_{k,l}.
\end{align}

Writing down the relation~\eqref{AMMM_AMMMA} for the $q$-Manin matrix $M=L(z)q^{2z\frac{\partial}{\partial_z}}$ we obtain
\begin{align}
  A^q_m L^{(1)}(z)L^{(2)}(q^2z)\cdots L^{(m)}(q^{2(m-1)}z)=A^q_m L^{(1)}(z)L^{(2)}(q^2z)\cdots L^{(m)}(q^{2(m-1)}z)A^q_m. \label{ALLL_ALLLA}
\end{align}
This implies
\begin{align} \label{AklLkl}
 A^q_{k,l}\mathbb L_{k,l}(z,w)=A^q_{k,l}\mathbb L_{k,l}(z,w)A^q_{k,l},
\end{align}
where
\begin{multline}
\mathbb L_{k,l}(z,w)= \\
 =L^{(1)}(z)L^{(2)}(q^2z)\cdots L^{(k)}(q^{2(k-1)}z)L^{(k+1)}(w)L^{(k+2)}(q^2w)\cdots L^{(k+l)}(q^{2(l-1)}w). \label{LzwDef}
\end{multline}

\begin{Rem} \normalfont
 If $L(z)$ is an $L$-operator satisfying~\eqref{RLL21} (not usual ~\eqref{RLL}) then the  substitution of the $q^{-1}$-Manin matrix $M=L(z)q^{-2z\frac{\partial}{\partial z}}$ to~\eqref{AMMM_AMMMA21} yields
\begin{multline}
  A^q_m L^{(m)}(z)L^{(m-1)}(q^{-2}z)\cdots L^{(1)}(q^{-2(m-1)}z)= \\
=A^q_m L^{(m)}(z)L^{(m-1)}(q^{-2}z)\cdots L^{(1)}(q^{-2(m-1)}z)A^q_m. \label{ALLL_ALLLA21}
\end{multline}
\end{Rem}

\begin{Prop}
If $L(z)$ is an $L$-operator over $\mathfrak R$ satisfying~\eqref{RLL} then the $\mathfrak R$-valued functions
\begin{align} \label{tkz_def}
 t_k(z)=\tr_{1,\ldots,k}\big(A^q_k L^{(1)}(z)L^{(2)}(q^2z)\cdots L^{(k)}(q^{2(k-1)}z)\big)
\end{align}
commute among themselves for different values of the spectral parameter:
\begin{align} \label{tkztlw}
 t_k(z)t_l(w)&=t_l(w)t_k(z), & &k,l=1,2,3,\ldots
\end{align}
\end{Prop}

\noindent{\bf Proof.} Due to the $RLL$-relation~\eqref{RLL} the products of $L$-operators~\eqref{LzwDef} and
\begin{align*}
\wt{\mathbb L}_{k,l}(z,w)&=L^{(k+l)}(q^{2(l-1)}w)\cdots L^{(k+2)}(q^2w)L^{(k+1)}(w)L^{(k)}(q^{2(k-1)}z)\cdots L^{(2)}(q^2z)L^{(1)}(z)
\end{align*}
are related as
\begin{align} \label{RLLL}
 \mathbb R_{k,l}(z,w)\mathbb L_{k,l}(z,w)
  =\wt{\mathbb L}_{k,l}(z,w)\mathbb R_{k,l}(z,w).
\end{align}
Multiplying this equation by $\mathbb R_{k,l}(z,w)^{-1}$ from the right, by $A^q_{k,l}$ from the left and taking the trace over all the spaces we obtain the equality
\begin{align} \label{tr_RAALLL}
 \tr_{1,\ldots,k+l}\big(A^q_{k,l}\mathbb R_{k,l}(z,w)\mathbb L_{k,l}(z,w)\mathbb R_{k,l}(z,w)^{-1}\big)
 =\tr_{1,\ldots,k+l}\big(A^q_{k,l}\wt{\mathbb L}_{k,l}(z,w)\big).
\end{align}
Using the relations~\eqref{AklRkl}, \eqref{AklLkl}, \eqref{AklRklinv} and periodicity of the trace we rewrite the left hand side of~\eqref{tr_RAALLL} as
\begin{multline}
\tr_{1,\ldots,k+l}\big(A^q_{k,l}\mathbb R_{k,l}(z,w)A^q_{k,l}\mathbb L_{k,l}(z,w)\mathbb R_{k,l}(z,w)^{-1}\big)= \\
=\tr_{1,\ldots,k+l}\big(A^q_{k,l}\mathbb R_{k,l}(z,w)A^q_{k,l}\mathbb L_{k,l}(z,w)A^q_{k,l}\mathbb R_{k,l}(z,w)^{-1}\big)= \\
=\tr_{1,\ldots,k+l}\big(\mathbb R_{k,l}(z,w)A^q_{k,l}\mathbb L_{k,l}(z,w)A^q_{k,l}\mathbb R_{k,l}(z,w)^{-1}A^q_{k,l}\big)= \\
=\tr_{1,\ldots,k+l}\big(\mathbb R_{k,l}(z,w)A^q_{k,l}\mathbb L_{k,l}(z,w)A^q_{k,l}\mathbb R_{k,l}(z,w)^{-1}\big)= \\
=\tr_{1,\ldots,k+l}\big(A^q_{k,l}\mathbb L_{k,l}(z,w)A^q_{k,l}\big)
=\tr_{1,\ldots,k+l}\big(A^q_{k,l}\mathbb L_{k,l}(z,w)\big)=t_k(z)t_l(w).
\end{multline}
Since the right hand side of~\eqref{tr_RAALLL} equals $t_l(w)t_k(z)$ we obtain~\eqref{tkztlw}. \qed

\begin{Rem} \normalfont
 The relations~\eqref{AklRkl}, \eqref{AklRklinv} follow from the relation~\eqref{ALLL_ALLLA} for the $L$-operators $L^{(0)}(z)=\mathop{\overleftarrow\prod}\nolimits_{1\le j\le l}
 R^{(0,j)}(q^{2(1-j)}z/w)$,
 $L^{(0)}(w)=\mathop{\overleftarrow\prod}\nolimits_{1\le i\le k}R^{(i,0)}(q^{2(i-1)}z/w)^{-1}$ (satisfying~\eqref{RLL}) and from~\eqref{ALLL_ALLLA21} for the $L$-operators $L^{(0)}(w)=\mathop{\overrightarrow\prod}\nolimits_{1\le i\le k}
 R^{(i,0)}(q^{2(i-l)}z/w)$, $L^{(0)}(z)=\mathop{\overrightarrow\prod}\nolimits_{1\le j\le l}
 R^{(0,j)}(q^{2(k-j)}z/w)^{-1}$ (satisfying~\eqref{RLL21}). Note that in this way we do not need Lemma~\ref{lemRA} since we use the properties of $q$-Manin matrices instead.
\end{Rem}

The functions~\eqref{tkz_def} are related to the sums of its principal $q$-minors of the \qMM{} $M=L(z)q^{2z\frac{\partial}{\partial z}}$ by the formula
\begin{align} \label{tkz_M}
 t_k(z)=\tr_{1,\ldots,k}\big(A^q_k M^{(1)}\cdots M^{(k)}\big)
  q^{-2kz\frac{\partial}{\partial z}}=
  \sum\limits_{I=(i_1<\ldots<i_m)}\det_q(M_{II})\,
  q^{-2kz\frac{\partial}{\partial z}}.
\end{align}

The functions~\eqref{tkz_def} can be regarded as a generated functions for the integrals of motion for the system defined by the Hamiltonians~\eqref{HkHl}. They can be in turn gathered in the two-variable function
\begin{align}
g(z,u)=\sum_{m=0}^n(-1)^m t_m(z) u^{n-m}. \label{CharPol_g}
\end{align}
It is related with the characteristic polynomial:
\begin{align}
\Char_M(u)=\sum_{m=0}^n(-1)^m t_m(z) q^{2mz\frac{\partial}{\partial z}}u^{n-m}=g(z,uq^{-2z\frac{\partial}{\partial z}})q^{2nz\frac{\partial}{\partial z}}. \label{CharPol_tm}
\end{align}
Comparing it with~\eqref{CharPol} we obtain
\begin{align} \label{emtm}
 \mathfrak e_m=t_m(z)q^{2mz\frac{\partial}{\partial z}}.
\end{align}
Hence the functions~\eqref{tkz_def} equal to sums of the corresponding quantum minors of $L(z)$:
\begin{align}
 t_m(z)=\sum_{I=(i_1<\ldots<i_m)}\qdet L_{II}(z).
\end{align}

\subsection{Quantum powers for $L$-operators}

Here we apply the Cayley-Hamilton theorem and Newton identities to the 
$q$-Manin matrix $M=L(z)q^{2z\frac{\partial}{\partial z}}$. In particular, the Newton identities give us a 
new family of generating functions of integrals of motions, which provide a suitable quantum version of 
the powers of classical $L$-operators.

Let us use the notation
$L^{[m]}(z)$ defined iteratively by the formulae
\begin{align}
L^{[0]}(z)&=1, & &L^{[m]}(z)=L^{[m-1]}(z)\qstar L(q^{2(m-1)}z)= \\
& & &=\big(\ldots\big(\big(L(z)\qstar L(q^2z)\big)\qstar L(q^4z)\big)\qstar\ldots\qstar L(q^{2(m-1)}z)\big), \notag
\end{align}
( the notion of the $\qstar$-product was introduced in Definition \ref{def8}).
This definition is related with the $q$-powers of $M$:
\begin{align} \label{Mm_Lm}
M^{[m]}=L^{[m]}(z)q^{2mz\frac{\partial}{\partial z}}.
\end{align}
Applying Theorems~\ref{propCharPol_qPower} and \ref{propNewton} for $M=L(z)q^{2z\frac{\partial}{\partial z}}$ and taking into account the formulae~\eqref{emtm}, \eqref{Mm_Lm} one obtains
\begin{gather}
 \sum_{m=0}^n(-1)^m t_m(z) L^{[n-m]}(q^{2m}z)=0, \\
 mt_m(z)=\sum_{k=0}^{m-1}(-1)^{m+k+1} t_k(z)\tr\big(L^{[m-k]}(q^{2k}z)\big). \label{L_Newton}
\end{gather}

Let us consider the $\mathfrak R$-valued functions
\begin{align} \label{Ikz_def}
 I_k(z)&=\tr\big(L^{[k]}(z)\big), & &k=1,2,\ldots.
\end{align}
Then the Newton identities~\eqref{L_Newton} can be rewritten in terms of~\eqref{Ikz_def} as
\begin{gather}
 mt_m(z)=\sum_{k=0}^{m-1}(-1)^{m+k+1} t_k(z) I_{m-k}(q^{2k}z). \label{L_Newton_Ikz}
\end{gather}
As we can recurrently express the elements $\mathfrak e_m$ and $\tr(M^{[k]})$ through each other using the Newton identities~\eqref{Newton_th} the formula~\eqref{L_Newton_Ikz} allows us to express $t_m(z)$ through $I_k(z)$ and vice versa. This leads to the following.

\begin{Th}
The functions~\eqref{Ikz_def} commute with themselves, with each other and with the functions~\eqref{tkz_def} for different values of parameters:
\begin{align}
 I_k(z)I_l(w)=I_l(w)I_k(z),  \label{Ikz_Ilw}\\
 I_k(z)t_l(w)=t_l(w)I_k(z), \label{Ikz_tlw}
\end{align}
where $k,l\in\mathbb Z_{>0}$.
\end{Th}

\noindent{\bf Proof.} Let us prove this theorem by induction. First we note that $I_1(z)=t_1(z)$ and $I_1(z)$ commute with itself (for different values of parameters) and with $t_l(w)$. Further suppose that $I_1(z)$, $I_2(z)$, \ldots, $I_{m-1}(z)$ commute with itself, each other and with $t_l(w)$.
Since $t_0(z)=1$, this allows to express the function $I_m(z)$ through $I_1(z)$, \ldots, $I_{m-1}(z)$ and $t_l(w)$, $l=1,\ldots,m$:
\begin{gather}
 I_{m}(z) =(-1)^{m+1}mt_m(z)-\sum_{k=1}^{m-1}(-1)^{k} t_k(z) I_{m-k}(q^{2k}z). \label{Imz_expr}
\end{gather}
This implies that $I_1(z)$, $I_2(z)$, \ldots, $I_{m}(z)$ commute with the functions $t_l(w)$ and with itself and each other for different $z$. \qed

The functions $I_k(z)$ can be considered as generating another set of quantum integrals of motion for the system defined by $t_l(z)$. 
They form a set of integrals alternative to the set of integrals generated by the functions $t_l(z)$.


\appendix

\section{Proof of Lemma~\ref{lemRA}}
\label{AppLemRA}

For each $s=1,\ldots,m-1$ one can move the factor $R^{(s,s+1)}$ in the product~\eqref{RR_def} to the right using the Yang-Baxter equation~\eqref{YBE}. Since $R(q^{-2})=2A^q_2$ this implies $\mathbb R_m=-\mathbb R_m\pi_q(\sigma_s)$, where $\pi_q$ is the representation of symmetric group defined by~\eqref{piqDef}. This means that for any $\sigma\in \sg_m$ we have the identity $\mathbb R_m=(-1)^\sigma\mathbb R_m\pi_q(\sigma)$. In the other hand, the operator $A^q_m$ satisfy the same identity $A^q_m=(-1)^\sigma A^q_m\pi_q(\sigma)$. Since $\pi_q(\sigma_{sr})e_{\ldots,i,\ldots,i,\ldots}=0$ (with subscripts $i$ placed onto the $s$-th and $r$-th sites) we have $\mathbb R_m e_{\ldots,i,\ldots,i,\ldots}=0=m!A^q_m e_{\ldots,i,\ldots,i,\ldots}$. By the same reason taking into account the formula~\eqref{pi_q_tau_e} we conclude that it is sufficient to check the equality~\eqref{lemRAkm} on the vectors $e_{i_1,\ldots,i_m}$ with $i_1<\ldots<i_m$. To do it we suppose by induction that the equality $\mathbb R_{m-1}=(m-1)!A^q_{m-1}$ is already proven. Let us substitute the explicit expression for the $R$-matrices to the result of the action of $\mathbb R_m$ written in the form
\begin{multline} \label{Rm_e}
 \mathbb R_m e_{i_1,\ldots,i_m}=\mathbb R_{m-1} R^{(1m)}(q^{2(1-m)})R^{(2m)}(q^{2(2-m)})\cdots R^{(m-1,m)}(q^{-2}) e_{i_1,\ldots,i_m}= \\
 =(m-1)!A^q_{m-1}R^{(1m)}(q^{2(1-m)})R^{(2m)}(q^{2(2-m)})\cdots R^{(m-1,m)}(q^{-2}) e_{i_1,\ldots,i_m}.
\end{multline}
Due to the inequality $i_1<\ldots<i_m$ each vector $R^{(s+1,m)}(q^{2(s+1-m)})\cdots R^{(m-1,m)}(q^{-2}) e_{i_1,\ldots,i_m}$ is a linear combinations of the vectors $e_{j_1,\ldots,j_m}$ where $j_1,\ldots,j_m$ are pairwise different and $j_1<\ldots<j_s<j_m$. Then we can proceed as follows. Substituting the expression~\eqref{Rmatry} for $R^{(sm)}(q^{2(s-m)})$ we see that the first sum and first term of third sum in~\eqref{Rmatry} act by zero, the second sum acts identically and the second term of the third sum acts as $\dfrac{q-q^{-1}}{q^{2(s-m)}-1}\pi_1(\sigma_{sm})$, where $\pi_1$ is a representation $\pi_q$ at $q=1$: $$ \pi_1(\sigma)e_{j_1,\ldots,j_m}=e_{j_{\sigma^{-1}(1)},\ldots,j_{\sigma^{-1}(m)}}.$$ Thus,  the expression~\eqref{Rm_e} takes the form
\begin{multline} \label{Rm_e_}
  (m-1)!A^q_{m-1}
\Big(1+\frac{q-q^{-1}}{q^{2(1-m)}-1}\pi_1(\sigma_{1m})\Big)\cdots\Big(1+\frac{q-q^{-1}}{q^{-2}-1}\pi_1(\sigma_{m-1,m})\Big)
e_{i_1,\ldots,i_m}.
\end{multline}
Opening the big parentheses we obtain
\begin{multline} \label{Rm_e__}
  (m-1)!A^q_{m-1}\bigg(1+\sum_{g=1}^{m-1}\sum_{1\le k_1<\ldots<k_g\le m-1}\prod_{l=1}^g\frac{q-q^{-1}}{q^{2(k_l-m)}-1}\;\pi_1(\sigma_{k_1m}\cdots\sigma_{k_gm})\bigg)
e_{i_1,\ldots,i_m}.
\end{multline}
Taking into account $\pi_1(\sigma)e_{i_1,\ldots,i_m}=(q)^{-\inv(\sigma)}\pi_q(\sigma)e_{i_1,\ldots,i_m}$ and $\inv(\sigma_{k_1m}\cdots\sigma_{k_gm})=2m-2k_1-g$ and $A^q_{m-1}\pi_q(\sigma_{k_1m}\cdots\sigma_{k_gm})=(-1)^{g-1}A^q_{m-1}\pi_q(\sigma_{k_1m})$ we obtain
\begin{multline} \label{Rm_e___}
  (m-1)!A^q_{m-1}\bigg(1-\sum_{g=1}^{m-1}\sum_{1\le k_1<\ldots<k_g\le m-1}q^{2(k_1-m)}\prod_{l=1}^g\frac{1-q^2}{q^{2(k_l-m)}-1}\;\pi_q(\sigma_{k_1m})\bigg)
e_{i_1,\ldots,i_m}.
\end{multline}
The sum over $g$ can be calculated as follows:
\begin{multline} \label{Rm_e_sum_g}
  \sum_{g=1}^{m-1}\sum_{1\le k_1<\ldots<k_g\le m-1}q^{2(k_1-m)}\prod_{l=1}^g\frac{1-q^2}{q^{2(k_l-m)}-1}\;\pi_q(\sigma_{k_1m})= \\
    =\sum_{k_1=1}^{m-1}\pi_q(\sigma_{k_1m})q^{2(k_1-m)}\frac{1-q^2}{q^{2(k_1-m)}-1}\sum_{g=1}^{m-k_1}\sum_{k_1+1\le k_2<\ldots<k_g\le m-1}\prod_{l=2}^g\frac{1-q^2}{q^{2(k_l-m)}-1}= \\
    =\sum_{k_1=1}^{m-1}\pi_q(\sigma_{k_1m})q^{2(k_1-m)}\frac{1-q^2}{q^{2(k_1-m)}-1}\prod_{s=k_1+1}^{m-1}\Big(1+\frac{1-q^2}{q^{2(s-m)}-1}\Big)=\sum_{k_1=1}^{m-1}\pi_q(\sigma_{k_1m}).
\end{multline}
Finally, we have
\begin{align} \label{Rm_e_fin}
  \mathbb R_m e_{i_1,\ldots,i_m}=(m-1)!A^q_{m-1}\bigg(1-\sum_{k=1}^{m-1}\pi_q(\sigma_{km})\bigg)e_{i_1,\ldots,i_m}=m!A^q_m e_{i_1,\ldots,i_m}.
\end{align}
Thus we have proved the formula~\eqref{lemRAkm} on all the basis $\{e_{i_1,\ldots,i_m}\}$.

\section{An alternative proof of the Lagrange-Desnanot-Jacobi-Lewis Caroll formula}
\label{AppPropLDJLC}

Here we present an alternative proof of the Lagrange-Desnanot-Jacobi-Lewis Caroll formula expressed in the form~\eqref{remLDJLC1}--\eqref{remLDJLC3} that demands the slightly weaker condition of right-invertibility on the \qMM{} $M$.

\begin{Prop} \label{PropLDJLC} 
Let $M$ be a \qMM{} right -invertible. Then
\begin{align}
 M^{adj}_{i_1j}M^{-1}_{i_2j}-qM^{adj}_{i_2j}M^{-1}_{i_1j}&=0, \label{lemLDJLC1} \\
 M^{adj}_{i_1j_1}M^{-1}_{i_2l_2}-qM^{adj}_{i_2j_1}M^{-1}_{i_1l_2}&=(-q)^{j_1+l_2-i_1-i_2} \det_q\big(M_{\backslash (j_1l_2)\backslash (i_1i_2)}\big), \label{lemLDJLC2} \\
 M^{adj}_{i_1l_2}M^{-1}_{i_2j_1}-qM^{adj}_{i_2l_2}M^{-1}_{i_1j_1}&=(-q)^{j_1+l_2-i_1-i_2+1} \det_q\big(M_{\backslash(j_1l_2)\backslash(i_1i_2)}\big), \label{lemLDJLC3}
\end{align}
where $1\le i_1<i_2\le n$,~~~$1\le j\le n$ and $1\le j_1<j_2\le n$.
\end{Prop}

\noindent{\bf Proof.} Consider $q$-Grassmann variables $\psi_i$. Since the matrix $M$ is $q$-Manin the variables $\tpsi_j$ are also $q$-Grassmann. The right invertibility of $M$ implies the equality $\sum_{a=1}^n\tpsi_a M^{-1}_{as}=\psi_s$. Multiplying it by $\psi_r\tpsi_1\cdots\tpsi_{i-1}\tpsi_{i+1}\cdots\tpsi_{k-1}\tpsi_{k+1}\cdots\tpsi_n$ from the left we obtain
\begin{align} \label{lemLDJLC_R1}
 \sum_{a=1}^n\psi_r\tpsi_1\cdots\tpsi_{i-1}\tpsi_{i+1}\cdots\tpsi_{k-1}\tpsi_{k+1}\cdots\tpsi_n\tpsi_a M^{-1}_{as}=\psi_r\tpsi_1\cdots\tpsi_{i-1}\tpsi_{i+1}\cdots\tpsi_{k-1}\tpsi_{k+1}\cdots\tpsi_n\psi_s.
\end{align}
The left hand side of this equation contains only two non-vanishing items, for $a=k$ and $a=i$:
\begin{multline} \label{LDJLCR2}
 \sum_{a=1}^n\psi_r\tpsi_1\cdots\tpsi_{i-1}\tpsi_{i+1}\cdots\tpsi_{k-1}\tpsi_{k+1}\cdots\tpsi_n\tpsi_a M^{-1}_{as}= \\
 =(-q)^{k-n}\psi_r\tpsi_1\cdots\tpsi_{i-1}\tpsi_{i+1}\cdots\tpsi_n M^{-1}_{ks}
 +(-q)^{i-n+1}\psi_r\tpsi_1\cdots\tpsi_{k-1}\tpsi_{k+1}\cdots\tpsi_n M^{-1}_{is}.
\end{multline}
Due to the relation
\begin{align}
 \psi_r\tpsi_1\cdots\tpsi_{i-1}\tpsi_{i+1}\cdots\tpsi_n&=
 \psi_r\psi_1\cdots\psi_{r-1}\psi_{r+1}\cdots\psi_n\det_q M^{\backslash r}_{\backslash i}= \\
 &=(-q)^{i-r}\psi_r\psi_1\cdots\psi_{r-1}\psi_{r+1}\cdots\psi_n M^{adj}_{ir}
\end{align}
we can rewrite~\eqref{LDJLCR2} as
\begin{align} \label{LDJLCR4}
 (-q)^{i+k-r-n}\psi_r\psi_1\cdots\psi_{r-1}\psi_{r+1}\cdots\psi_n \big(M^{adj}_{ir} M^{-1}_{ks} -q M^{adj}_{kr} M^{-1}_{is}\big).
\end{align}
From the other hand, if $r=s$ the right hand side of~\eqref{lemLDJLC_R1} vanish and, therefore, we obtain~\eqref{lemLDJLC1}. Otherwise the right hand side of~\eqref{lemLDJLC_R1} is equal to
\begin{multline} \label{lemLDJLC_R5}
 \psi_r\tpsi_1\cdots\tpsi_{i-1}\tpsi_{i+1}\cdots\tpsi_{k-1}\tpsi_{k+1}\cdots\tpsi_n\psi_s= \\
 =\psi_r\psi_1\cdots\psi_{r-1}\psi_{r+1}\cdots\psi_{s-1}\psi_{s+1}\cdots\psi_n\psi_s\det_q M^{\backslash rs}_{\backslash ik}.
\end{multline}
If $r<s$ then
\begin{align}
 \psi_r\psi_1\cdots\psi_{r-1}\psi_{r+1}\cdots\psi_{s-1}\psi_{s+1}\cdots\psi_n\psi_s=(-q)^{s-n}\psi_r\psi_1\cdots\psi_{r-1}\psi_{r+1}\cdots\psi_n,
\end{align}
and substituting $r=j$, $s=l$ we obtain~\eqref{lemLDJLC2}.
If $r>s$ then
\begin{align}
 \psi_r\psi_1\cdots\psi_{r-1}\psi_{r+1}\cdots\psi_{s-1}\psi_{s+1}\cdots\psi_n\psi_s=(-q)^{s-n+1}\psi_r\psi_1\cdots\psi_{r-1}\psi_{r+1}\cdots\psi_n,
\end{align}
and substituting $r=l$, $s=j$ we obtain~\eqref{lemLDJLC3}. \qed


\begin{thebibliography}{99}

\bibitem{ArnaudonCrampeDoikouFrappRagoucy05}
D. Arnaudon, N. Crampe, A. Doikou, L. Frappat, E. Ragoucy,
{\em Spectrum and Bethe ansatz equations for the $U_q(gl(N))$
closed and open spin chains in any representation},
 Ann. Henri Poincar\'e {\bf 7} (2006), 1217--1268.\\
\href{http://arxiv.org/abs/math-ph/0512037v3}{math-ph/0512037}
\bibitem{Bres99}
D. Bressoud,
{\em Proofs and Confirmations: The Story of the Alternating Sign Matrix Conjecture},
Cambridge University Press, 1999
\bibitem{AC12} A. Chervov, {\em Decomplexification of the Capelli identities and holomorphic factorization}, arXiv 1203.5759.
\bibitem{CF07}
A. Chervov, G. Falqui,
{\em Manin matrices and Talalaev's formula},
J. Phys. A: Math. Theor. \textbf{41} (2008) 194006\\
\href{http://arxiv.org/abs/0711.2236}{arXiv:0711.2236}
\bibitem{CFR08}
A. Chervov, G. Falqui, V. Rubtsov,
{\em Algebraic properties of Manin matrices 1},
Adv. Appl. Math. 43 (2009) 239–315.\\
\bibitem{CFRy}
A. Chervov, G. Falqui, L. Rybnikov,
{\em
Limits of Gaudin algebras, quantization of bending flows, Jucys--Murphy elements and Gelfand--Tsetlin bases},
 Lett. Math. Phys. 91 (2010), 129–150.\\
\href{http://lanl.arxiv.org/abs/0710.4971}{arXiv:0710.4971}

{\em Limits of Gaudin systems: classical and quantum
cases}, SIGMA Symmetry Integrability Geom. Methods Appl. 5 (2009), Paper 029, 17 pp.
\bibitem{CM}
A. Chervov, A. Molev, {\em On higher order Sugawara operators},
IMRN, no.9, (2009), 1612-1635.\\
 \href{http://lanl.arxiv.org/abs/0808.1947}{arXiv:0808.1947}, 
\bibitem{CP94}
V. Chari, A. Pressley,
{\em A guide to quantum groups},
Cambridge University Press, Cambridge, 1994.
\href{http://books.google.com/books?id=YQSKPnFzDOEC&printsec=frontcover&dq=A+guide+to+quantum+groups&sig=0wLtF19lwtBDVu0lHFwd7kKXspk}
{http://books.google.com/\ldots}

\bibitem{CSS08}
S. Caracciolo, A. Sportiello, A. D. Sokal,
{\em Noncommutative determinants, Cauchy-Binet formulae, and Capelli-type identities. I.
Generalizations of the Capelli and Turnbull identities},
 Electron. J. Combin. 16 (2009), Research Paper 103, 43 pp. 
\\
\href{http://arxiv.org/abs/0809.3516}{arXiv:0809.3516}
\bibitem{CT06-1}
A.Chervov, D. Talalaev {\em Quantum spectral curves, quantum integrable systems  and
the geometric Langlands correspondence },
\href{http://arxiv.org/abs/hep-th/0604128}{hep-th/0604128}


\bibitem{D1866}
C.L. Dodgson, {\em Condensation of Determinants},
Proceedings of the Royal Society of London
\textbf{15}, (1866), 150-155.\\
\href{http://www.jstor.org/pss/112607}{http://www.jstor.org/pss/112607}
\bibitem{DI}
D. Zeilberger,
{\em Reverend Charles to the aid of Major Percy and Fields-Medalist Enrico},
 Amer. Math. Monthly  \textbf{103},  (1996),  no. 6, 501--502.
\href{http://arxiv.org/abs/math/9507220}{arXiv:math/9507220}

D. Zeilberger
{\em Dodgson's Determinant-Evaluation Rule proved by Two-Timing Men and Women},
 Electron. J. Combin. 4 (1997), no. 2, Research Paper 22, approx. 2 pp.\\
\href{http://arxiv.org/abs/math/9808079}{arXiv:math/9808079}

 A.N.W. Hone,
{\em Dodgson condensation, alternating signs and square ice},
 Philos. Trans. R. Soc. Lond. Ser. A Math. Phys. Eng. Sci. \textbf{364}, (2006), no. 1849, 3183--3198.\\
\href{http://dx.doi.org/10.1098/rsta.2006.1887}{http://dx.doi.org/10.1098/rsta.2006.1887}

R. E. Schwartz
{\em Discrete Monodromy, Pentagrams, and the Method of Condensation},
 J. Fixed Point Theory Appl. 3 (2008), 379–409.\\
\href{http://arxiv.org/abs/0709.1264}{arXiv:0709.1264}

K. Said, A. Salem, R. Belgacem
{\em  A Mathematical Proof of Dodgson's Algorithm},
\href{http://arxiv.org/abs/0712.0362}{arXiv:0712.0362}
\bibitem{FRT}
N.Yu, Reshetikhin, L.A. Takhtadzhyan, L.D, Faddeev, 
{\em Quantization of Lie groups and Lie algebras}.
Leningr. Math. J. 1,  193-225 (1990); translation from Algebra Anal. 1, No.1, 178-206 (1989).
\bibitem{FT87}
Faddeev, L. D.; Takhtajan, L. A.
{\em Hamiltonian methods in the theory of solitons},
Springer-Verlag, Berlin, 1987. x+592 pp.



%





\bibitem{GIOPS}
P. Pyatov, P. Saponov,
{\em Characteristic Relations for Quantum Matrices},
J. Phys. A, \textbf{28} (1995),  no. 15, 4415--4421.\\
\href{http://arxiv.org/abs/q-alg/9502012}{q-alg/9502012}\\

A. P. Isaev, O. V. Ogievetsky, P. N. Pyatov, {\it Generalized Cayley-Hamilton-Newton identities,}
 Czechoslovak J. Phys., \textbf{48}, (1998), no. 11, 1369--1374, 
\\
\href{http://arxiv.org/abs/math/9809047}{arXiv:math.QA/9809047}.

A. P. Isaev, O. V. Ogievetsky, P. N. Pyatov,
{\it On quantum matrix algebras satisfying the Cayley-Hamilton-Newton identities,} 
J. Phys. A, \textbf{32}, (1999), no. 9, L115--L121,\\
\href{http://arxiv.org/abs/math/9809170}{arXiv:math.QA/9809170}.\\

D. Gurevich, P.  Pyatov, P.  Saponov,
{\em The Cayley-Hamilton theorem for quantum matrix algebras of ${\rm GL}(m\vert n)$ type.},
St. Petersburg Math. J.,\textbf{17}  (2006),  no. 1, 119--135.\\
\href{http://arxiv.org/abs/math.QA/0508506}{math.QA/0508506}

\bibitem{GLZ}
S. Garoufalidis, T. TQ. L$\hat{\text{e}}$, D. Zeilberger,
{\em The quantum MacMahon Master Theorem}, Proc. Natl. Acad. Sci. USA 103 (2006), no. 38, 13928–13931 
\\
\href{http://arxiv.org/abs/math.QA/0303319}{math.QA/0303319}.
\bibitem{GR91} I. M. Gelfand, V. S. Retakh, {\em Determinants
of matrices over non commutative rings}, Funct. Anal. Appl., \textbf{25}
(1991), 91--102.
\bibitem{GR92}
 I. M. Gelfand, V. S. Retakh, {\it Theory of noncommutative determinants, and characteristic functions of graphs,} 
Funct. Anal. Appl., \textbf{26}, (1992), no. 4, 231--246.
\bibitem{GR97}
I. Gelfand, V. Retakh, 
{\em Quasideterminants, I}, Selecta Math. (N.S.), \textbf{3}, (1997), no. 4, 517--546,
\\
\href{http://arxiv.org/abs/q-alg/9705026}{q-alg/9705026}.
\bibitem{GGRW02}
I. Gelfand, S. Gelfand, V. Retakh, R. Wilson,
{\em Quasideterminants}, Adv. Math., \textbf{193}, (2005), no. 1, 56--141,\\
\href{http://arxiv.org/abs/math.QA/0208146}{math.QA/0208146}
\bibitem{HM06}
M. J. Hopkins, A. I. Molev,
{\em  A $q$-Analogue of the Centralizer Construction and Skew Representations of the Quantum Affine Algebra},
 SIGMA Symmetry Integrability Geom. Methods Appl. 2 (2006), Paper 092, 29 pp.
\\
\href{http://arxiv.org/abs/math.QA/0606121v2}{math.QA/0606121}
\bibitem{KL94}
D. Krob and B. Leclerc, {\em Minor identities for
quasi-determinants and quantum determinants}, Comm. Math. Phys.
\textbf{169}, 1995, pages 1--23,\\
{\href{http://arxiv.org/abs/hep-th/9411194}{hep-th/9411194}}
\bibitem{Konvalinka07-1}
M. Konvalinka,
{\em A generalization of Foata's fundamental transformation and its applications to the right-quantum algebra},
\href{http://arxiv.org/abs/math.CO/0703203}{math.CO/0703203}
\bibitem{Konvalinka07-2}
M. Konvalinka,
{\em Non-commutative Sylvester's determinantal identity},
 Electron. J. Combin.  \textbf{14}  (2007),  no. 1, Research Paper 42, 29 pp. (electronic).\\
\href{http://arxiv.org/abs/math.CO/0703213}{math.CO/0703213}
\bibitem{KoPa07}
M. Konvalinka, I. Pak, {\em Non-commutative extensions of the MacMahon
Master Theorem}. Adv. Math. 216 (2007) 29–61.
\bibitem{Lauve04}
A. Lauve, {\em Quantum- and quasi-Pl\"ucker coordinates}. Glasg. Math. J. 52 (2010), 663–675. \\
(ArXiv Math.QA 0406062).
\bibitem{LT07}
A. Lauve and E.J. Taft, {\em A class of left quantum groups modeled after SLq(r)},
Journal of Pure and Applied Algebra,
\textbf{208},  March 2007, Pages 797-803.\\
\href{http://dx.doi.org/10.1016/j.jpaa.2006.03.017}
{http://dx.doi.org/10.1016/j.jpaa.2006.03.017}
\bibitem{Maillet90}
J.-M. Maillet, {\em  Lax equations and quantum groups},
 Phys. Lett. B \textbf{245}, (1990), no. 3-4, 480--486.\\
\bibitem{Manin87}
Y. E. Manin, {\em Some remarks on Koszul algebras and quantum groups}
Ann. de l'Inst. Fourier, \textbf{37}, no. 4 (1987), p. 191-205\\
\bibitem{Manin}
Y.  Manin, {\em Quantum Groups and Non Commutative Geometry},  University de Montreal,
Centre de Recherches Mathématiques, Montreal, QC, 1988, 91 pp.
\bibitem{ManinBook91}
Y. Manin, {\em Topics in  Non Commutative Geometry}, M. B. Porter Lectures.
Princeton University Press,  1991; 164 pp.
\bibitem{Manin92}
Y.  Manin, {\em Notes on  Quantum Groups and Quantum de Rham Complexes},
Theoretical and Mathematical Physics, \textbf{92},  (1992).\\
\href{http://www.springerlink.com/content/uv1q60145574t428/}{http://www.springerlink.com/content/uv1q60145574t428/}
\bibitem{ManinDemid}
E.E. Demidov, Yu. I. Manin, E.E. Mukhin,  D.V. Zhdanovich,
{\em Non-Standard Quantum Deformations of GL(n) and Constant Solutions of the Yang-Baxter Equation},
Progress of theoretical physics. Supplement
\textbf{102} (19910329) pp. 203-218\\
\href{http://ci.nii.ac.jp/naid/110001207947/en/}{http://ci.nii.ac.jp/naid/110001207947/en/}
\bibitem{MR09}
A.I. Molev, E. Ragoucy, {\em The MacMahon Master Theorem for right quantum
superalgebras and higher Sugawara operators for $\hat{\mathfrak{g}}_{m|n}$},
ArXiv:0911.3447.
\bibitem{OPS07}
A. Oskin, S. Pakuliak, A. Silantyev,
{\em
On the universal weight function for the quantum affine algebra $U_q(\hat{\mathfrak{gl}}_N)$
},  Algebra i Analiz 21 (2009), 196--240; translation in St. Petersburg Math. J. 21 (2010), 651Ð680. 
\href{http://arxiv.org/abs/0711.2821v2}{arXiv:0711.2821}

\bibitem{ReRu10}V. S. Retakh,  V. N. Rubtsov,  {\em Noncommutative Toda Chains, Hankel Quasideterminants and Painlev\'e II Equation,}
 J. Phys. A {\bf 43} (2010),  505204, 13 pp.


\bibitem{RRT02}
S. Rodriguez-Romo, E. Taft, {\em Some quantum-like Hopf algebras which remain noncommutative when q = 1},
Lett. Math. Phys. 61 (2002), 41–50.
\href{http://www.springerlink.com/content/m6842k7r0u0q461u}
{http://www.springerlink.com/content/m6842k7r0u0q461u}

\bibitem{RRT05}
S. Rodriguez-Romo, E. Taft, {\em A left quantum group}, J. Algebra \textbf{286} (2005), 154 160.
\href{http://dx.doi.org/10.1016/j.jalgebra.2005.01.002}
{http://dx.doi.org/10.1016/j.jalgebra.2005.01.002}
\bibitem{RST}
V. Rubtsov, A. Silantyev, D. Talalaev, 
{\it Manin Matrices, Quantum Elliptic Commutative Families and Characteristic Polynomial of Elliptic Gaudin model}, 
SIGMA 5 (2009), 110, 22 pages.\\
\href{http://lanl.arxiv.org/abs/0908.4064}{arXiv:0908.4064}. 
\bibitem{TV07}
V. Tarasov, A. Varchenko,
{\em Combinatorial Formulae for Nested Bethe Vectors},
\href{http://arxiv.org/abs/math.QA/0702277v1}{math.QA/0702277}
\end{thebibliography}
\end{document}